\documentclass{amsart}
\usepackage{amsmath}
\usepackage{tikz}
\usetikzlibrary{calc}

\usepackage{
  amsmath,
  amssymb,
  amsthm,
  bbm,
  hyperref,
  paralist,
  thmtools,
  tikz,
  tikz-cd,
  pgfplots,
  subcaption,
  verbatim,
  mathrsfs
}
\pgfplotsset{compat=1.18}
\usepackage{booktabs}
\usepackage[normalem]{ulem}
\usepackage[margin=1in]{geometry}
\usetikzlibrary{angles, arrows, automata, backgrounds, calc, fit, positioning, shapes, shapes.geometric, decorations.markings, decorations.pathreplacing}
\usepackage[nocompress]{cite}
\usepackage[textsize=scriptsize]{todonotes}
\presetkeys{todonotes}{inline}{}

\makeatletter
\providecommand\@dotsep{5}
\makeatother

\renewcommand{\k}{\mathbf{k}}

\newcommand{\from}{\colon}
\newcommand{\heart}{\heartsuit}

\DeclareMathOperator{\id}{id}

\DeclareMathOperator{\dbl}{dbl}

\DeclareMathOperator{\Hom}{Hom}

\DeclareMathOperator{\std}{std}
\DeclareMathOperator{\Stab}{Stab}
\DeclareMathOperator{\supp}{supp}

\DeclareMathOperator{\PMF}{PMF}
\DeclareMathOperator{\MF}{MF}
\DeclareMathOperator{\occ}{occ}

\DeclareMathOperator{\dgmod}{dgmod-}

\declaretheorem[parent=section]{theorem}
\declaretheorem[sibling=theorem]{lemma}
\declaretheorem[sibling=theorem]{corollary}

\declaretheorem[sibling=theorem]{proposition}
\declaretheorem[sibling=theorem]{question}
\declaretheorem[sibling=theorem, style=definition]{remark}

\declaretheorem[sibling=theorem, style=definition]{definition}

\usepackage[capitalize]{cleveref}

\newcommand{\BigConf}{\operatorname{Conf}}

\newcommand{\FDomain}{\ensuremath{\overline{U}^+}}
\newcommand{\BigFDomain}{\ensuremath{U^+}}

\DeclareMathOperator{\FP}{FP}
\newcommand{\PS}{\operatorname{PS}}
\newcommand{\degen}{\operatorname{degen}}
\newcommand{\deform}{\operatorname{deform}}
\DeclareMathOperator{\TLJ}{TLJ}

\title{A sphere of spherical objects}
\author{Asilata Bapat}
\address[Asilata Bapat]{Mathematical Sciences Institute, Australian National University}
\email{asilata.bapat@anu.edu.au}
\author{Anand Deopurkar}
\address[Anand Deopurkar]{Mathematical Sciences Institute, Australian National University}
\email{anand.deopurkar@anu.edu.au}
\author{Anthony M. Licata}
\address[Anthony M. Licata]{Mathematical Sciences Institute, Australian National University}
\email{anthony.licata@anu.edu.au}

\begin{document}
\begin{abstract}
  Given a Bridgeland stability condition on a 2-Calabi--Yau category, we define a simplicial complex that encodes the Harder--Narasimhan filtrations of spherical objects.
  For 2-Calabi--Yau categories of type \(A\), we relate this complex to the complex of pointed pseudo-triangulations on configurations of points on the plane.
  Using this connection, we prove that the complex undergoes piecewise-linear wall-crossings as we vary the stability condition, and is piecewise-linearly homeomorphic to a sphere.

  Additionally, we prove that for a generic stability condition on a 2-Calabi--Yau category, a spherical object is determined by the ordered list of its Harder--Narasimhan factors.  
\end{abstract}

\maketitle
\section{Introduction}
The goal of this paper is to introduce a geometric structure on the collection of spherical objects in a 2-CY category.

Results of~\cite{bap.deo.lic:20} suggest that the collection of spherical objects should have a natural geometric structure.
Indeed,~\cite{bap.deo.lic:20} shows that the spherical objects appear in the closure of the Bridgeland stability manifold under a suitable embedding into an infinite projective space.
Our main construction in this paper associates to a Bridgeland stability condition \(\tau\) a simplicial complex \(\Sigma_\tau\), such that the spherical objects themselves are rational points on \(\Sigma_\tau\).
Although the precise combinatorial structure of the simplicial complex depends on the stability condition, we expect the complexes associated to different stability conditions to be piecewise-linearly homeomorphic to each other.
In this way, taken together, the simplicial complexes for various stability conditions should give rise to a piecewise linear structure on the set of spherical objects of \(\mathcal{C}\).

It is particularly fruitful to study the complexes \(\Sigma_{\tau}\) for the 2-CY category associated to the Coxeter system of type \(A\).
This setting connects to the geometry of curves in the plane.
We use this connection to make precise and prove the above expectations.

The simplicial complexes of interest in type \(A\) arise from two different sources:
\begin{enumerate}
\item Harder--Narasimhan filtrations of spherical objects, and
\item the geometry and combinatorics of pointed pseudo-triangulations on configurations of points in \(\mathbf{R}^2\).
\end{enumerate}
In the remainder of this introduction, we first give the main definitions and results, starting with type \(A\).
Results, questions, and conjectures that arise for other 2-CY categories are described towards the end of the introduction.
Let \(\mathcal{C}_n\) be the 2-CY category associated to the Coxeter system of type 
\(A_n\), and let \(\tau\) be a stability condition on \(\mathcal{C}_n\).
The simplicial complex \(\Sigma_{\tau}\) associated to \(\tau\) is defined as follows.
\begin{definition}\label{def:tau-simplicial-complex}
  The vertices of \(\Sigma_{\tau}\) are the \(\tau\)-semistable spherical objects of \(\mathcal{C}_n\) up to triangulated shift.
  The simplices of \(\Sigma_{\tau}\) are collections of \(\tau\)-semistable spherical objects that appear together in the Harder--Narasimhan filtration of some spherical object of \(\mathcal{C}\).
\end{definition}
The following proposition motivates the definition of \(\Sigma_{\tau}\).
In the main text, the first assertion is~\cref{cor:HN-pieces-spherical} and the second assertion is~\cref{prop:densearcs}.
\begin{proposition}\label{prop:spherical-objects-dense}
   Let \(\tau\) be a stability condition on \(\mathcal{C}_n\).
  Let \(S\) be the collection of spherical objects in \(\mathcal{C}_n\) up to triangulated shift.
  We have the following results.
  \begin{enumerate}
  \item The \(\tau\)-HN factors of a spherical object are direct sums of \(\tau\)-semistable spherical objects.
  \item Consider the map \(\supp \colon S \to \Sigma_{\tau}\) that sends a spherical object \(x\) to the formal sum of the indecomposable summands of its \(\tau\)-HN factors.
    The map \(\supp\) is injective and has dense image.
  \end{enumerate}
\end{proposition}
Thus the map \(\supp\) identifies the spherical objects in \(\mathcal{C}_n\) with a dense subset of \(\Sigma_{\tau}\).
Let \(\Stab(\mathcal{C}_n)\) be the space of stability conditions on \(\mathcal{C}_n\).
The following theorem summarises our main results about \(\Sigma_{\tau}\).
\begin{theorem}\label{thm:main}
In the setup above, we have the following results.
  \begin{enumerate}
  \item \(\Sigma_{\tau}\) is piecewise-linearly homeomorphic to the sphere of dimension \(2n-3\).
  \item\label{it:def-homeo} A continuous path from \(\tau\) to \(\tau'\) in \(\Stab(\mathcal{C}_n)\) induces a piecewise-linear homeomorphism from \(\Sigma_{\tau}\) to \(\Sigma_{\tau'}\), which depends only on the homotopy class of the path.
  \item The \((n+1)\)-strand braid group acts on \(\Sigma_{\tau}\) by piecewise-linear homeomorphisms.
    This action is compatible with its action on \(S\).
    That is, the map \(\supp \colon S \to \Sigma_{\tau}\) is equivariant with respect to the braid group.
  \end{enumerate}
\end{theorem}
The sphere \(\Sigma_{\tau}\) in~\cref{thm:main} is the ``sphere of spherical objects'' in the title.
\cref{thm:main} is a combination of~\cref{thm:pl-homeomorphisms,cor:K=sigma} in the main text.

\subsection{The simplicial complex of pointed pseudo-triangulations}
The proof of~\cref{thm:main} uses the strong relationship between \(\mathcal{C}_n\) and the geometry of the plane with \((n+1)\)-marked points, initiated in~\cite{kho.sei:02}.
We expand on this relationship to include stability conditions and HN filtrations, extending the observations in~\cite{tho:06}.
As a consequence, we rephrase the construction of \(\Sigma_{\tau}\) in purely geometric terms as the simplicial complex of pointed pseudo-triangulations.
This rephrasing is interesting in its own right and the resulting geometry is independent of the categorical considerations.

Let \(\mathbf{a}\) be a configuration of \((n+1)\) distinct points on the plane.
First let us restrict to the generic case, when no three points are collinear.
\begin{definition}
  We say that a collection of line segments with endpoints in \(\mathbf{a}\) is:
  \begin{enumerate}
  \item \emph{non-crossing} if no two segments intersect except at the endpoints, and
  \item \emph{pointed} if the segments incident at every marked point lie in a half-space.
  \end{enumerate}
  Maximal collections of segments that are non-crossing and pointed are called \emph{pointed pseudo-triangulations} or \emph{ppts}.
\end{definition}
As a special case, if the marked points form the vertices of a convex polygon, the pointedness condition is automatically satisfied at every vertex, and thus a ppt is simply a triangulation.
Ppts on generic point configurations have previously been studied in other contexts (see, e.g.~\cite{rot.san.str:08,sno.str:05,rot.san.str:03}).
We extend this definition to include point configurations that possibly have collinearities, by considering \emph{pseudo-straight} segments instead of straight segments (see~\cref{def:pseudo-straight-segment,fig:examples-pseudo-straight}, and~\cite{bap.pil:25} for the case of maximally degenerate configurations).
Further, we say that a collection of (pseudo-straight) segments is a ppt* if it is a ppt minus an external edge.
We are thus led to the following.
\begin{definition}\label{def:pptcomplex}
  Let \(\mathbf{a}\) be a configuration of distinct points in \(\mathbf{R}^2\).
  Let \(K(\mathbf{a})\) be the simplicial complex whose vertices are pseudo-straight segments joining two points of \(\mathbf{a}\).
  A collection of segments forms a simplex if it is non-crossing and pointed.
  Let \(K^{*}(\mathbf{a}) \subset K(\mathbf{a})\) be the sub-complex consisting of simplices that exclude at least one external edge.
\end{definition}
There is a dictionary between geometric features of the plane with \((n+1)\) marked points and categorical features of \(\mathcal{C}_n\).
(See also~\cref{tab:geometry-vs-category} for a summary.)
A configuration \(\mathbf{a}\) of \((n+1)\)-ordered points on the plane corresponds to a standard stability condition \(\tau\) on \(\mathcal{C}_n\).
An arc (non-crossing curve joining two marked points up to isotopy) corresponds to a spherical object on \(\mathcal{C}_n\).
Among these, the pseudo-straight (resp.~straight) segments correspond precisely to the \(\tau\)-semistable (resp.~\(\tau\)-stable) spherical objects.
An arbitrary arc can be ``pulled tight'' around the points, and it then becomes a concatenation of pseudo-straight segments.
This procedure corresponds to decomposing the associated spherical object into semistable \(\tau\)-HN factors.
\begin{table}[ht]
  \begin{tabular}[h]{l l l }
    \toprule
    Geometry & 2-CY \(A_n\) category & Reference\\
    \midrule
    configuration of \((n+1)\) points& standard stability condition \(\tau\) & \cref{prop:lifting}\\
    arc on \((n+1)\)-punctured plane& spherical object & \cref{prop:generalised-ks}\\
    straight segments & \(\tau\)-stable spherical objects & \cref{prop:curve-stables-and-semistables}\\ 
    pseudo-straight segments & \(\tau\)-semistable objects & \cref{prop:curve-stables-and-semistables}\\
    decomposition of arc by pulling tight & \(\tau\)-HN filtration & \cref{thm:hn-multiplicity-count}\\
    ppt* & \(\tau\)-HN support of a spherical object & \cref{prop:traintrack}\\
    Dehn twist & spherical twist&\cite{kho.sei:02}\\    
    \bottomrule
  \end{tabular}
  \caption{A dictionary between geometric objects and corresponding categorical objects.}\label{tab:geometry-vs-category}
\end{table}

Consequently, a simplex in \(\Sigma_{\tau}\) corresponds to a collection of pseudo-straight segments that can appear together (possibly with other pseudo-straight segments) by pulling an arc tight.
Such a collection is automatically non-crossing and pointed, and misses at least one external edge.
Thus we obtain an isomorphism between \(K^{\ast}(\mathbf{a})\) and \(\Sigma_{\tau}\).
We do not yet know a compelling categorical interpretation of \(K(\mathbf{a})\).

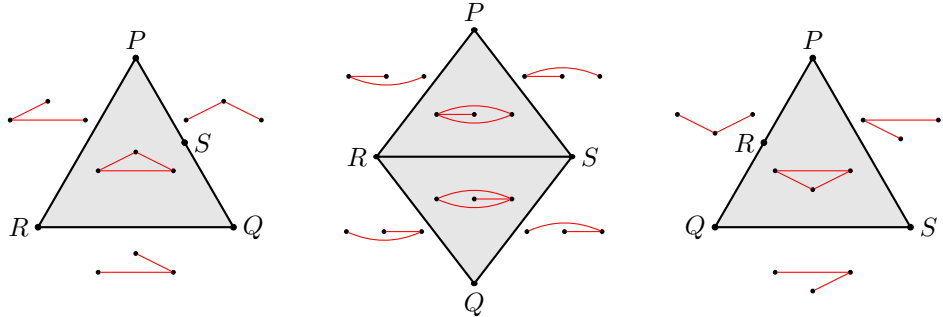
\begin{figure}
   \begin{tikzpicture}[scale=0.75]
  \tikzset{
      upper/.pic={
        \draw[fill, scale=0.5] (-1,0) node[inner sep=0em] (A) {} circle (0.05) (1,0) node[inner sep=0em] (B) {} circle (0.05) (0,0.5) node[inner sep=0em] (C) {} circle (0.05);
      }}
    \tikzset{
     lower/.pic={
        \draw[fill, scale=0.5] (-1,0) node[inner sep=0em] (A) {} circle (0.05) (1,0) node[inner sep=0em] (B) {} circle (0.05) (0,-0.5) node[inner sep=0em] (C) {} circle (0.05);
      }}
    \tikzset{
    middle/.pic={
        \draw[fill, scale=0.5] (-1,0) node[inner sep=0em] (A) {} circle (0.05) (1,0) node[inner sep=0em] (B) {} circle (0.05) (0,0) node[inner sep=0em] (C) {} circle (0.05);
      }
    }
    \draw[fill] (-30:2) circle (0.05) node[right] (Q) {\(Q\)} (90:2) circle (0.05) node[above] (P) {\(P\)} (210:2) circle (0.05) node[left] (R) {\(R\)};
    \draw[fill] (30:1) circle (0.05) node[below,right] {\(S\)};
    \draw[fill, opacity=0.1, thick] (-30:2) -- (90:2) -- (210:2) -- (-30:2);
    \draw[thick] (-30:2) -- (90:2) -- (210:2) -- (-30:2);
    \draw (30:1.8) pic {upper} (C) edge[red] (A) (C) edge [red] (B);
    \draw (-90:1.8) pic {upper} (B) edge[red] (A) (C) edge [red] (B);
    \draw (150:1.8) pic {upper} (B) edge[red] (A) (C) edge [red] (A);
    \draw (0,0) pic {upper} (B) edge[red] (A) (C) edge [red] (A) (B) edge[red] (C);
    \begin{scope}[xshift=6cm, yshift=1cm, yscale=0.75]
      \draw[fill] (-30:2) circle (0.05) node[right] (S) {\(S\)} (90:2) circle (0.05) node[above] (P) {\(P\)} (210:2) circle (0.05) node[left] (R) {\(R\)} (-90:4) circle (0.05) node[below] {\(Q\)};
     \draw[fill, opacity=0.1, thick] (-30:2) -- (90:2) -- (210:2) -- (-30:2) -- (-90:4) -- (210:2);
    \draw[thick] (-30:2) -- (90:2) -- (210:2) -- (-30:2)  -- (-90:4) -- (210:2);
    \draw (30:1.8) pic {middle} (C) edge[red] (A) (A) edge [bend left=30, red] (B);
    \draw (-60:3.2) pic {middle} (C) edge[red] (B) (A) edge [bend left=30, red] (B);
    \draw (-90:2) pic {middle} (B) edge[bend left=30, red] (A) (B) edge[bend right=30, red] (A) (C) edge [red] (B);
    \draw (150:1.8) pic {middle} (B) edge[bend left=30, red] (A) (C) edge [red] (A);
    \draw (0,0) pic {middle} (B) edge[bend left=30, red] (A) (B) edge[bend right=30, red] (A) (C) edge [red] (A);
    \draw (240:3.2) pic {middle} (B) edge[bend left=30, red] (A) (C) edge [red] (B);
  \end{scope}
  \begin{scope}[xshift=12cm]
    \draw[fill] (-30:2) circle (0.05) node[right] (S) {\(S\)} (90:2) circle (0.05) node[above] (P) {\(P\)} (210:2) circle (0.05) node[left] (Q) {\(Q\)};
    \draw[fill] (150:1) circle (0.05) node[below,left] {\(R\)};
    \draw[fill, opacity=0.1, thick] (-30:2) -- (90:2) -- (210:2) -- (-30:2);
    \draw[thick] (-30:2) -- (90:2) -- (210:2) -- (-30:2);
    \draw (30:1.8) pic {lower} (C) edge[red] (A) (A) edge [red] (B);
    \draw (-90:1.8) pic {lower} (B) edge[red] (A) (C) edge [red] (B);
    \draw (150:2) pic {lower} (C) edge[red] (A) (C) edge [red] (B);
    \draw (0,0) pic {lower} (B) edge[red] (A) (C) edge [red] (A) (B) edge[red] (C);
  \end{scope}
\end{tikzpicture}
\caption{The piecewise-linear wall-crossing transformations induced on the complex \(K\).
  On the left is the complex \(K\) for three non-collinear points; in the middle for three collinear points obtained by specialising the middle point onto the line segment joining the other two; on the right for the three non-collinear points as the middle point deforms downwards.
}\label{fig:wall-crossing}
\end{figure}

We prove statements about \(\Sigma_{\tau}\) by proving the analogous geometric statements about \(K^{\ast}(\mathbf{a})\).
Specifically,~\cref{prop:spherical-objects-dense} follows by observing that arcs pull tight to chains of arcs and from~\cref{thm:hn-multiplicity-count} as well as~\cref{prop:densearcs}.
\cref{thm:main} is a consequence of~\cref{thm:pl-homeomorphisms,cor:braid-pl-action,cor:ball}.
These statements about \(K^{\ast}(\mathbf{a})\) and \(K(\mathbf{a})\) are proved by studying the wall-crossing between the complexes as the point configurations change.
\Cref{fig:wall-crossing} shows the simplest example of a wall-crossing for \(K(\mathbf{a})\) as \(\mathbf{a}\) degenerates from a triangle to three collinear points, and then deforms again to a triangle.
The subcomplex \(K^{*}(\mathbf{a})\), which in this case is the \(1\)-sphere of spherical objects, is the boundary circle.

\begin{remark}\label{rem:pmf}
  The geometric realisation \(|K^{\ast}(\mathbf{a})|\) of \(K^{*}(\mathbf{a})\) is naturally homeomorphic to the space of projective measured foliations (PMFs) on the \((n+1)\)-punctured disk (see~\cref{sec:pmf}).
  The piecewise-linear wall-crossings and braid group action on \(|K^{*}(\mathbf{a})|\) described in \cref{thm:main} are equivalent to the corresponding properties of the space of PMFs.
  However, our proof of \cref{thm:main} is independent of this connection.
\end{remark}

\begin{remark}
  When the configuration of \((n+1)\) points forms the vertices of a convex polygon, the complex of ppts \(K(\mathbf{a})\) is just the complex of triangulations of the polygon.
  This complex \(K(\mathbf{a})\) is an iterated cone over the cluster complex in type \(A\) (see, e.g.~\cite{fom.zel:03}), which is the polar dual of the associahedron.
  It would be interesting to know if complexes \(K(\mathbf{a})\) and \(K^{*}(\mathbf{a})\) for non-convex configurations also appear independently in cluster theory or other parts of mathematics which study the combinatorics of triangulations.
\end{remark}

\begin{remark}\label{rem:ahat}
  For the 2-CY category associated to the Coxeter system of type \(\widehat A_n\), there is a dictionary, similar to~\cref{tab:geometry-vs-category}, between spherical objects/stability conditions and the topology/geometry of arcs on an on an annulus with 
\((n+1)\)-marked points ~\cite{gad.thi.wag:17}.
  As a result, we expect analogues of~\cref{prop:spherical-objects-dense,thm:main} to follow from methods similar to this paper.
  Note, however, that the complex \(\Sigma_{\tau}\) is more complicated to study (for example, it has infinitely many simplices).
  In particular, results of~\cite{bap.deo.lic:20} show that when \(n=1\), the complex \(\Sigma_\tau\) is not a circle, but is instead homeomorphic to \(\mathbf{R}\).
\end{remark}

\subsection{The complex of sphericals in other 2-CY categories}
Let \(\mathcal{C}\) be any 2-CY category, linear over an infinite field \(\k\), of finite type, and which admits a dg enhancement.  A version of~\cref{prop:spherical-objects-dense} continues to hold.  Namely, suppose \(\tau\) is a stability condition that is generic in the following sense: all \(\tau\)-semi-stable spherical objects are \(\tau\)-stable.
Then the \(\tau\)-HN factors of a spherical object are spherical (this follows from the Mukai lemma; see~\cref{sec:graded-spherical-support}).
As a result, we can define the \(\tau\)-HN support of a spherical object as before.
Thus, for a generic stability condition \(\tau\),~\cref{def:tau-simplicial-complex} gives a simplicial complex whose points include HN supports of spherical objects.
We do not know whether the HN support always characterises spherical objects, but we prove a close result: for spherical objects in \(\mathcal{C}\), the \emph{ordered} list of HN factors determines the spherical object (\cref{thm:object-determined-by-hn-filtration}).

\subsection{The sphere of sphericals for finite-type Coxeter systems}
In addition to the 2-CY categories of type \(A\), we also completely understand the complex of sphericals for 2-CY categories associated to rank two Coxeter systems.
Let \(\mathcal{C}\) be a 2-CY category associated to a rank two Coxeter system.
Then \(\mathcal{C}\) carries an additional action of a fusion category, and it is natural to study stability conditions that are equivariant with respect to this action.
Likewise, it is natural to consider the collection of spherical objects in \(\mathcal{C}\) up to triangulated shift and the action of the fusion category.

When \(\mathcal{C}\) is the 2-CY category of a rank two Coxeter system, fusion-equivariant stability conditions on \(\mathcal{C}\) are well-understood (see~\cite{hen.lic:24,hen:24,del.hen.lic:23,hen:22}) via the machinery of Harder--Narasimhan automata.  We use these results to establish analogues of~\cref{prop:spherical-objects-dense,thm:main} in this setting in~\cref{sec:rank-2}.  In particular, for any finite Coxeter group \(W\) of rank two, the simplicial complexes associated to any stability condition are all homeomorphic to \(S^1\).

\subsection{Further questions and conjectures}\label{sec:questions-and-conjectures}
The paper raises several questions---on the categorical side, on the geometric side, and on the interface.
We gather a few with some remarks.

\cref{def:tau-simplicial-complex} applies to generic stability conditions on arbitrary 2-CY categories, not just to the ones covered by our main theorems.
This includes 2-CY categories for more general graphs or Coxeter systems, as well as 2-CY categories arising in geometry such as the derived categories of coherent sheaves on K3 surfaces.
Let \(\mathcal{C}\) be an arbitrary 2-CY category.
\begin{question}
  Let \(\tau\) be a generic stability condition on \(\mathcal{C}\).
  What can we say about the complex \(\Sigma_{\tau}\)?
  Does it remain a manifold?
  If it does, of what dimension?
\end{question}
The key to answering these questions is understanding the constraints on the stable factors of a spherical object in a 2-CY category.
In type \(A\), we characterise these constraints geometrically in terms of ppt*s.
We conjecture that for generic stability conditions on 2-CY categories- associated to finite-type Coxeter systems of rank \(r\), the simplicial complex \(\Sigma_\tau\) is homeomorphic to the sphere of dimension \((2r-3)\).
\begin{question}
  Let \(\tau\) be a generic stability condition on \(\mathcal{C}\).
  What are the constraints on the collections of stable spherical objects that can appear together in the \(\tau\)-HN filtration of a spherical object?
\end{question}

Now let \(\tau\) be any stability condition on \(\mathcal{C}\), not necessarily generic.
In general, the \(\tau\)-HN factors of spherical objects are not necessarily direct sums of spherical objects (a counterexample can be found in type \(D_4\)).
So we cannot define the support of a spherical object as a point of \(\Sigma_{\tau}\).
In this setting, we should replace \(\Sigma_{\tau}\) by a different simplicial complex.
\begin{question}
  What is the correct definition of \(\Sigma_{\tau}\) for non-generic stability conditions \(\tau\)?
\end{question}

Let \(\mathcal{C}\) be the 2-CY category associated to a Coxeter system.
Let \(B\) be the Artin-Tits group of the Coxeter system.
Set \(\Sigma = \Sigma_{\tau}\) for some generic stability condition \(\tau\).
The group \(B\) acts on \(\Sigma\) by homeomorphisms.  
As we establish in type \(A\) in our main theorem (\cref{thm:main}), we expect this action to be  piecewise-linear.  We end with a somewhat open ended question.
\begin{question}
  To what extent can \(B\) be recovered from the piecewise-linear structure on \(\Sigma\)?
\end{question}
Our motivation comes from analogous rigidity statements in Teichm\"uller theory.
The space of projective measured foliations on a surface admits a piecewise-linear (PL) structure using the train-track atlas.
One can define a distinguished collection of functions with respect to this atlas.
Under suitable hypotheses, the mapping class group consists precisely of PL automorphisms that preserve this collection (see, e.g.~\cite[Theorem~3.11]{pap:15}).
In this vein, consider functions \(f \colon \Sigma \to \mathbf{R}\) that are linear combinations of Harder--Narasimhan mass functions of a finite number of stability conditions.  The braid group preserves this collection of functions.  Is the braid group precisely the group of PL automorphisms of \(\Sigma\) that preserves this collection?     

\subsection{Organisation}
We have organised the paper so that the geometry precedes the category theory, and is self-contained. In \cref{sec:basic-geometry}, we define the main geometric players.  In~\cref{sec:simplicial-complex}, we define and study the simplicial complex of pointed pseudo-triangulations.
In~\cref{sec:pmf}, we relate this complex to with the space of projective measured foliations.

In~\cref{sec:categorical-setup}, we recall the definitions of the categories, and give a brief reminder on Bridgeland stability conditions.  We also recall results of Khovanov--Seidel~\cite{kho.sei:02} that relate these categories in type \(A\) to the geometric setup discussed earlier.
In~\cref{sec:stab-via-point-configurations} we give a geometric description of semistable objects and Harder--Narasimhan filtrations.

In~\cref{sec:rank-2}, we describe the simplicial complex \(\Sigma_{\tau}\) for (fusion-equivariant) stability conditions for 2-Calabi--Yau categories of rank two Coxeter systems.

In~\cref{sec:graded-spherical-support}, we return to the setting of stability conditions on a general 2-Calabi--Yau category.
  We prove that in a generic stability condition, a spherical object is determined by the ordered list of its Harder--Narasimhan factors.

In Appendix~\ref{appendixA}, we give a self-contained account of a spectral sequence that computes morphisms between two objects in terms of morphisms between their filtered pieces.
  
\subsection*{Acknowledgements}
We thank Edmund Heng, Ian Le, Vincent Pilaud, Nicholas Proudfoot, and Hoel Queffelec for several helpful conversations.
We are supported by the Australian Research Council Award DP240101084, and A. B. is additionally supported by the Australian Research Council award DE240100447.
We acknowledge support from the International Research Lab France-Australia Mathematical Sciences and Interactions (IRL-FAMSI).

\section{Geometric setup and background}\label{sec:basic-geometry}
In this section we introduce the geometric setup we need, including arcs on point configurations (\cref{sec:configurations}) as well as the notion of pointed-pseudo-triangulations (\cref{sec:ppt}).
\subsection{Configurations, arcs, segments, and spines}\label{sec:configurations}
Fix a non-negative integer \(n\).
A \emph{point configuration} is a sequence of \(n+1\) distinct points in \(\mathbf{C}\).
Let  \(\mathbf{a} = (a_0, \dots, a_n)\) be a point configuration and let \(i,j \in \{0,\dots,n\}\) be distinct.
\begin{definition}
  An \emph{oriented arc} from \(a_i\) to \(a_j\) is an isotopy class of a continuous map \(\gamma \colon [0,1] \to \mathbf{C}\) with \(\gamma(0) = a_i\) and \(\gamma(1) = a_j\) with \(\gamma((0,1)) \subset \mathbf{C} \setminus \{a_0, \dots, a_n\}\).
  
  Fix an orientation on \(S^1\).
  An \emph{oriented curve} on \(\mathbf{a}\) is the isotopy class of a simple closed curve \(\gamma \colon S^1 \to \mathbf{C} \setminus \{a_0, \dots, a_n\}\) that encloses at least 2 marked points.
\end{definition}
\begin{definition}
  An \emph{arc} (resp.~\emph{curve}) is an oriented arc (resp.~curve) modulo the orientation.
  The \emph{length} \(\ell(\gamma)\) of an arc (resp.~curve) \(\gamma\) is the infimum of the lengths of the representatives in its isotopy class.
\end{definition}
The \emph{boundary-parallel curve} on a configuration \(\mathbf{a}\) is the curve represented by a circle that contains all points of \(\mathbf{a}\) in its interior.
\begin{definition}
  An \emph{(oriented) multi-curve} is a multi set \(\{c_1, \dots, c_k\}\) whose elements \(c_{i}\) are either (oriented) curves or arcs, and such that we can choose representatives \(\gamma_1, \dots, \gamma_k\) of \(c_1, \dots, c_k\) that are pairwise disjoint.
  Such a multi-curve is \emph{admissible} if no component \(c_i\) is the boundary parallel curve.
\end{definition}
Let \(\BigConf_{n+1}\) be the space of point configurations modulo translations.
We can identify \(\BigConf_{n+1}\) with the complement of the (complexified) hyperplane arrangement of type \(A_n\).
Indeed, consider the quotient of \(\mathbf{C}^{n+1}\) by the diagonally embedded copy of \(\mathbf{C}\).
Then \(\BigConf_{n+1}\) is the complement of the hyperplanes \(x_i = x_j\) in this quotient.

The symmetric group \(S_{n+1}\) acts freely on \(\BigConf_{n+1}\) by permuting the points.
The fundamental group of \(\BigConf_{n+1}/S_{n+1}\) is the \((n+1)\)-strand braid group \(B_{n+1}\).
The braid group acts on the set of (oriented) arcs.

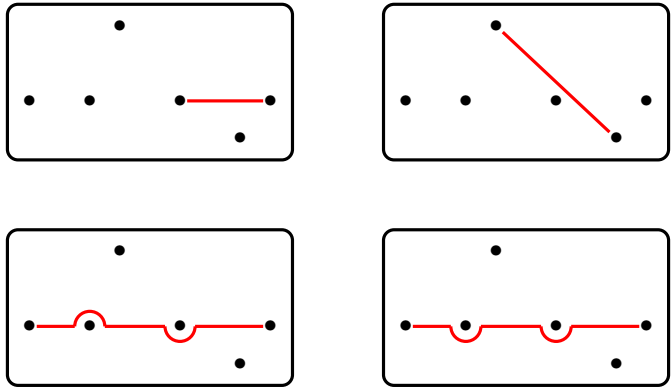
\begin{figure}[ht]
  \centering
  \begin{tikzpicture}
    \tikzset{
      config/.pic={
        \node[inner sep=0pt] (-pt1) at (-0.5,0) {\(\bullet\)};
        \node[inner sep=0pt] (-pt2) at (0.3,0) {\(\bullet\)};
        \node[inner sep=0pt] (-pt3) at (1.5,0) {\(\bullet\)};
        \node[inner sep=0pt] (-pt4) at (2.7,0) {\(\bullet\)};
        \node[inner sep=0pt] (-pt5) at (0.7,1) {\(\bullet\)};
        \node[inner sep=0pt] (-pt6) at (2.3,-0.5) {\(\bullet\)};
        \coordinate (-pt2left) at (0.1,0.0);
        \coordinate (-pt2right) at (0.5,0.0);
        \coordinate (-pt3left) at (1.3,0.0);
        \coordinate (-pt3right) at (1.7,0.0);        
        \node[draw, solid, very thick, rounded corners, inner sep=2mm, fit=(-pt1) (-pt5) (-pt4) (-pt6)] {};
      }}
    \node[matrix] at (0,0) {\pic (c1) {config};\\};
    \node[matrix] at (5,0) {\pic (c2) {config};\\};
    \node[matrix] at (0,-3) {\pic (c3) {config};\\};
    \node[matrix] at (5,-3) {\pic (c4) {config};\\};
    \path[very thick, red]
    (c1-pt3) edge (c1-pt4)
    ;

    \path[very thick, red]
    (c2-pt5) edge (c2-pt6)
    ;
    
    \draw[very thick, red]
    (c3-pt1) edge (c3-pt2left)
    (c3-pt2left) arc (180:0:0.2)
    (c3-pt2right) edge (c3-pt3left)
    (c3-pt3left) arc (180:360:0.2) 
    (c3-pt3right) edge (c3-pt4)        
    ;

     \draw[very thick, red]
    (c4-pt1) edge (c4-pt2left)
    (c4-pt2left) arc (180:360:0.2)
    (c4-pt2right) edge (c4-pt3left)
    (c4-pt3left) arc (180:360:0.2) 
    (c4-pt3right) edge (c4-pt4)        
    ;

  \end{tikzpicture}
  \caption{Several examples of pseudo-straight segments}
  \label{fig:examples-pseudo-straight}
\end{figure}

\begin{figure}[ht]
    \centering
    \begin{tikzpicture}
      \tikzset{
        config/.pic={
          \node[inner sep=0pt] (-pt1) at (0,1) {\(\bullet\)};
          \node[inner sep=0pt] (-pt2) at (-1,0) {\(\bullet\)};
          \node[inner sep=0pt] (-pt3) at (2,0) {\(\bullet\)};
          \node[inner sep=0pt] (-pt4) at (1.1,-0.8) {\(\bullet\)};
          \node[inner sep=0pt] (-pt5) at (0.7,-1.4) {\(\bullet\)};
          \node[inner sep=0pt] (-pt6) at (0.2,0) {\(\bullet\)};          
          \coordinate (-pt2bot) at (-1,-0.2);
          \coordinate (-pt2left) at (-1.2,0);
          \coordinate (-pt1top) at (0, 1.2);
          \coordinate (-pt3right) at (2.2, 0);
          \coordinate (-pt4left) at (0.9, -0.8);
          \coordinate (-pt5bot) at (0.7, -1.6);
          \coordinate (-pt6top) at (0.2,0.2);
          \coordinate (-pt6left) at (0.0,0.0);
          \coordinate (-pt6right) at (0.4,0.0);
          \node[draw, solid, very thick, rounded corners, inner sep=2mm, fit=(-pt1top) (-pt2left) (-pt3right) (-pt5bot)] {};          
        }}
      \node[matrix] (f1) at (0,0) {\pic (c1) {config};\\};
      \node[matrix] (f2) at (6,0) {\pic (c2) {config};\\};

      \draw[very thin, ->] (f1) -- node[above]{spine} (f2);
      
      \path[thick,red,decoration={markings,mark=at position 0.5 with {\arrow{>}}}]
      (c1-pt3) edge[out=180, in=0,postaction=decorate] (c1-pt6top)
      (c1-pt6top) edge[out=180, in=0] (c1-pt2bot)
      (c1-pt2bot) edge[out=180, in=-90] (c1-pt2left)
      (c1-pt2left) edge[out=90,in=180] (c1-pt1top)
      (c1-pt1top) edge[out=0,in=90,postaction=decorate] (c1-pt3right)
      (c1-pt3right) edge[out=-90, in=40] (c1-pt4left)
      (c1-pt4left) edge[out=-140, in=180] (c1-pt5bot)
      (c1-pt5bot) edge[out=0, in=-90, postaction=decorate] (c1-pt4);
      
      \draw[thick,blue,decoration={markings,mark=at position 0.4 with {\arrow{>}}}, font=\scriptsize]
      (c2-pt3) edge[postaction=decorate] (c2-pt6right)
      (c2-pt6right) arc (0:180:0.2)
      (c2-pt6left) edge[postaction=decorate] node[below]{\(\gamma_1\)} (c2-pt2)
      (c2-pt2) edge[postaction=decorate] node[above left]{\(\gamma_2\)} (c2-pt1)
      (c2-pt1) edge[postaction=decorate] node[above right]{\(\gamma_3\)} (c2-pt3)
      (c2-pt3) edge[postaction=decorate] node[below right]{\(\gamma_4\)}(c2-pt4)
      (c2-pt4) edge[bend right=30, postaction=decorate] node[above left]{\(\gamma_5\)}(c2-pt5)
      (c2-pt5) edge[bend right=30, postaction=decorate] node[below right]{\(\gamma_6\)} (c2-pt4);      
    \end{tikzpicture}
    \caption{The spine of the arc on the left (red) is the chain \((\gamma_1, \dots, \gamma_6)\).}
    \label{fig:spine}
  \end{figure}
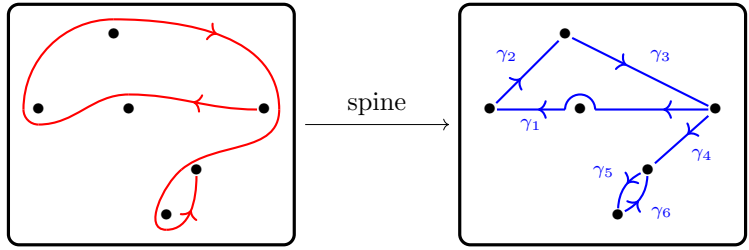

Fix a point configuration \(\mathbf{a} = (a_0,\dots, a_n)\).
Given a curve \(\gamma\), we decompose it as a sequence of ``pseudo-straight segments'' that we call its ``spine.''
We first describe this process informally, and then give the formal definitions.
Let \(\gamma\) have length \(\ell\).
There may not be a curve isotopic to \(\gamma\) that has length \(\ell\).
Nevertheless, consider curves whose lengths approach \(\ell\).
These curves approach a sequence of straight line segments that bend at the marked points.
This is the shape that the curve would take if it were made of a physical string that was pulled tight while constrained to avoid the marked points (see~\cref{fig:spine}).
The ``spine'' of \(\gamma\), roughly speaking, is the sequence of these straight line segments.
If the point configuration contains triples of collinear points, we refine the notion of a spine by taking it to be a sequence of ``pseudo-straight segments'' instead of straight line segments.
A pseudo-straight segment from \(a_i\) to \(a_j\), roughly speaking, is a curve that follows the straight line segment from \(a_i\) to \(a_j\) except possibly swerving to avoid points of the configuration that are in the way (see~\cref{fig:examples-pseudo-straight}).

\begin{definition}\label{def:pseudo-straight-segment}
  Given a configuration \(\mathbf{a} = (a_0, \dots, a_n)\) and \(i \neq j\), a \emph{pseudo-straight} segment on \(\mathbf{a}\) from \(a_i\) to \(a_j\) is an arc from \(a_i\) to \(a_j\) whose length is \(|a_i-a_j|\), the euclidean distance between \(a_i\) and \(a_j\).

  We denote by \(\PS(\mathbf{a})\) the set of (unoriented) pseudo-straight segments on \(\mathbf{a}\).
\end{definition}

Let \(s \colon [0,1] \to \mathbf{C}\) be the standard parametrisation of the straight line segment from \(a_i\) to \(a_j\).
Given any \(\epsilon > 0\), a pseudo-straight segment from \(a_i\) to \(a_j\) has a representative \(\gamma\) that is \(\epsilon\)-close to \(s\); that is, satisfies \(|\gamma(t) - s(t)| < \epsilon\) for all \(t \in [0,1]\).
Therefore, the image of \(\gamma\) lies in an \(\epsilon\) neighbourhood of the segment from \(a_i\) to \(a_j\).
If no other point \(a_k\) lies on the segment joining \(a_i\) and \(a_j\), then \(s\) is the unique pseudo-straight segment from \(a_i\) to \(a_j\), up to isotopy.
Otherwise, there are several (see~\cref{fig:examples-pseudo-straight}).
To represent such a pseudo-straight segment, we may choose the curve that follows \(s\) except in a small neighbourhood of an intermediate point \(a_k\), where it goes in a small semicircle around \(a_k\), either on the left or the right, and then rejoins \(s\).
If there are \(r\) points of \(\mathbf{a}\) on the open segment from \(a_i\) to \(a_{j}\), then there are \(2^r\) pseudo-straight segments from \(a_i\) to \(a_j\).
As a result, \(\PS(\mathbf{a})\) is a finite set.

\begin{remark}
  The paper~\cite{bap.pil:25} studies the combinatorics of pseudo-straight segments and their pointed pseudotriangulations (which we define in the next section), in the case of a configuration that consists entirely of collinear points.
  In \cite{bap.pil:25}, pseudo-straight segments are called wiggly arcs.
\end{remark}

\begin{definition}[Chains and collinear chains]
  \label{def:chain}
  By a \emph{chain of arcs}, we mean a finite sequence of oriented arcs \((\gamma_1, \dots, \gamma_k)\) such that the end point of \(\gamma_i\) is the starting point of \(\gamma_{i+1}\).
  Let the start and end points of \(\gamma_i\) be \(a_{p_{i-1}}\) and \(a_{p_{i}}\) respectively.
  We say that the chain is \emph{collinear} if all \(\gamma_i\) are pseudo-straight segments and the points \(a_{p_0}, \dots, a_{p_k}\) lie on a line in order, that is, \(a_{p_i}\) lies between \(a_{p_{i-1}}\) and \(a_{p_{i+1}}\).
\end{definition}

The spine of a curve is the chain obtained by ``pulling the curve tight'' around the marked points.
A more precise definition follows.
  \begin{definition}[Spine]
  \label{def:spine}
  Let \(\mathbf{a} = (a_0, \dots, a_n)\) be a configuration and let \(\gamma\) be an oriented arc on \(\mathbf{a}\).
  The \emph{spine} of \(\gamma\), denoted by \(\operatorname{spine}(\gamma)\), is the chain of pseudo-straight segments \((\gamma_1, \dots, \gamma_k)\) satisfying the following properties.
  \begin{enumerate}
  \item \(\ell(\gamma) = \sum_i \ell(\gamma_i)\).
  \item For any \(i\), the chain \((\gamma_i, \gamma_{i+1})\) is not collinear.
  \item For any representatives \(s_i\) of \(\gamma_i\) and any \(\epsilon > 0\), the arc \(\gamma\) has a representative that is \(\epsilon\)-close to the concatenation of \(s_1, \dots, s_k\).
  \end{enumerate}
\end{definition}

  See~\cref{fig:spine} for an illustration of the spine of an arc.
  The spine of an unoriented arc is the spine of one of its orientations, up to reversal.
Implicit in~\cref{def:spine} is the existence and uniqueness of spines, whose proof we omit.

The definition of a \emph{spine} has an obvious analogue for curves and multi-curves.
Observe that for an curve, the spine is unique only up to a cyclic rotation.
We say that a pseudo-straight segment is an \emph{external edge} of \(\mathbf{a}\) if it occurs as a component of the spine of the boundary-parallel curve.

\subsection{Pointed pseudo-triangulations}\label{sec:ppt}
Fix a point configuration \(\mathbf{a}\).
The set of pseudo straight segments that appear in the spine of a multi-curve satisfies some constraints.
Identifying these constraints leads to the notion of a pointed pseudo-triangulation (ppt).
The goal of this section is to define ppts and construct a simplicial complex associated to every ppt that acts as a parameter space for curves.

The notion of a ppt has appeared before in a different context---rigidity theory and motion planning in robotics \cite{rot.san.str:08}.
The simplest situation is when the points of \(\mathbf{a}\) form the vertices of a convex polygon.
Then, a ppt is the same as a triangulation of the polygon.
A more general situation is when \(\mathbf{a}\) forms a non-degenerate (no three collinear points) but not necessarily convex configuration.
In this case, the notion of ppts is developed in \cite{rot.san.str:08}.
We extend the theory to incorporate collinear points.
We note that the most extreme case of collinearity, namely when all points of \(\mathbf{a}\) lie on a single line, has also been studied previously \cite{bap.pil:25}.

 Fix a configuration \(\mathbf{a} = (a_0, \ldots, a_n)\) of distinct points in \(\mathbf{C}\).
 Recall that there are finitely many pseudo-straight curves (as in~\cref{def:pseudo-straight-segment}) with endpoints in \(\mathbf{a}\).
 \begin{definition}[non-crossing, pointed]
   Let $T$ be a set of pseudo-straight segments with endpoints in \(\mathbf{a}\).
   \begin{enumerate}
   \item We say that \(T\) is \emph{non-crossing} if there exists a set of representatives for elements of \(T\) that are pairwise non-intersecting except at the endpoints.
   \item We say that \(T\) is \emph{pointed} at \(a_i \in \mathbf{a}\) if all the straight line segments underlying the pseudo-straight segments incident to \(a_i\) lie in an open half plane at \(a_i\).
     We say that \(T\) is \emph{pointed} if it is pointed at each \(a_i \in \mathbf{a}\).
   \end{enumerate}
\end{definition}
 \begin{definition}[Pointed pseudo-triangulations---ppt and ppt*]
   \label{def:ppt}
   Let \(\mathbf{a}\) be a point configuration in \(\mathbf{C}\).

   A \emph{pointed pseudo-triangulation} or \emph{ppt} on \(\mathbf{a}\) is a maximal collection of pseudo-straight segments on $\mathbf{a}$ that is non-crossing and pointed.

   A \emph{ppt*} is a maximal collection of pseudo-straight segments on \(\mathbf{a}\) that is non-crossing and pointed and does not contain all external edges.
 \end{definition}

 Owing to the maximality of a ppt, every ppt on \(\mathbf{a}\) must contain all the external edges of \(\mathbf{a}\).
 So a ppt* is simply a ppt minus one external edge.

We prove the following lemma for later use.
\begin{lemma}\label{lem:edges-upper-bound}
  Let \(n \geq 1\) and fix a point configuration \(\mathbf{a}\) of \(n+1\) points, and let \(T\) be a set of pseudo-straight segments on \(\mathbf{a}\) that is planar and pointed.
  Let \(e = |T|\).
  Then \(e \leq 2n - 1\).
  In particular, a ppt on \(n+1\) points has at most \(2n-1\) edges.
\end{lemma}
\begin{remark}
  In the case of non-degenerate point configurations (no three points collinear) as well as maximally degenerate configurations (all points collinear) of \(n+1\) points, it is known that every ppt has exactly \(2n-1\) edges (see, e.g.~\cite{rot.san.str:08,bap.pil:25}).
  We will prove in~\cref{prop:ppt-cardinality} that the same is true for all point configurations.
  However, an upper bound is sufficient for the moment.
\end{remark}
\begin{proof}
  Let \(G\) be the undirected graph with vertex set \(V(G) = \{a_0, \dots, a_n\}\) and edge set \(E(G) = T\).
  We use \(v(G)\) and \(e(G)\) to denote the sizes of \(V(G)\) and \(E(G)\), respectively, dropping \(G\) if it is clear from the context.
  We want to prove that
  \[ e(G) \leq 2v(G) -3.\]

  Fix a set of representatives of elements of \(T\) that are non-crossing and pointed.
  Then we have a embedding of \(G\) in \(\mathbf{C}\).
    The complement of \(G\) in \(\mathbf{C}\) is a disjoint union of regions, which we call \emph{faces}.

  Assume that \(G\) has no isolated points.
  We induct on the number of connected components of \(G\).
  As the base case, assume that \(G\) is connected.

  Exactly one face is unbounded, and suppose that there are \(f\) bounded faces.
    We use a variant of the Euler characteristic argument used, e.g., in~\cite[Theorem 2.4]{rot.san.str:08}.
  Euler's formula gives the equation
  \begin{equation}\label{eq:euler}
    v - e + f = 1.
  \end{equation}
  Since \(G\) is connected, every vertex is incident to at least one edge.
  Then the number of angles incident to each vertex is exactly the valence of the vertex.
  Each angle is either convex (at most \(\pi\)) or reflex (greater than \(\pi\)).
  Let us now count the angles in the graph in two different ways.

  Summing over the vertices, we see that the total number of angles is the sum of the valences of the vertices.
  On the other hand, the sum of the valences of the vertices is twice the number of edges.
  Thus
  \begin{equation}\label{eq:angles}
    \text{total number of angles} = 2e.
  \end{equation}
  The pointedness condition implies that exactly one angle incident to each point is a reflex angle, and that all other angles are strictly convex.
  Thus
  \begin{equation}\label{eq:reflex-angles}
    \text{number of reflex angles} = v.
  \end{equation}
  The total number of convex angles is equal to the sum over the convex angles of every bounded face.
  Since \(G\) is connected, a bounded face must be simply connected and therefore must have at least three edges.
  Indeed, if it had only two edges, they would be isotopic to each other.
  Thus
  \begin{equation}\label{eq:convex-angles}
    \text{number of convex angles} \geq 3f.
  \end{equation}
  Now we combine~\cref{eq:angles,eq:reflex-angles,eq:convex-angles} with Euler's formula~\eqref{eq:euler} and simplify to see that
  \[e(G) \leq 2v(G) - 3,\]
  as desired.

  Suppose \(G\) has more than one connected component.
  Let \(H \subset G\) be a connected component that is innermost in the embedding of \(G\) in \(\mathbf{C}\); that is, no face of \(H \subset \mathbf{C}\) contains another component of \(G\).
  Let \(\mathbf{a}' = \mathbf{a} - V(H)\) and \(T' = T - E(H)\).
  Then \(T'\) consists of pointed and non-crossing pseudo-straight segments on \(\mathbf{a}'\).
  Moreover, the elements of \(T'\) represent distinct pseudo-straight segments on \(\mathbf{a}'\), except for the following possibility.
  If a bounded face of \(G-H \subset \mathbf{C}\) contains the image of \(H\) and only has two edges, the two edges represent the same pseudo-straight segment on \(\mathbf{a}'\).
  In any case, we have
  \[ e(G) \leq e(G-H) + e(H) + 1.\]
  We wish to apply the inductive hypothesis to \(H\) and \(G-H\), but we must account for the possibility that \(v(H)\) or \(v(G-H)\) is \(1\).
  If both \(v(H)\) and \(v(G-H)\) are at least \(2\), we apply the inductive hypothesis to both and get
  \[ e(G) \leq 2v(G-H) - 3 + 2v(H) - 3 + 1 < 2v(G) - 3.\]
  If exactly one of \(v(H)\) or \(v(G-H)\) is \(1\), say the first, we apply the inductive hypothesis to \(G-H\) and get
  \[ e(G) \leq 2(v (G)-1)-3 + 1 < 2 v(G) - 3.\]
  If both \(v(H)\) and \(v(G-H)\) are \(1\), then \(0 = e(G) < 2v(G) - 3.\)
\end{proof}

\section{The simplicial complex of pointed pseudo-triangulations}\label{sec:simplicial-complex}
The aim of this section is to study the simplicial complex consisting of ppt*s or ppts on any point configuration.
We prove in~\cref{thm:pl-homeomorphisms} that there are piecewise-linear homeomorphisms between the complexes associated to any two point configurations.
Thus the braid group acts on these complexes in a piecewise-linear way (\cref{cor:braid-pl-action}).
Moreover, we show in~\cref{cor:ball} that the complex of ppt*s (resp.~ppts) is a sphere (resp.~a closed ball).

\subsection{Construction and basic properties}
We refer the reader to~\cref{sec:basic-geometry} for the definitions of the objects used in this section.
\begin{definition}[The complexes \(K^{*}(\mathbf{a})\) and \(K(\mathbf{a})\)]
  \label{def:pptconf}
  Let \(\mathbf{a}\) be a point configuration in \(\mathbf{C}\).
  The simplicial complex \(K^{*}(\mathbf{a})\) is the abstract simplicial complex whose maximal simplices are the ppt*'s on \(\mathbf{a}\).
  The simplicial complex \(K(\mathbf{a})\) is the abstract simplicial complex whose maximal simplices are the ppts on \(\mathbf{a}\).  
\end{definition}
Recall that \(\PS(\mathbf{a})\) is the set of pseudo-straight segments on \(\mathbf{a}\).
Let \(\Delta(\PS(\mathbf{a})) = (\mathbf{R}_{\geq 0}^{\PS(\mathbf{a})} - \{0\})/ \mathrm{scaling}\) be the standard simplex on \(\PS(\mathbf{a})\).
We think of points of \(\Delta(\PS(\mathbf{a}))\) as non-negative real-valued functions on \(\PS(\mathbf{a})\) modulo simultaneous scaling.
The geometric realisation of \(K^{\ast}(\mathbf{a})\) (resp.~\(K(\mathbf{a})\)), denoted \(|K^{\ast}(\mathbf{a})|\) (resp.~\(|K(\mathbf{a})|\)), is then the subspace of \(\Delta(\PS(\mathbf{a}))\) consisting of the functions whose support is contained in a ppt* (resp.~ppt).
\begin{remark}\label{rem:unproj}
  We occasionally also need the unprojectivised spaces \(\widetilde K^{*}(\mathbf{a})\) and \(\widetilde{K}(\mathbf{a})\), defined simply as the pre-images of \(|K^{*}(\mathbf{a})|\) and \(|K(\mathbf{a})|\) in \(\mathbf{R}_{\geq 0}^{\PS(\mathbf{a})} - \{0\}\).
\end{remark}

\begin{definition}[Support of a multi-curve]
  Let \(\gamma\) be a curve or an arc.
  The \emph{support} of \(\gamma\), denoted by \(\supp(\gamma)\), is the formal sum of the pseudo-straight segments that appear in the spine of \(\gamma\).
  The \emph{support} of a multi-curve \(\gamma\) is the sum of the supports of its components.
  The \emph{set-theoretic support}, denoted by \(| \supp(\gamma)|\), is the set of pseudo-straight segments that appear with a positive coefficient in \(\supp(\gamma)\).
\end{definition}
\begin{proposition}\label{prop:traintrack}
  Let \(\gamma\) be a multi-curve on a configuration \(\mathbf{a}\).
  \begin{enumerate}
  \item The set \(|\supp(\gamma)|\) is non-crossing and pointed.
    That is, \(|\supp(\gamma)|\) is contained in a ppt.
  \item Given any non-negative integer linear combination \(s\) of a pointed and non-crossing collection of pseudo-straight segments, there exists a unique multi-curve \(\gamma\) such that \(\supp(\gamma) = s\).
  \item Conversely, the curve \(\gamma\) is admissible (does not have the boundary parallel curve as a component) if and only if \(|\supp(\gamma)|\) is contained in a ppt*.
  \end{enumerate}
\end{proposition}
\begin{proof}
  Let \(\gamma = \{c_1, \dots, c_k\}\).
  Choose representatives \(\gamma_i\) of \(c_i\) that are pairwise disjoint.
  Observe that the support of each \(c_i\) is non-crossing and pointed.
  Then (1) follows from the pairwise disjointness.

  We now prove (2).
  Given \(s\), we explicitly construct a \(\gamma\) such that \(\supp(\gamma) = s\).
  The construction is very similar to constructing a multi-curve on a surface from its train track (see, e.g.,~\cite[\S~1.2]{pen.har:92}).
  Recall that \(\PS(\mathbf{a})\) is the set of pseudo-straight segments on \(\mathbf{a}\).
  Suppose
  \[s = \sum_{\alpha \in \PS(\mathbf{a})} s_{\alpha} \cdot \alpha.\]
  For every \(\alpha\),  we begin by drawing \(s_{\alpha}\) many curves representing \(\alpha\).
  We assume that these curves are in close proximity to the segment underlying \(\alpha\), and are pairwise disjoint except at the end points.
  We will construct \(\gamma\) by joining these curves around the end-points (see~\cref{fig:multicurves-from-support}).
  Consider a point \(a \in \mathbf{a}\).
  If there is only one curve incident at \(a\), we do nothing.
  Suppose there is more than one curve incident at \(a\).
  Since the segments appearing in \(s\) are pointed at \(a\), we can choose two extremal curves \(\alpha\) and \(\beta\) such that, locally near \(a\), all the other curves are contained in the conical region spanned by \(\alpha\) and \(\beta\).
  We connect \(\alpha\) and \(\beta\) to form a longer curve that follows \(\alpha\) almost up to \(a\) (stops a little short), takes a circular detour around \(a\) (leaving the conical region), and then joins and follows \(\beta\) (see~\cref{fig:multicurves-from-support}).
  We repeat the procedure for the next pair of extremal curves at \(a\), and for all vertices \(a \in \mathbf{a}\).
  The end result is a multi-curve \(\gamma\) such that \(\supp(\gamma) = s\).
  We leave it to the reader to convince themselves that this is the only multi-curve \(\gamma\) whose support is \(s\).

  To prove (3), assume that \(\gamma\) contains the boundary parallel curve as a component.
  Then \(|\supp(\gamma)|\) contains all external edges, and hence does not lie in a ppt*.
   Conversely, assuming that \(|\supp(\gamma)|\) does not lie in a ppt*, it must contain all external edges.
  From re-construction of \(\gamma\) from \(\supp(\gamma)\) described above, we see that \(\gamma\) contains the boundary parallel curve as a component.
\end{proof}
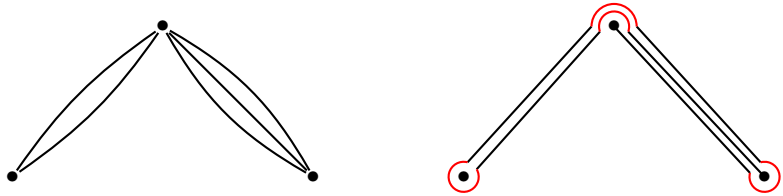
\begin{figure}[ht]
  \centering
  \begin{tikzpicture}
    \node[inner sep=0pt](a) at (0,0) {\(\bullet\)};
    \node[inner sep=0pt](b) at (-2,-2) {\(\bullet\)};
    \node[inner sep=0pt](c) at (2,-2) {\(\bullet\)};
    \coordinate (b1) at ([shift=(75:0.2cm)]-2,-2);
    \coordinate (b2) at ([shift=(30:0.2cm)]-2,-2);    
    \coordinate (a1) at ([shift=(180:0.3cm)]0,0);
    \coordinate (a2) at ([shift=(200:0.2cm)]0,0);
    \coordinate (a3) at ([shift=(-20:0.2cm)]0,0);
    \coordinate (a4) at ([shift=(0:0.3cm)]0,0);    
    \coordinate (c1) at ([shift=(165:0.2cm)]2,-2);
    \coordinate (c2) at ([shift=(105:0.2cm)]2,-2);    
    \draw[thick] (a1) -- (b1);
    \draw[thick] (a2) -- (b2);
    \draw[thick] (a.center) -- (c1);
    \draw[thick] (a3) -- (c.center);
    \draw[thick] (a4) -- (c2);        
    \draw[thick, red] (b1) arc (75:390:0.2cm);
    \draw[thick, red] (c1) arc (165:465:0.2cm);
    \draw[thick, red] (a3) arc (-20:200:0.2cm);
    \draw[thick, red] (a4) arc (0:180:0.3cm);
    \begin{scope}[xshift=-6cm]
          \node[inner sep=0pt](a) at (0,0) {\(\bullet\)};
    \node[inner sep=0pt](b) at (-2,-2) {\(\bullet\)};
    \node[inner sep=0pt](c) at (2,-2) {\(\bullet\)};
    \coordinate (b1) at ([shift=(75:0.2cm)]-2,-2);
    \coordinate (b2) at ([shift=(30:0.2cm)]-2,-2);    
    \coordinate (a1) at ([shift=(180:0.3cm)]0,0);
    \coordinate (a2) at ([shift=(200:0.2cm)]0,0);
    \coordinate (a3) at ([shift=(-20:0.2cm)]0,0);
    \coordinate (a4) at ([shift=(0:0.3cm)]0,0);    
    \coordinate (c1) at ([shift=(165:0.2cm)]2,-2);
    \coordinate (c2) at ([shift=(105:0.2cm)]2,-2);
    \draw[thick] (a) edge [bend left=10] (b) (a) edge [bend right=10] (b);
    \draw[thick] (a) edge [bend left=15] (c)
                 (a) edge [bend right=15] (c) (a) edge (c);
    \end{scope}
  \end{tikzpicture}
  \caption[Constructing a multi-curve from its support]{We construct a multi-curve from its support by first drawing the correct number of underlying pseudo-straight segments (black) and then joining them by circular detours (red) at the end-points.}
  \label{fig:multicurves-from-support}
\end{figure}

As a consequence of this proposition, the support of a multi-curve is naturally an element of \(\widetilde{K}(\mathbf{a})\).
We use the same notation to denote its image in \(|K(\mathbf{a})|\).
\begin{remark}\label{rem:ppt}
Let \(\partial \in |K(\mathbf{a})|\) be the support of the boundary parallel curve.
  Then it is not hard to see that \(|K(\mathbf{a})|\) is piecewise-linearly homeomorphic to the cone over \(|K^{*}(\mathbf{a})|\) via a homeomorphism that sends \(\partial\) to the cone point.
\end{remark}

Given a configuration \(\mathbf{a}\), let \(M(\mathbf{a})\) (resp. \(M^{*}(\mathbf{a})\)) be the set of all multi-curves (resp. admissible multi-curves) on \(\mathbf{a}\).
The map \(\gamma \to \supp(\gamma)\) induces maps
\( M(\mathbf{a}) \to |K(\mathbf{a})| \text{ and } M^{*}(\mathbf{a}) \to |K^{*}(\mathbf{a})|.\)

\begin{proposition}\label{prop:densecurves}
  The image of \(M(\mathbf{a})\) is a dense subset of \(|K(\mathbf{a})|\) (and likewise for \(M^{*}\) and \(K^{*}\)).
\end{proposition}
\begin{proof}
  The image of \(M(\mathbf{a})\) is the set of points of \(|K(\mathbf{a})|\) with rational coordinates, which is dense.
\end{proof}

\begin{remark}\label{prop:densearcs}
  Recall that \(S(\mathbf{a})\) is the set of arcs on \(\mathbf{a}\).
  One can prove that the map \(S(\mathbf{a}) \to |K^{*}(\mathbf{a})|\) is injective and has a dense image.
  Both injectivity and density follow from analysing the curve-joining procedure in the proof of \cref{prop:traintrack}. 
  Injectivity is easy; it follows from the observation that if \(\gamma\) is an arc then for any \(n \geq 2\), the curve with support \(n \supp(\gamma)\) is not an arc.
  Density is more involved, so we omit its proof.
\end{remark}

We now describe how \(K^{*}(\mathbf{a})\) and \(K(\mathbf{a})\) vary with \(\mathbf{a}\).
More precisely, we show that a deformation from \(\mathbf{a}\) to \(\mathbf{b}\) induces a piecewise linear homeomorphism from \(|K(\mathbf{a})|\) to \(|K(\mathbf{b})|\) that restricts to a piecewise linear homeomorphism from \(|K^\ast(\mathbf{a})|\) to \(|K^\ast(\mathbf{b})|\).

Let \(\Gamma \colon [0,1] \to \BigConf_{n+1}\) be a continuous path with \(\Gamma(0) = \mathbf{a}\) and \(\Gamma(1) = \mathbf{b}\).
We can deform arcs, curves, and multi-curves along \(\Gamma\).
Such a deformation induces a bijection from arcs, curves, and multi-curves on \(\mathbf{a}\) to arcs, curves, multi-curves on \(\mathbf{b}\), respectively.

\begin{figure}[h]
  \centering
  \begin{tikzpicture}[inner sep=0em, outer sep=0em]
      \draw[dashed, <->] (-2.5,0) -- (2.5,0);
      \draw (0, -0.5) node (a) {\(\widetilde{\mathbf{a}}\)};
      \draw (-2, 0) node(p1) {\(\bullet\)} (-1, 0) node(p2) {\(\bullet\)} (0, 0.5) node(p3) {\(\bullet\)} (1,0) node(p4) {\(\bullet\)} (2,0) node(p5) {\(\bullet\)};
    \draw[green!70!black,thick] (p5) edge (p3);
    \draw[orange, thick] (p1) edge[in=180,out=-20] (-1,-0.2) (-1,-0.2) edge [out=0,in=220] (p3);
    \draw[red, thick] (p5) edge [in=0, out=150] (0,0.7) (0,0.7) edge[out=180, in=50] (p2);
    
    \begin{scope}[xshift=-7cm]
      \draw[dashed, <->] (-2.5,0) -- (2.5,0);
      \draw (0, -0.5) node (b) {\(\mathbf{a}\)};
      \draw (-2, 0) node(q1) {\(\bullet\)} (-1, 0) node(q2) {\(\bullet\)} (0, 0) node(q3) {\(\bullet\)} (1,0) node(q4) {\(\bullet\)} (2,0) node(q5) {\(\bullet\)}; 
      \draw[green!70!black,thick] (q5) edge[bend right=15] (q3);
      \draw[orange, thick] (q1) edge [bend right=15] (q3);
      \draw[red, thick] (q5) edge [bend right=30] (q2);
    \end{scope}
  \end{tikzpicture}
  \caption{Under an elementary deformation from \(\mathbf{a}\) to \(\widetilde{\mathbf{a}}\), pseudo-straight segments may (green) or may not (orange, red) deform to pseudo-straight segments. }
  \label{fig:deformation}
\end{figure}
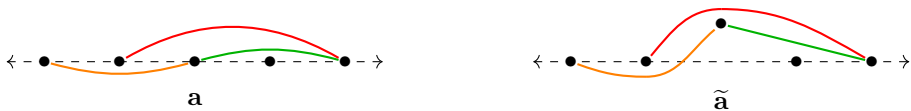

We say that \emph{relative positions of points remain unchanged} along \(\Gamma\) if for all indices \(i,j,k\) and every \(t \in [0,1]\), the points \(\Gamma(t)_i, \Gamma(t)_j, \Gamma(t)_{k}\) are collinear if and only if they are collinear for \(t = 0\).
If relative positions of points remain unchanged along \(\Gamma\), then a pseudo-straight segment in \(\mathbf{a}\) deforms along \(\Gamma\) to a pseudo-straight segment in \(\mathbf{b}\).
This deformation induces a simplicial isomorphism \(K(\mathbf{a}) \to K(\mathbf{b})\).
Let \(M(\mathbf{a}) \to M(\mathbf{b})\) be the bijection obtained by deforming (multi)-curves along \(\Gamma\).
In this case, it is easy to see that the following diagram commutes:
\[
  \begin{tikzcd}
    M(\mathbf{a})\ar{d}{\supp}\ar{r}{\sim} & M(\mathbf{b}) \ar{d}{\supp} \\
    {|K(\mathbf{a})|}\ar{r}{\sim} & {|K(\mathbf{b})|}.
  \end{tikzcd}
\]  
Our goal is to obtain a similar picture for paths that change relative positions.
To do so, we analyse particular paths called elementary deformations.
An elementary deformation only moves one point in a way that does not create any additional collinearities.
The precise definition follows.
\begin{definition}
  Let \(\ell \in \{0,\dots,n\}\).
  A continuous path \(\Gamma \colon [0,1] \to \BigConf_{n+1}\) is an \emph{elementary deformation} at \(\ell\) if
  \begin{enumerate}
  \item \(\Gamma(t)_i\) is constant for all \(i \neq \ell\), and
  \item for all distinct \(i, j \in \{0,\dots, n\} - \{\ell\}\),  if \(\Gamma(0)_{\ell}\) does not lie on the line spanned by \(\Gamma(0)_i\) and \(\Gamma(0)_j\), then neither does \(\Gamma(t)_{\ell}\) for any \(t \geq 0\).
  \end{enumerate}
  The reverse of an elementary deformation is an \emph{elementary degeneration}.
\end{definition}
Fix an elementary deformation \(\Gamma\) with \(\mathbf{a} = \Gamma(0)\) and \(\widetilde{\mathbf{a}} = \Gamma(1)\) (see, e.g.~\cref{fig:deformation}).
A pseudo-straight segment on \(\widetilde{\mathbf{a}}\) degenerates to a pseudo-straight segment on \(\mathbf{a}\), yielding a map \[\degen \colon \PS(\widetilde{\mathbf{a}}) \to \PS({\mathbf{a}}).\]
To understand what happens under deformation, we introduce the notion of a critical line.
For distinct \(i,j \in \{0, \dots, n\}- \{\ell\}\), we say that the line spanned by \(a_i\) and \(a_j\) is \emph{critical} if it contains \(a_\ell\).
Let \(\gamma\) be a pseudo-straight segment on \(\mathbf{a}\).
Unless both end points of \(\gamma\) lie on the same critical line, \(\gamma\) deforms to a pseudo-straight segment on \(\widetilde{\mathbf{a}}\).
If both end points of \(\gamma\) lie on a critical line, it is possible that \(\gamma\) does not deform to a pseudo-straight segment on \(\widetilde{\mathbf{a}}\) (see~\cref{fig:deformation}).

For \(\gamma \in \PS(\mathbf{a})\), let \(\deform(\gamma)\) be the curve obtained by deforming \(\gamma\).
Extend the map \(\gamma \mapsto \supp \deform(\gamma)\) linearly to obtain
\begin{equation}\label{eq:kappa-tilde}
  \kappa \colon \mathbf{R}^{\PS(\mathbf{a})} \to \mathbf{R}^{\PS(\widetilde{\mathbf{a}})}.
\end{equation}

We make three observations.
First, the following diagram commutes
\begin{equation}\label{eq:kappa-preserves-support}
  \begin{tikzcd}
    M(\mathbf{a})\ar{d}{\supp}\ar{r}{\deform} & M(\widetilde{\mathbf{a}}) \ar{d}{\supp} \\
    {\mathbf{R}_{\geq 0}^{\PS(\mathbf{a})}}\ar{r}{\kappa} & {\mathbf{R}_{\geq 0}^{\PS(\widetilde{\mathbf{a}})}}.
  \end{tikzcd}
\end{equation}
Second, if \(S \subset \PS(\mathbf{a})\) is non-crossing and pointed, then the union of the set-theoretic supports of \(\{\deform(\gamma) \mid \gamma \in S\}\) is also non-crossing and pointed.
In particular, if \(S\) is a ppt on \(\mathbf{a}\), then there exists a ppt \(T\) on \(\widetilde{\mathbf{a}}\) such that \(\kappa\) maps \(\mathbf{R}_{\geq 0}^S\) to \(\mathbf{R}_{\geq 0}^T\).
Third, if \(S\) is a ppt* on \(\mathbf{a}\), then the same holds for a ppt* \(T\) on \(\widetilde{\mathbf{a}}\).
The last statement follows from chasing the boundary parallel curve in the diagram \eqref{eq:kappa-preserves-support}.
As a result of the three observations, we see that the map \(\kappa\) induces maps
\begin{equation}\label{def:kappa}
  \begin{tikzcd}
    {|K(\mathbf{a})|} \ar{r}{\kappa}& {|K(\widetilde{\mathbf{a}})|}\\
    {|K^{*}(\mathbf{a})|} \ar{r}{\kappa}\ar[sloped, phantom]{u}{\subset}& {|K^{*}(\widetilde{\mathbf{a}})|\ar[sloped, phantom]{u}{\subset}}.
  \end{tikzcd}
\end{equation}

\subsection{Homeomorphisms between simplicial complexes of different point configurations}
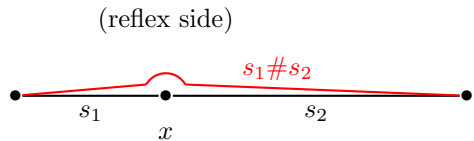
\begin{figure}[ht]
  \centering
  \begin{tikzpicture}
    \node[inner sep=0pt] (p1) at (-2,0) {\(\bullet\)};
    \node[inner sep=0pt] (p2) at (0,0) {\(\bullet\)};
    \draw (0, 1.0) node {(reflex side)};
    \draw (0,-0.5) node {\(x\)};
    \node[inner sep=0pt] (p3) at (4,0) {\(\bullet\)};
    \coordinate(q1) at ([shift=(150:0.3cm)]p2);
    \coordinate(q2) at ([shift=(30:0.3cm)]p2);    
    \draw[thick] (p1) -- node[below]{\(s_1\)} (p2);
    \draw[thick] (p2) -- node[below]{\(s_2\)} (p3);
    \draw[thick,red] (p1) -- (q1);
    \draw[thick, red] (q2) arc (30:150:0.3cm);
    \draw[thick,red] (q2) -- node[above left]{\(s_1 \# s_2\)} (p3);
  \end{tikzpicture}
  \caption{We join two opposing segments \(s_1\) and \(s_2\) to get a longer segment \(s_1 \# s_2\).}
  \label{fig:hash-join}
\end{figure}

Let \(\Gamma\) be an elementary deformation from \(\Gamma(0) = \mathbf{a}\) to \(\Gamma(1) = \widetilde{\mathbf{a}}\).
Our goal is to show that the maps \(\kappa\) defined above are piecewise-linear homeomorphisms.
(See~\Cref{fig:wall-crossing} for the simplest example of an elementary degeneration followed by an elementary deformation.)
For the proof, we need a construction to combine opposing pseudo-straight segments.
\begin{definition}[The \(\#\) construction]
  Suppose \(s_1\) and \(s_2\) are pseudo-straight segments on \(\mathbf{a}\) that lie on the same critical line \(L\), share a common end-point \(x\), and approach \(x\) from opposite sides.
  Let \(\widetilde x \in \widetilde{\mathbf{a}}\) be the deformation of \(x \in \mathbf{a}\), and set \(\widetilde s_i = \deform(s_i)\).
  Assume that the spines of \(\widetilde{s}_1\) and \(\widetilde s_2\) do not remain collinear around \(\widetilde x\).
  The support of \(\widetilde s_1\) and \(\widetilde s_2\) spans a convex cone at \(\widetilde x\) whose limit under the degeneration to \(x\) is a half-plane bounded by \(L\).
  By the \emph{reflex side} of \(L\) at \(x\), we mean the opposite half-plane.
  Let \(s_1 \# s_2\) be the pseudo-straight segment on \(\mathbf{a}\) obtained by following \(s_1\) almost up to \(x\), taking a small semi-circular detour around \(x\) on the reflex side, and then following \(s_2\) (see \cref{fig:hash-join}).
\end{definition}
Note that the spine of the deformation of \(s_1 \# s_2\) is the concatenation of the spines of the deformations of \(s_1\) and \(s_2\).

Let \(\gamma \in \PS(\mathbf{a})\), and let \((\widetilde\gamma_1, \dots, \widetilde\gamma_k)\) be the spine of the deformation of \(\gamma\) to \(\widetilde{\mathbf{a}}\).
Let \(\gamma_i \in \PS(\mathbf{a})\) be  the degeneration of \(\widetilde\gamma_i \in \PS(\mathbf{a})\).
Then all \({\gamma}_i\) must lie in order on a single critical line \(L\).
Furthermore, we can reconstruct \(\gamma\) from \(({\gamma}_1, \dots, {\gamma}_k)\) as follows:
\[ \gamma = {\gamma}_1 \# {\gamma}_2 \# \cdots \# {\gamma}_k.\]
It is easy to see that the \(\#\) construction is associative.
But we caution the reader that the reflex side of \(L\) may be different for different intermediate points.

We now take up the task of proving that \(\kappa\) from \eqref{def:kappa} is a piecewise-linear homeomorphism.
We fix a ppt \(T\) on \(\widetilde{\mathbf{a}}\), and argue that the ppts \(S\) on \(\mathbf{a}\) such that \(\kappa(S) \subset T\) give a linear subdivision of \(|T|\).
To do so, we need an auxiliary notion of a faux-ppt.

\begin{definition}
  Let \(T \subset \PS (\widetilde{\mathbf{a}})\) be a ppt.
  A \emph{faux-ppt} for \(T\) is a subset \(U \subset \PS(\mathbf{a})\) satisfying the following conditions.
  \begin{enumerate}
  \item The elements of \(U\) are non-crossing.
  \item The cardinality of \(U\) is the same as the cardinality of \(T\).
  \item For every \(u \in U\), the curve \(\deform(u)\) is supported on \(T\).
  \item The map \(\kappa \colon \mathbf{R}^U \to \mathbf{R}^T\) induced by \(u \mapsto \supp \deform(u)\) is injective.
  \end{enumerate}
\end{definition}
Recall the degeneration map \(\degen \colon \PS(\widetilde{\mathbf{a}}) \to \PS(\mathbf{a})\).
Then \(U = \degen(T)\) is a faux-ppt for \(T\).

Note that \(\kappa \colon \mathbf{R}^U \to \mathbf{R}^T\) maps \(\mathbf{R}_{\geq 0}^U\) to \(\mathbf{R}_{\geq 0}^T\).
Recall that 
\[ |U| = \left(\mathbf{R}_{\geq 0}^U - \{0\}\right) / \textrm{scaling}, \text{ and } |T| = \left(\mathbf{R}_{\geq 0}^{T} - \{0\}\right) / \textrm{scaling}. \]
Thus \(\kappa\) induces a map \(|U| \to |T|\), which we also call \(\kappa\).
The image of \(|U|\) is full-dimensional in \(|T|\) by the injectivity assumption on \(\kappa\).
So, \(U\) cannot also be a faux-ppt for a different ppt \(T' \subset \PS(\mathbf{a})\).

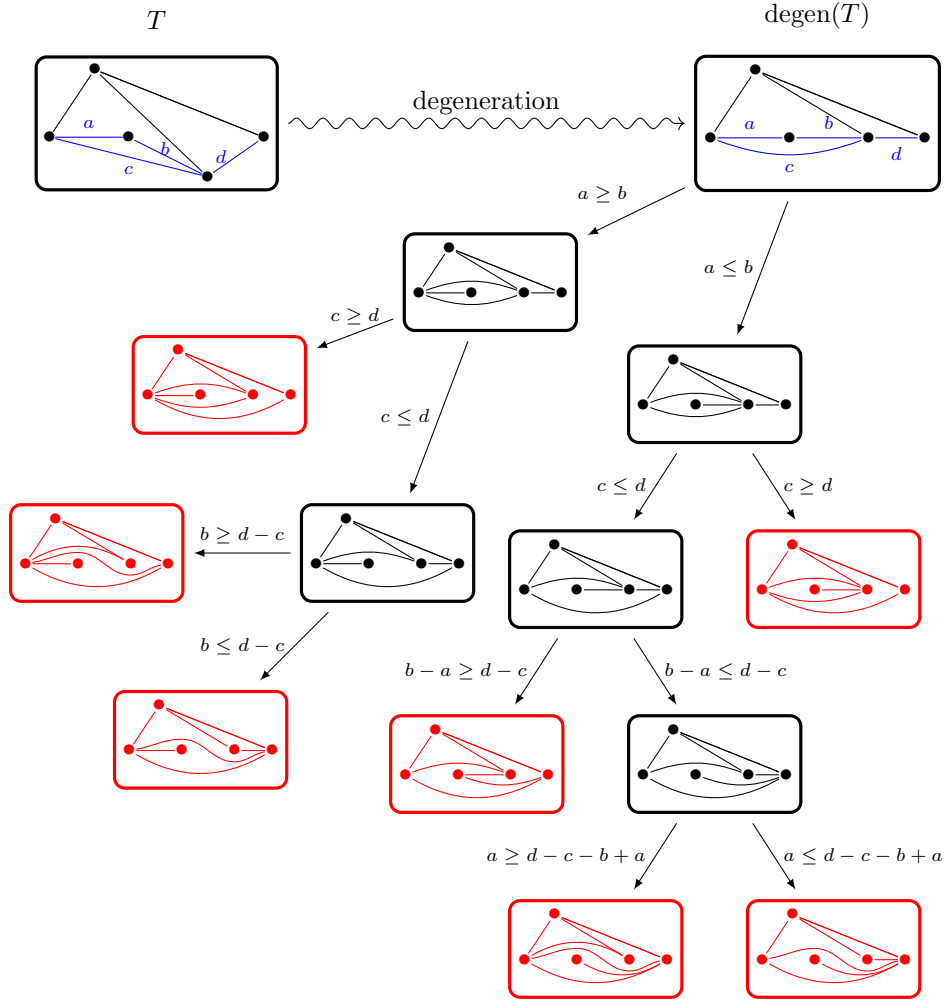
\begin{figure}[h]
    \centering
  \begin{tikzpicture}
    [sibling distance=9em,
    edge from parent/.style={draw,-latex},
    grow=south west, level distance=10em      
    ]
    \tikzset{
      ppt/.pic={
        \node[inner sep=0pt] (-pt1) at (0,0) {\(\bullet\)};
        \node[inner sep=0pt] (-pt2) at (-0.4,-0.6) {\(\bullet\)};
        \node[inner sep=0pt] (-pt3) at (0.3,-0.6) {\(\bullet\)};
        \node[inner sep=0pt] (-pt4) at (1,-0.6) {\(\bullet\)};
        \node[inner sep=0pt] (-pt5) at (1.5,-0.6) {\(\bullet\)};
        \coordinate (-pt3top) at (0.3, -0.45);
        \coordinate (-pt3bot) at (0.3, -0.75);
        \coordinate (-pt4top) at (1, -0.45);
        \coordinate (-pt4bot) at (1, -0.75);                
        \coordinate(phantom) at (1.5,-1);
        \path
        (-pt1) edge (-pt2)
        (-pt1) edge (-pt4)
        (-pt1) edge (-pt5)
        (-pt1) edge (-pt5);                                        
        \node[draw, solid, very thick, rounded corners, inner sep=1mm, fit=(-pt1) (-pt2) (phantom) (-pt5)] {};
      }}
    \tikzset{
      nondegppt/.pic={
        \node[inner sep=-0.5pt] (-pt1) at (0,0) {\(\bullet\)};
        \node[inner sep=-0.5pt] (-pt2) at (-0.4,-0.6) {\(\bullet\)};
        \node[inner sep=-0.5pt] (-pt3) at (0.3,-0.6) {\(\bullet\)};
        \node[inner sep=-0.5pt] (-pt4) at (1,-0.95) {\(\bullet\)};
        \node[inner sep=-0.5pt] (-pt5) at (1.5,-0.6) {\(\bullet\)};
        \coordinate(phantom) at (1.5,-1);
        \path
        (-pt1) edge (-pt2)
        (-pt1) edge (-pt4)
        (-pt1) edge (-pt5)
        (-pt1) edge (-pt5)
        (-pt2) edge[blue] node[above, midway, font=\scriptsize]{\(a\)} (-pt3)
        (-pt2) edge[blue] node[below, midway, font=\scriptsize]{\(c\)} (-pt4)
        (-pt3) edge[blue] node[right, near start, font=\scriptsize]{\(b\)} (-pt4)
        (-pt4) edge[blue] node[left, midway, font=\scriptsize]{\(d\)} (-pt5)
        ;                                        
        \node[draw, very thick, rounded corners, inner sep=1mm, fit=(-pt1) (-pt2) (phantom) (-pt5)] {};
      }}
    \node[matrix] (nondeg) {\pic[scale=1.5] {nondegppt};\\};
    \node[above=0.3em of nondeg] {\(T\)};
\node[matrix,right=15em of nondeg] (deg) {\pic[scale=1.5] (top) {ppt};\\}
child[level distance=13em] {node[matrix] {\pic (r2l) {ppt};\\}
      child[level distance=10em]
      {node[matrix] {\pic[red] (r2l3l) {ppt};\\}
        edge from parent node[above, font=\scriptsize] {\(c \geq d\)}
      }
      child[level distance=10em]
      {node[matrix] {\pic (r2l3r) {ppt};\\}
        child [grow=west,level distance=11em] {node[matrix]{\pic[red](r2l3r3l) {ppt};\\}
          edge from parent node[above,font=\scriptsize] {\(b \geq d-c\)}}
        child [grow=south west] {node[matrix]{\pic[red](r2l3r3r) {ppt};\\}
          edge from parent node[left,font=\scriptsize] {\(b \leq d-c\)}}          
        edge from parent node[left,font=\scriptsize] {\(c \leq d\)}
      }
      edge from parent node[above left, font=\scriptsize] {\(a \geq b\)}
    }
child {node[matrix] {\pic (r2r) [inner sep=0pt] {ppt};\\}
      [grow=270, level distance=7em]
child {node[matrix] {\pic (r2r3r) {ppt};\\}
        [grow=270, level distance=7em]
        child
        {node[matrix] {\pic[red] (r2r3r3l) {ppt};\\}
          edge from parent node[left, font=\scriptsize] {\(b-a \geq d-c\)}
        }
        child
        {node[matrix] {\pic (r2r3r3r) {ppt};\\}
          child {node[matrix]{\pic[red](r2r3r3r3l) {ppt};\\}
            edge from parent node[left, font=\scriptsize] {\(a \geq d-c-b+a \)}}
          child {node[matrix]{\pic[red](r2r3r3r3r) {ppt};\\}
            edge from parent node[right, font=\scriptsize] {\(a \leq d-c-b+a \)}}
          edge from parent node[right, font=\scriptsize] {\(b-a \leq d-c \)}            
        }
        edge from parent node[left, font=\scriptsize] {\(c \leq d\)}}
child
      {node[matrix] {\pic[red] (r2r3l) {ppt};\\}
        edge from parent node[right, font=\scriptsize] {\(c \geq d\)}}
      edge from parent node[left, font=\scriptsize] {\(a \leq b\)}
    }
    ;
    \node[above=0.3em of deg] {\(\degen(T)\)};
    \draw[decorate,decoration={snake,amplitude=2,segment length=10,post length=1mm},->] 
    (nondeg) -- (deg) node[midway, above]{degeneration};
    
    \path[blue]
    (top-pt2) edge node[above, font=\scriptsize]{\(a\)} (top-pt3)
    (top-pt3) edge node[above, font=\scriptsize]{\(b\)} (top-pt4)
    (top-pt2) edge[bend right=20] node[below, font=\scriptsize]{\(c\)} (top-pt4)
    (top-pt4) edge node[below, font=\scriptsize]{\(d\)} (top-pt5);

    \path[black]
    (r2l-pt2) edge (r2l-pt3)
    (r2l-pt2) edge[bend right=20] (r2l-pt4)
    (r2l-pt2) edge[bend left=20] (r2l-pt4)      
    (r2l-pt4) edge (r2l-pt5);

    \path[black]
    (r2r-pt3) edge (r2r-pt4)
    (r2r-pt2) edge[bend left=20] (r2r-pt4)    
    (r2r-pt2) edge[bend right=20] (r2r-pt4)
    (r2r-pt4) edge (r2r-pt5);

    \path[red]
    (r2l3l-pt2) edge (r2l3l-pt3)
    (r2l3l-pt2) edge[bend left=20] (r2l3l-pt4)    
    (r2l3l-pt2) edge[bend right=20] (r2l3l-pt4)
    (r2l3l-pt2) edge[bend right] (r2l3l-pt5);

    \path[black]
    (r2l3r-pt2) edge (r2l3r-pt3)
    (r2l3r-pt2) edge[bend left=20] (r2l3r-pt4)
    (r2l3r-pt4) edge (r2l3r-pt5)
    (r2l3r-pt2) edge[bend right] (r2l3r-pt5);

    \path[red]
    (r2l3r3l-pt2) edge (r2l3r3l-pt3)
    (r2l3r3l-pt2) edge[out=10,in=180] (r2l3r3l-pt3top)
    (r2l3r3l-pt3top) edge[in=180,out=0] (r2l3r3l-pt4bot)
    (r2l3r3l-pt4bot) edge[in=190,out=0] (r2l3r3l-pt5)
    (r2l3r3l-pt2) edge[bend left=30] (r2l3r3l-pt4)
    (r2l3r3l-pt2) edge[bend right] (r2l3r3l-pt5);

    \path[red]
    (r2l3r3r-pt2) edge (r2l3r3r-pt3)
    (r2l3r3r-pt2) edge[out=10,in=180] (r2l3r3r-pt3top)
    (r2l3r3r-pt3top) edge[in=180,out=0] (r2l3r3r-pt4bot)
    (r2l3r3r-pt4bot) edge[in=190,out=0] (r2l3r3r-pt5)
    (r2l3r3r-pt4) edge (r2l3r3r-pt5)
    (r2l3r3r-pt2) edge[bend right] (r2l3r3r-pt5);
    
    \path[red]
    (r2r3l-pt3) edge (r2r3l-pt4)
    (r2r3l-pt2) edge[bend left=20] (r2r3l-pt4)    
    (r2r3l-pt2) edge[bend right=20] (r2r3l-pt4)
    (r2r3l-pt2) edge[bend right] (r2r3l-pt5);

    \path
    (r2r3r-pt3) edge (r2r3r-pt4)      
    (r2r3r-pt4) edge (r2r3r-pt5)
    (r2r3r-pt2) edge[bend left=20] (r2r3r-pt4)
    (r2r3r-pt2) edge[bend right] (r2r3r-pt5);
    
    \path[red]
    (r2r3r3l-pt3) edge (r2r3r3l-pt4)
    (r2r3r3l-pt2) edge[bend left=20] (r2r3r3l-pt4)    
    (r2r3r3l-pt3) edge[bend right=20] (r2r3r3l-pt5)
    (r2r3r3l-pt2) edge[bend right] (r2r3r3l-pt5);
    
    \path
    (r2r3r3r-pt4) edge (r2r3r3r-pt5)
    (r2r3r3r-pt2) edge[bend left=20] (r2r3r3r-pt4)    
    (r2r3r3r-pt3) edge[bend right=20] (r2r3r3r-pt5)
    (r2r3r3r-pt2) edge[bend right] (r2r3r3r-pt5);

    \path[red]
    (r2r3r3r3l-pt3) edge[bend right=30] (r2r3r3r3l-pt5)
    (r2r3r3r3l-pt2) edge[out=10,in=180] (r2r3r3r3l-pt3top)
    (r2r3r3r3l-pt3top) edge[in=180,out=0] (r2r3r3r3l-pt4bot)
    (r2r3r3r3l-pt4bot) edge[in=190,out=0] (r2r3r3r3l-pt5)
    (r2r3r3r3l-pt2) edge[bend left=30] (r2r3r3r3l-pt4)
    (r2r3r3r3l-pt2) edge[bend right] (r2r3r3r3l-pt5);
    
    \path[red]
    (r2r3r3r3r-pt3) edge[bend right=30] (r2r3r3r3r-pt5)
    (r2r3r3r3r-pt2) edge[out=10,in=180] (r2r3r3r3r-pt3top)
    (r2r3r3r3r-pt3top) edge[in=180,out=0] (r2r3r3r3r-pt4bot)
    (r2r3r3r3r-pt4bot) edge[in=190,out=0] (r2r3r3r3r-pt5)
    (r2r3r3r3r-pt4) edge (r2r3r3r3r-pt5)
    (r2r3r3r3r-pt2) edge[bend right] (r2r3r3r3r-pt5);
    
  \end{tikzpicture}
  \caption[PL map of a degeneration]{We split a faux-ppt in two faux-ppts by fusing a pair of extremal edges.
    Starting with a ppt \(T\) on \(\widetilde{\mathbf{a}}\), we apply this procedure applied repeatedly to the faux-ppt \(\degen(T)\).
    The result is a tree of faux-ppts whose leaves are the ppts of \(\mathbf{a}\) whose deformations are supported on \(T\).}
    \label{fig:pl-resolution-tree}
  \end{figure}

A faux-ppt may or may not be pointed.
If it is not pointed, we describe a procedure that ``splits'' it in two.
The branching of each node in \cref{fig:pl-resolution-tree} is an example of such a split.

Fix a ppt \(T \subset \PS(\widetilde{\mathbf{a}})\) and a faux-ppt \(U \subset \PS(\mathbf{a})\) for \(T\).
Let \(s_1\) and \(s_2\) be two pseudo-straight segments in \(U\) that lie on the same critical line \(L\), share an endpoint \(x\), and come to \(x\) from opposite directions.
Let \(\widetilde x \in \widetilde{\mathbf{a}}\) be the deformation of \(x \in \mathbf{a}\).
The supports of \(\deform(s_1)\) and \(\deform(s_2)\) are pointed at \(\widetilde x\).
We say that \(s_1\) and \(s_2\) are \emph{extremal opposing edges} if all other edges of \(T\) incident to \(\widetilde x\) lie in the cone spanned by \(\deform(s_1)\) and \(\deform(s_2)\) at \(\widetilde x\).
Suppose \(s_1\) and \(s_2\) are extremal opposing edges at \(x\).
Set \(U_1 = U \cup \{s_1 \# s_2\} - \{s_2\}\) and \(U_2 = U \cup \{s_1 \# s_2\} - \{s_1\}\).

\begin{lemma}\label{lem:faux-ppt-split}
  For \(i = 1,2\), the set \(U_i\) is a faux-ppt for \(T\).
  Furthermore, \(\kappa(|U_1|)\) and \(\kappa(|U_2|)\) give a linear subdivision of \(\kappa(|U|)\).
\end{lemma}
\begin{proof}
  We check the four conditions in the definition of a faux-ppt.
  By the extremality of the pair, \(s_1 \# s_2\) does not cross any element of \(U\).
  So the first condition holds.
  Note that \(s_1 \# s_2\) cannot already lie in \(U\); if it did, then \(s_1 \# s_2 - s_1 - s_2\) would be in the kernel of \(\kappa \colon \mathbf{R}^U \to \mathbf{R}^{T}\).
  As a result, both \(U_i\) have the same cardinality as \(U\).
  So the second condition holds.
  The spine of the deformation of \(s_1 \# s_2\) is the concatenation of the spines of the deformations of \(s_{1}\) and \(s_2\).
  So the third condition holds.
  Let \(V = U - \{s_1, s_2\}\).
  Note that
  \[\mathbf{R}^{U_i} = \mathbf{R}^V \times \mathbf{R}^{\{s_i, s_1 \# s_2\}} \quad \text{and} \quad \mathbf{R}^U = \mathbf{R}^V \times \mathbf{R}^{\{s_1, s_2\}}.\]
  We have the map
  \(\phi_i \colon \mathbf{R}^{U_i} \to \mathbf{R}^U\)
  that is the identity on \(\mathbf{R}^V\), sends \(s_i\) to \(s_i\), and sends \(s_1 \# s_2\) to \(s_1 + s_2\).
  Then \(\kappa \colon \mathbf{R}^{U_i} \to \mathbf{R}^{T}\) is the composite of \(\phi_i \colon \mathbf{R}^{U_i} \to \mathbf{R}^U\) and \(\kappa \colon \mathbf{R}^U \to \mathbf{R}^T\).
  Since \(\phi_i\) is injective, the fourth condition holds.

  Finally, it is clear that the images of \(\phi_{i} \colon |U_i| \to |U|\) already give a linear subdivision of \(|U|\).
  The last statement follows by applying \(\kappa\).
\end{proof}

\begin{definition}
  Let \(U \subset \PS(\mathbf{a})\) be a faux-ppt.
  We say that a pseudo-straight segment \(v \in \PS(\mathbf{a})\) \emph{resolves} \(U\) if there exists a chain \((\gamma_1, \dots, \gamma_k)\) of edges in \(U\) such that
  \[ v = \gamma_1 \# \cdots \# \gamma_k.\]
We say that a subset \(V\subset \PS(\mathbf{a})\) resolves \(U\) if each \(v \in V\) resolves \(U\).
\end{definition}
\begin{remark}\label{rem:resolve-consequences}
Implicit in the definition above is that the \(\gamma_i\) are pseudo-straight segments that lie in order on the same line, and at each intermediate end-point, their deformations bend.
Observe the following.
\begin{enumerate}
\item Suppose that \(v\) resolves a faux-ppt \(U\) for \(T\).
  Since \(\kappa \colon \mathbf{R}^U \to \mathbf{R}^{T}\) is injective and \(\kappa(v) = \sum \kappa(\gamma_i)\), the chain \((\gamma_1, \dots, \gamma_k)\) is unique.
\item If \(U\) is already a ppt, then it contains no opposing pseudo-straight segments, so there are no non-trivial \(\#\) joins of pseudo-straight segments of \(U\).
  Thus any \(v \in \PS(\mathbf{a})\) that resolves \(U\) must be an element of \(U\).
\end{enumerate}
\end{remark}

\begin{lemma}\label{lem:resolve-split}
  Let \(U\) be a faux-ppt and let \(U_1\) and \(U_2\) be the faux-ppts obtained by splitting it as in~\cref{lem:faux-ppt-split}.
  Let \(V \subset \PS(\mathbf{a})\) be non-crossing and pointed (for example, a ppt).
  If \(V\) resolves \(U\), then \(V\) resolves \(U_1\) or \(V\) resolves \(U_2\).
  Conversely, if \(V\) resolves \(U_1\) or \(U_2\) then \(V\) resolves \(U\).
\end{lemma}
\begin{proof}
  Take a \(v \in V\), and let \(\gamma = (\gamma_i)\) be the chain on \(U\) with \(v = \gamma_1 \# \cdots \# \gamma_k\).
  If for all \(i\), we have \(\gamma_i \not \in \{s_{1}, s_2\}\), then \(v\) resolves both \(U_1\) and \(U_2\).
  If two consecutive edges \(\{\gamma_i, \gamma_{i+1}\}\) are \(\{s_1, s_2\}\), remove them and replace them by \(s_1 \# s_2\).
  The resulting chain has the same \(\#\)-join, namely \(v\).
  So \(v\) resolves both \(U_1\) and \(U_2\).

  Suppose \(\gamma_i = s_1\), but neither \(\gamma_{i-1}\) nor \(\gamma_{i+1}\) is \(s_2\).
  In this case, we say that \(v\) contains \(s_1\) ``alone'' (and likewise for \(s_2\)).
  If no \(v \in V\) contains \(s_1\) alone, then \(V\) resolves \(U_2\) (and likewise for \(s_2\)).
  We show that \(V\) cannot contain \(v, v'\) such that \(v\) contains \(s_1\) alone and \(v'\) contains \(s_2\) alone.

  Let us analyse \(v \in V\) that contains \(s_1\) alone.
  Let \(\gamma\) be as before, and orient it so that in the induced orientation, \(\gamma_i = s_1\) ends at \(x\).
  There are two possibilities: either \(\gamma\) ends with \(s_1\) at \(x\), or it continues with \(\gamma_{i+1} \neq s_2\) that starts at \(x\).
  In this case, by the extremality of the pair \((s_1, s_2)\), the segment \(\gamma_{i+1}\) must be on the non-reflex side at \(x\).
  Likewise, if \(v'\) contains \(s_2\) alone, then either the corresponding \(\gamma'\) ends at \(x\) or continues with \(\gamma'_{i+1} \neq s_1\) on the non-reflex side at \(x\).

  Suppose for contradiction that \(V\) contains \(v\) and \(v'\) such that the corresponding deformed chains contain \(s_1\) alone and \(s_2\) alone, respectively.
  Then \(v\) and \(v'\) must cross or form an opposing pair at \(x\) (see~\cref{fig:crossing-or-pointed}).
  Since \(V\) is non-crossing and pointed, this is impossible.
  
  The converse is straightforward.
  Suppose \(V\) resolves \(U_1\).
  For a \(v \in V\), if the corresponding chain \(\gamma\) contains \(s_1 \# s_2\), we simply break it into \((s_1, s_2)\).
  This procedure yields a chain in \(U\) with the same \(\#\) join.
\end{proof}
  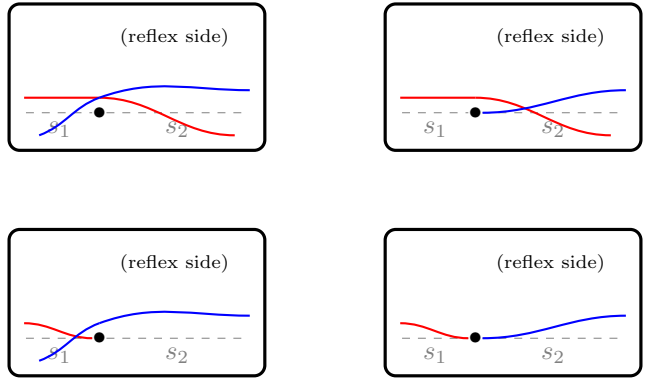
\begin{figure}[ht]
    \centering
    \begin{tikzpicture}
      \tikzset{
        config/.pic={
          \node[inner sep=0pt] (-pt1) at (-1,0) {};
          \node[inner sep=0pt] (-pt2) at (0,0) {\(\bullet\)};
          \node[inner sep=0pt] (-pt3) at (2,0) {};
          \node[font=\scriptsize] (text) at (1,1) {(reflex side)};

          \path[very thin, gray, dashed]
          (-pt1) edge node[midway, below]{\(s_1\)} (-pt2)
          (-pt2) edge node[midway, below]{\(s_2\)} (-pt3);
          \coordinate (-leftend) at (-1,0.2);
          \coordinate (-leftbot) at (-0.8,-0.3);          
          \coordinate (-pt2top) at (0,0.2);
          \coordinate (-rightend) at (2,0.3);                  
          \coordinate (-rightbot) at (1.8,-0.3);
          \node[draw, solid, very thick, rounded corners, inner sep=2mm, fit=(-leftbot) (-leftend) (-rightbot) (-rightend) (text)] {};
        }}
      \node[matrix] at (0,0) {\pic (c1) {config};\\};
      \node[matrix] at (5,0) {\pic (c2) {config};\\};
      \node[matrix] at (0,-3) {\pic (c3) {config};\\};
      \node[matrix] at (5,-3) {\pic (c4) {config};\\};

      \path[thick, red]
      (c1-leftend) edge[out=0, in=180] (c1-pt2top)
      (c1-pt2top) edge[out=0, in=180] (c1-rightbot);
      \path[thick, blue]
      (c1-rightend) edge[out=180, in=20] (c1-pt2top)
      (c1-pt2top) edge[in=20, out=200] (c1-leftbot);

      \path[thick, red]
      (c2-leftend) edge[out=0, in=180] (c2-pt2top)
      (c2-pt2top) edge[out=0, in=180] (c2-rightbot);      
      \path[thick, blue]
      (c2-rightend) edge[out=180, in=0] (c2-pt2);

      \path[thick, red]
      (c3-leftend) edge[out=0, in=180] (c3-pt2);      
      \path[thick, blue]
      (c3-rightend) edge[out=180, in=20] (c3-pt2top)
      (c3-pt2top) edge[in=20, out=200] (c3-leftbot);
      
      \path[thick, red]
      (c4-leftend) edge[out=0, in=180] (c4-pt2);
      \path[thick, blue]
      (c4-rightend) edge[out=180, in=0] (c4-pt2);
    \end{tikzpicture}
    \caption{If two pseudo-straight segments (red and blue) on a critical line are \(\#\) joins of chains that contain \(s_1\) alone and \(s_2\) alone, respectively, then they must cross or oppose.}
    \label{fig:crossing-or-pointed}
  \end{figure}

  \begin{lemma}\label{lem:faux-cardinality}
    Let \(U\) be a faux-ppt on \(\widetilde{\mathbf{a}}\) of cardinality \(2n-1\).
    If \(S\) is a ppt on \(\mathbf{a}\) that resolves \(U\), then the cardinality of \(S\) is also \(2n-1\).
  \end{lemma}
  \begin{proof}
    We induct on the number of pairs of opposing edges of \(U\).
    If \(U\) has no such edges, then \(U\) is non-crossing, pointed, and of maximum possible cardinality (see \cref{lem:edges-upper-bound}).
    So \(U\) is a ppt on \(\mathbf{a}\).
    The only ppt that resolves \(U\) is \(U\) itself.
    The assertions follow.

    If \(U\) has a pair of opposing edges, pick an extremal such pair and let \(U_1\) and \(U_2\) be the faux-ppts obtained by splitting \(U\) as in~\cref{lem:faux-ppt-split}.
    By \cref{lem:resolve-split}, the ppt \(S\) resolves \(U_1\) or \(U_2\).
    We conclude by applying the inductive hypothesis.
  \end{proof}

  As a consequence of \cref{lem:faux-cardinality}, we get the following.
  \begin{proposition}\label{prop:ppt-cardinality}
    Let \(\mathbf{a}\) be any point configuration.
    Every ppt on \(\mathbf{a}\) contains \(2n-1\) edges.
  \end{proposition}
  \begin{proof}
    Consider an elementary deformation \(\widetilde{\mathbf{a}}\) of \(\mathbf{a}\).
    Let \(S\) be a ppt on \(\mathbf{a}\) whose deformation is supported on a ppt \(T\) of \(\widetilde{\mathbf{a}}\).
    \cref{lem:faux-cardinality} applied to \(U = \degen(T)\), we see that if \(T\) has \(2n-1\) edges, so does \(S\).
    By a sequence of elementary deformations, we may deform \(\mathbf{a}\) to a point configuration with no three collinear points.
    For such configurations, it is known that every ppt has \(2n-1\) edges.
    The result follows.
  \end{proof}

  \begin{lemma}\label{lem:subdivision}
    Let \(U\) be a faux-ppt on \(\mathbf{a}\) for a ppt \(T\) on \(\widetilde{\mathbf{a}}\).
    \begin{enumerate}
    \item Let \(p \in |K(\mathbf{a})|\), and let \(V \subset \PS(\mathbf{a})\) be the support of \(p\).
      Then we have \(\kappa(p) \in \kappa(|U|)\) if and only \(V\) resolves \(U\).
    \item Let \(S_1, \dots, S_m\) be ppts on \(\mathbf{a}\) that resolve \(U\).
      Then \(\kappa(|S_1|), \cdots, \kappa(|S_m|)\) are full-dimensional and give a simplicial subdivision of \(\kappa(|U|)\).
    \end{enumerate}
  \end{lemma}
  \begin{proof}
    We induct on the number of pairs of opposing edges in \(U\).

    Let us prove the base case, in which \(U\) has no opposing edges, and is thus a ppt.
    We prove the first statement.
    For one direction, suppose that \(V\) resolves \(U\).
    Then for each \(v \in V\) there is a chain \((\gamma_1, \ldots, \gamma_k)\) in \(U\) such that \(\kappa(v) = \sum_i\kappa(\gamma_i)\), and thus \(\kappa(v) \in \kappa(|U|)\).
    Now \(p\) is a convex linear combination of the elements \(v \in V\), and \(|U|\) is a convex set.
    Thus \(\kappa(p) \in \kappa(|U|)\) as well.
    Let us prove the other direction.
    Recall that \(\kappa\) is induced by a linear map with rational coordinates.
    So it suffices to prove the statement assuming that \(p\) has rational coordinates.
    Let \(\widetilde p \in \mathbf{Z}_{\geq 0}^{\PS(\mathbf{a})}\) be an integer lift of \(p\).

    Let \(q \in |U|\) be such that \(\kappa(q) = \kappa(p)\). 
    Let \(\widetilde{q} \in \mathbf{Z}^{\PS(\mathbf{a})}_{\geq 0}\) be an integer lift of \(q\).
    Recall that \(\kappa \colon |K(\mathbf{a})| \to |K(\widetilde{\mathbf{a}})|\) is induced by a linear map \(\mathbf{R}^{\PS(\mathbf{a})} \to \mathbf{R}^{\PS(\widetilde{\mathbf{a}})}\).
    In~\eqref{eq:kappa-tilde} it was denoted \(\kappa\), but for clarity, we denote that map here as \(\widetilde{\kappa}\).
    
    Since \(\kappa(q) = \kappa(p)\), there exists positive integers \(m\) and \(n\) such that
    \[\widetilde{\kappa}(m\widetilde{q}) = \widetilde{\kappa}(n\widetilde{p}).\]
Let \(\gamma_q\) and \(\gamma_p\) be the unique multi-curves on \(\mathbf{a}\) (as in~\cref{prop:traintrack}) such that
    \[\supp(\gamma_q) = m\widetilde{q}, \text{ and }\supp(\gamma_p) = n\widetilde{p}.\]
    Since \(\widetilde{\kappa}\) commutes with supports of deformations (see~\eqref{eq:kappa-preserves-support}), we have
    \[ \supp \deform (\gamma_q) = \supp \deform(\gamma_p).\]
Since the support determines the multi-curve uniquely, we see that \(\deform (\gamma_q) = \deform (\gamma_p)\).
    But deformation is a bijection on (isotopy classes of) curves.
    So \(\gamma_q = \gamma_p\).
    Note that the set-theoretic support of \(\gamma_p\) is the same as the support of \(p\).
    On the other hand, this is also the same as the set-theoretic support of \(q\), which is a subset of \(U\).
    In particular, the support of \(p\) resolves \(U\).

    For the second statement, note that since \(U\) is a ppt, the only ppt that resolves \(U\) is \(U\) itself (see~\cref{rem:resolve-consequences}), and thus the subdivision asserted is trivial.

    Having done the base case, we do the induction step.
    Choose a pair of extremal opposing edges in \(U\), and split \(U\) as \(U_1\) and \(U_2\) as in \cref{lem:faux-ppt-split}.
    Note that each \(U_i\) has fewer pairs of opposing edges than \(U\).

    We prove the first statement.
    By \cref{lem:faux-ppt-split}, we get that \(\kappa(p) \in \kappa(|U_1|)\) or \(\kappa(p) \in \kappa(|U_2|)\).
    By the inductive hypothesis, the support of \(p\) resolves \(U_1\) or \(U_2\).
    So the support of \(p\) resolves \(U\).

    We now prove the second statement.
    By \cref{lem:resolve-split}, a ppt \(S\) that resolves \(U\) resolves \(U_1\) or \(U_2\).
    By the inductive hypothesis, \(\kappa(|S|)\) is full-dimensional.
    As a result, \(S\) cannot resolve both \(U_1\) and \(U_2\) (if it did, then \(\kappa(|S|) \subset \kappa(|U_1|) \cap \kappa(|U_2|)\), which is not full-dimensional).
    So the set of ppts that resolve \(U\) are partitioned into two disjoint sets: those that resolve \(U_1\) and those that resolve \(U_2\).
    Using the inductive hypothesis and that \(\kappa(|U_1|)\) and \(\kappa(|U_2|)\) form a simplicial subdivision of \(\kappa(|U|)\), we are done.
  \end{proof}
  
 We are now ready to prove the desired result for elementary deformations.
 \begin{proposition}\label{prop:pl-iso}
   Let \(\Gamma\) be an elementary deformation of \(\mathbf{a}\) to \(\widetilde{\mathbf{a}}\).
   The maps
   \[ \kappa \colon |K(\mathbf{a})| \to |K(\widetilde{\mathbf{a}})| \text{ and } \kappa \colon |K^{*}(\mathbf{a})| \to |K^{*}(\widetilde{\mathbf{a}})|\]
   defined by \(\gamma \mapsto \supp \deform \gamma\) are piecewise-linear homeomorphisms.
   In both cases, the images under \(\kappa\) of the simplices in the source give a simplicial subdivision of the target.
 \end{proposition}
 \begin{proof}
   We treat \(K\) first and then \(K^{*}\).
   We know that \(\kappa\) is continuous and sends the simplices of \(|K(\mathbf{a})|\) linearly to the simplices of \(|K(\widetilde{\mathbf{a}})|\).
   It suffices to prove that for every maximal simplex \(\Delta\) of \(|K(\widetilde{\mathbf{a}})|\), the map
   \[ \kappa \from \kappa^{-1}(\Delta) \to \Delta\]
   is a homeomorphism.
   In fact, it suffices to prove this after a simplicial subdivision on the target.

   Let \(T\) be a ppt on \(\widetilde{\mathbf{a}}\).
   Then \(|T|\) is a maximal simplex in \(K(\widetilde{\mathbf{a}})\).
   Let \(S_1, \dots, S_m\) be the ppts on \(\mathbf{a}\) that resolve \(\degen(T)\).
   By \cref{lem:subdivision}, \(\kappa(|S_1|), \dots, \kappa(|S_m|)\) give a simplicial subdivision of \(|T|\).

   We take \(\Delta = \kappa(S_i)\).
   By \cref{lem:subdivision}, \(v \in \PS(\mathbf{a})\) satisfies \(\kappa(v) \in \Delta\) if and only if \(v \in S_i\).
   Therefore, we get \(\kappa^{-1}(\Delta) = |S_i|\).
   Since \(\kappa(|S_i|)\) is full-dimensional, the map \(\kappa \colon |S_i| \to \kappa(|S_i|)\) is a linear isomorphism.

   To treat \(K^{*}\), recall that by chasing the boundary parallel curve in the commutative diagram \eqref{eq:kappa-preserves-support}, we get that \(\kappa\) sends \(|K^{*}(\mathbf{a})|\) to \(|K^{*}(\widetilde{\mathbf{a}})|\).
   The same consideration shows that \(|K^{*}(\mathbf{a})|\) is in fact the pre-image of \(|K^{*}(\widetilde{\mathbf{a}})|\).
   So the isomorphism on \(K\) restricts to the isomorphism on \(K^{*}\).
 \end{proof}

 Having treated elementary deformations, we now treat arbitrary motions in configuration space.
 \begin{theorem}\label{thm:pl-homeomorphisms}
  Let \(\Gamma \colon [0,1] \to \BigConf_{n+1}\) be a continuous map with \(\mathbf{a} = \Gamma(0)\) and \(\mathbf{b} = \Gamma(1)\).
  \begin{enumerate}
  \item \(\Gamma\) induces piecewise linear isomorphisms \(\kappa_{\Gamma} \colon |K(\mathbf{a})| \to |K(\mathbf{b})|\) and \(\kappa_{\Gamma} \colon |K^{*}(\mathbf{a})|\ \to |K^{*}(\mathbf{b})|\).
  \item Let \(\deform_{\Gamma} \colon M(\mathbf{a}) \to M(\mathbf{b})\) be the bijection on multi-curves obtained by deforming along \(\Gamma\).
    Then the following diagram commutes
\[
  \begin{tikzcd}[column sep=3em]
    M(\mathbf{a})\ar{d}{\supp}\ar{r}{\deform_{\Gamma}}[swap]{\sim} & M(\mathbf{b}) \ar{d}{\supp} \\
    {|K(\mathbf{a})|}\ar{r}{\kappa_{\Gamma}}[swap]{\sim} & {|K(\mathbf{b})|}.
  \end{tikzcd}
\]
\item If \(\Gamma'\) is a path homotopic to \(\Gamma\) by a homotopy that fixes the end-points, then \(\Gamma\) and \(\Gamma'\) induce the same piecewise linear isomorphism \(|K(\mathbf{a})| \to |K(\mathbf{b})|\).
  \end{enumerate}
\end{theorem}
See \cref{fig:wall-crossing} for a picture of the simplest wall-crossing isomorphisms.
\begin{proof}
  First let us prove all assertions for a path \(\Gamma\) that can be written as a concatenation of paths that do not change relative positions, or are elementary degenerations/deformations (we call such paths ``easy paths'').
  Then the first assertion follows by composing the isomorphisms given by \cref{prop:pl-iso} and their inverses.
  The second assertion is a result of the commutativity we have already observed in diagram~\eqref{eq:kappa-preserves-support}.
  For the third assertion, suppose that \(\Gamma'\) is another easy path homotopic to \(\Gamma\).
  Then the bijection \(\deform_{\Gamma'} \colon M(\mathbf{b}) \to M(\mathbf{a})\) is equal to \(\deform_{\Gamma}\), the one induced by \(\Gamma\).
  We also have a map \(\kappa_{\Gamma'}\colon |K(\mathbf{a})| \to |K(\mathbf{b})|\).
By the commutativity of the diagram in the second assertion, the map \(\kappa_{\Gamma'}\) takes the same values as the map \(\kappa_{\Gamma}\) on the image of \(M(\mathbf{b})\).
  Since the image of \(M(\mathbf{b})\) is dense (\cref{prop:densecurves}), we conclude that \(\kappa_{\Gamma'} = \kappa_{\Gamma}\).

  Now suppose that \(\Gamma\) is an arbitrary path.
  Then observe that it is homotopic (fixing endpoints) to an easy path \(\Gamma'\).
  Set \(\kappa_{\Gamma}\) to be \(\kappa_{\Gamma'}\).
  Then the first and second assertions are automatic by the previous arguments.
  The third assertion follows from the previous argument that all homotopic easy paths induce the same map \(\kappa_{\Gamma'}\), and thus the choice of \(\Gamma'\) does not matter.
\end{proof}
The definitions of the complexes \(K(\mathbf{a})\), \(K^{\ast}(\mathbf{a})\), and their geometric realisations only depend on the unordered point configuration underlying \(\mathbf{a}\).
In other words, these simplicial complexes are well-defined on \(\BigConf_{n+1}/S_{n+1}\).
Recall that the fundamental group of \(\BigConf_{n+1}/S_{n+1}\) is the \((n+1)\)-strand braid group \(B_{n+1}\).
Thus we obtain the following corollary.
\begin{corollary}\label{cor:braid-pl-action}
We get an action of the braid group \(B_{n+1} = \pi_1(\BigConf_{n+1}/S_{n+1}, \mathbf{a})\) on \(|K(\mathbf{a})|\) and \(|K^{*}(\mathbf{a})|\) by piecewise linear isomorphisms.
\end{corollary}
\begin{proof}
  Fix a base point \(\mathbf{a} \in \BigConf_{n+1}\), and its image (also called \(\mathbf{a}\)) in \(\BigConf_{n+1}/S_{n+1}\).
  Consider a path in \(\BigConf_{n+1}/S_{n+1}\) from \(\mathbf{a}\) to \(\mathbf{a}\).
  This is equivalent to considering a path in \(\BigConf_{n+1}\) from \(\mathbf{a}\) to another point \(\mathbf{a}'\) in its \(S_{n+1}\) orbit.
  Since \(K(\mathbf{a}') = K(\mathbf{a})\) and \(K^{*}(\mathbf{a}') = K^{*}(\mathbf{a})\),~\cref{thm:pl-homeomorphisms} implies that this path induces piecewise-linear isomorphisms \(|K(\mathbf{a})| \xrightarrow{\cong} |K(\mathbf{a})|\) and \(|K^{\ast}(\mathbf{a})| \xrightarrow{\cong} |K^{\ast}(\mathbf{a})|\).
  By the third assertion in the theorem, the isomorphism above is independent of the homotopy class of the path relative to its endpoints.
  Thus the proof is complete.
\end{proof}
We can deduce the global topology of \(|K(\mathbf{a})|\) and \(|K^{\ast}(\mathbf{a})|\) by using~\cref{thm:pl-homeomorphisms} to reduce to, e.g., the case of convex point configurations.
This case has been extensively studied in the literature, particularly in the context of cluster algebras.
The simplicial complexes are more accessible in this case, because they become the simplicial complexes associated to triangulations of convex polygons.
We obtain the following corollary.
\begin{corollary}\label{cor:ball}
  Let \(\mathbf{a}\) be a configuration of \((n+1)\) points in \(\mathbf{C}\).
  Then \(|K(\mathbf{a})|\) is homeomorphic to the \((2n-3)\)-ball.
  The subcomplex \(|K^{*}(\mathbf{a})|\) is the boundary of \(|K(\mathbf{a})|\), and hence homeomorphic to the \((2n-4)\)-sphere.
\end{corollary}
\begin{proof}
  By \cref{thm:pl-homeomorphisms}, we may take \(\mathbf{a}\) to be any configuration.
  Let us consider a convex configuration of \((n+1)\) points.
  That is, one where the points of \(\mathbf{a}\) form the vertices of a convex polygon.
  In this case, pointed pseudo-triangulations are simply triangulations.
  The subcomplex of \(K(\mathbf{a})\) whose elements only contain interior diagonals of the polygon is called the \emph{cluster complex of type \(A\)}, and its study was initiated in~\cite{fom.zel:03}.
  The case of an \((n+1)\)-gon corresponds to the cluster complex of type \(A_{n-2}\).
  It is proved in~\cite[Corollary 1.11]{fom.zel:03} that the geometric realisation of the cluster complex of type \(A_{n-2}\) is a sphere of dimension \(n-3\).
  
  It is not hard to see that the complex \(K(\mathbf{a})\) is an iterated cone over the cluster complex, in which we take the cone by adding one external edge of the polygon at a time.
  Since there are \(n\) external edges, we see that \(|K(\mathbf{a})|\) is a ball of dimension \((2n-3)\), and its boundary is a sphere of dimension \((2n-4)\), as desired.

  Finally, it is also not hard to see in the convex case that \(|K^{\ast}(\mathbf{a})|\) is precisely the boundary of \(|K(\mathbf{a})|\), as follows.
  The maximal simplices of \(K(\mathbf{a})\) have size \((2n-2)\).
  Consider a simplex \(S'\) of size \((2n-3)\), formed by deleting an edge \(e\) from a maximal simplex \(S\).
  If the edge \(e\) is an interior diagonal, then there is exactly one other interior diagonal that crosses \(e\), which is compatible with all elements of \(S'\).
  Thus \(S'\) is contained in exactly two maximal simplices, and \(|K(S')|\) is not on the boundary of \(|K(\mathbf{a})|\).
  On the other hand, if \(e\) is an external edge, then there is no other maximal simplex that contains \(S'\), and thus \(|K(S')|\) is on the boundary of \(|K(\mathbf{a})|\).
  The proof is complete.
\end{proof}

\section{Comparison of the ppt complex with projective measured foliations}\label{sec:pmf}
Fix a configuration \(\mathbf{a}\).
Let \(\overline D\) be a closed disk containing all points of \(\mathbf{a}\) and let \(D\) be the punctured disk obtained by deleting from \(\overline D\) small open disks around the marked points.
Let \(\MF(\mathbf{a})\) be the space of measured foliations on \(D\) up to Whitehead equivalence, as defined in \cite[\S~11.1]{fat.lau.poe:12}.
Let \(\PMF(\mathbf{a})\) be \(\MF(\mathbf{a})\) modulo scaling.
Then \(\PMF(\mathbf{a})\) is piecewise-linearly homeomorphic to the \((2n-3)\)-sphere.

Recall that \(M^{*}(\mathbf{a})\) is the set of multi-curves on \(\mathbf{a}\) that do not contain boundary parallel components.
We have natural maps \(M^{*}(\mathbf{a}) \to |K^{*}(\mathbf{a})|\) and \(M^{*}(\mathbf{a}) \to \PMF(\mathbf{a})\).
\begin{theorem}\label{prop:pmf}
  We have a piecewise-linear isomorphism
  \[ |K^{*}(\mathbf{a})| \to \PMF(\mathbf{a})\]
  compatible with the maps from \(M^{*}(\mathbf{a})\).
\end{theorem}

The proof relies on the piecewise-linearity of intersection and occurrence functions, which we now define.

Let \(\delta\) be either a simple closed curve on \(D\) or an arc on \(D\) whose end-points lie on the boundary.
Given a simple closed curve \(\gamma\), let \(I(\gamma, \delta)\) be the geometric intersection number of \(\gamma\) and \(\delta\), defined as the minimum of the size of \(\gamma \cap \delta\) over all representatives in the corresponding isotopy classes.
We extend the definition to arcs on \(\mathbf{a}\) as defined in \cref{sec:configurations} as follows.
Given an arc \(\gamma\) joining two distinct marked points \(\mathbf{a}\), we consider the simple closed curve \(2 \cdot \gamma\) and let \[I(\gamma, \delta) = \frac{1}{2}\cdot I(2\gamma, \delta).\]
Note that, for any such arc \(\gamma\), we have \(I(\gamma, \gamma) = 0\).
We extend \(I(-, \delta)\) to multi-curves by linearity.

\begin{proposition}\label{prop:pl-intersections}
  Let \(\delta\) be either a simple closed curve on \(D\) or an arc on \(D\) whose end-points lie on the boundary.
  The function \(I(-,\delta)\) induces a piecewise-linear function on the unprojectivised complex \(\widetilde K^{*}(\mathbf{a})\).
  That is, there is a continuous piecewise-linear function
  \[I_{\delta} \colon \widetilde K^{*}(\mathbf{a}) \to \mathbf{R}_{\geq 0}\]
  such that for every multi-curve \(\gamma \in M^{*}(\mathbf{a})\), we have
  \[ I(\gamma, \delta) = I_{\delta}(\supp \gamma).\]
\end{proposition}
\begin{proof}
  Consider a path in configuration space from \(\mathbf{a}\) to \(\mathbf{b}\).
  By~\cref{thm:pl-homeomorphisms}, this path induces a piecewise-linear homeomorphism \(\phi \colon \widetilde K^{*}(\mathbf{a}) \to \widetilde K^{*}(\mathbf{b})\).
  Deforming along this path also induces a bijection \(\psi \colon M^{*}(\mathbf{a}) \to M^{*}(\mathbf{b})\).
  So, the statement for \(\psi(\delta)\) on \(\mathbf{b}\) is equivalent to the statement for \(\delta\) on \(\mathbf{a}\).
  By choosing \(\psi\) correctly, we may bring \(\mathbf{a}\) to a configuration of \((n+1)\) collinear points and \(\delta\) to one of the following curves/arcs: the segment joining the first two marked points (red), a segment joining the first marked point to the boundary (blue), a segment joining two points on the boundary with \(k\) points on the left and \((n+1-k)\) points on the right (orange), a simple closed curve enclosing the first \(k\) marked points (green); see the following diagram.
  For each of these, the intersection number is a piecewise-linear combination of the occurrence functions associated to the pseudo-straight segments on \(\mathbf{a}\) (we omit the calculation).
  This combination defines the required continuous piecewise-linear function on \(\widetilde K^{*}(\mathbf{a})\). 
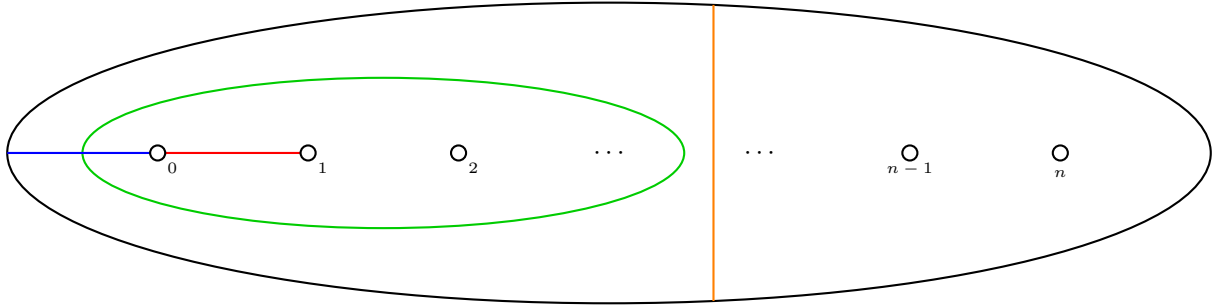
  \[
    \begin{tikzpicture}
      \begin{scope}[xscale=4]
        \draw[thick] (0,0) circle (2);
        \coordinate (a) at (80:2);
        \coordinate (b) at (280:2);
        \draw[thick, green!80!black] (-0.75,0) circle (1);
      \end{scope}
      \draw[thick, red] (-6,0) -- (-4,0);
      \draw[thick, blue] (-6,0) -- (-8,0);
      \draw[thick, orange] (a) -- (b);
                    \draw[thick, fill=white] (-6, 0) circle (0.1) node[below right]{\tiny \(0\)} (-4,0) circle (0.1) node[below right]{\tiny \(1\)} (-2, 0) circle (0.1) node[below right]{\tiny \(2\)} (0,0) node {\(\cdots\)} (2,0) node{\(\cdots\)} (4,0) circle (0.1) node[below]{\tiny \(n-1\)} (6,0) circle (0.1) node[below=0.1]{\tiny \(n\)};
    \end{tikzpicture}
  \]
\end{proof}

Fix a pseudo-straight segment \(s\) on \(\mathbf{a}\).
Given a simple closed curve \(\gamma\) on \(D\), let \(\operatorname{occ}(\gamma, s)\) be the coefficient of \(s\) in \(\supp \gamma \).
In other words, \(\operatorname{occ}(\gamma, s)\) is the number of occurencess of \(s\) in the curve obtained by pulling \(\gamma\) tight.
\begin{proposition}\label{prop:pl-occ}
  Let \(s\) be a pseudo-straight segment on \(\mathbf{a}\).
  The function \(\occ(-,s)\) induces a piecewise-linear function on \(\MF(\mathbf{a})\).
  That is, there is a continuous piecewise-linear function
  \[ \occ_s \colon \MF(\mathbf{a}) \to \mathbf{R}_{\geq 0}\]
  such that for every simple closed curve \(\gamma\), we have
  \[ \occ(\gamma, s) = \occ_s(\gamma).\]
\end{proposition}
\begin{proof}
  Consider a path in configuration space from \(\mathbf{a}\) to \(\mathbf{b}\).
  We identify \(M^{*}(\mathbf{a})\) with \(M^{*}(\mathbf{b})\) via the bijection given by deforming along the path. 
  
  \cref{thm:pl-homeomorphisms} gives a piecewise-linear isomorphism \(\widetilde K^{*}(\mathbf{a})\) to \(\widetilde K^{*}(\mathbf{b})\) compatible with the supports of multi-curves.
  In other words, the occurrence functions for pseudo-straight segments in \(\mathbf{b}\) are continuous piecewise-linear functions of the occurrence functions for pseudo-straight segments in \(\mathbf{a}\).
  Likewise, we also have a piecewise-linear isomorphism \(\MF(\mathbf{a}) \to \MF(\mathbf{b})\).
  So it suffices to prove the proposition for a single configuration.
    
  We choose \(\mathbf{a}\) to be the configuration of \(n+1\) points on a line.
  Let \(\alpha_1, \alpha_2, \dots, \alpha_{2n-2}\) be the dotted arcs and \(\beta_1, \dots, \beta_n\) the dashed arcs shown below.
  
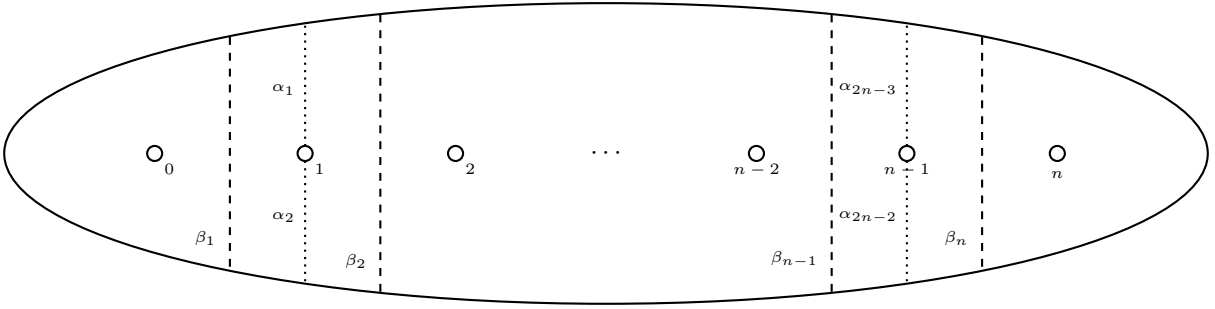
\[
    \begin{tikzpicture}
      \begin{scope}[xscale=4]
        \draw[thick] (0,0) circle (2);
      \end{scope}
      
      \draw[thick, dashed]
      (-5, 1.55) -- (-5,-1.55) node [above left=0.2 and 0.05]{\tiny \(\beta_{1}\)}
      (-3, 1.85) -- (-3,-1.85) node [above left=0.2 and 0.05]{\tiny \(\beta_{2}\)}
      (3, 1.85) -- (3,-1.85) node [above left=0.2 and 0.05]{\tiny \(\beta_{n-1}\)}
      (5, 1.55) -- (5,-1.55) node [above left=0.2 and 0.05]{\tiny \(\beta_{n}\)};
      \draw[thick, dotted]
      (-4,1.7) -- node[left] {\tiny \(\alpha_{1}\)} (-4,0)
      (-4,0) -- node[left] {\tiny \(\alpha_{2}\)} (-4,-1.7)
      (4,1.7) -- node[left] {\tiny \(\alpha_{2n-3}\)} (4,0)
      (4,0) -- node[left] {\tiny \(\alpha_{2n-2}\)} (4,-1.7);
              \draw[thick, fill=white] (-6, 0) circle (0.1) node[below right]{\tiny \(0\)} (-4,0) circle (0.1) node[below right]{\tiny \(1\)} (-2, 0) circle (0.1) node[below right]{\tiny \(2\)} (0,0) node {\(\cdots\)} (2,0) circle (0.1) node[below]{\tiny \(n-2\)} (4,0) circle (0.1) node[below]{\tiny \(n-1\)} (6,0) circle (0.1) node[below=0.1]{\tiny \(n\)};
    \end{tikzpicture}
  \]
  Pairing with these arcs gives continuous piecewise-linear functions \(a_1, \dots, a_{2n-2}\) and \(b_1, \dots, b_n\) on \(\PMF(\mathbf{a})\).
  The occurrence function of a pseudo-straight segment on \(\mathbf{a}\) can be written as continuous piecewise-linear combination of \(a_i\) and \(b_j\) (we omit the details of this calculation).
\end{proof}

We now have the tools to prove~\cref{prop:pmf}.
\begin{proof}[Proof of \cref{prop:pmf}]
  We explain the maps in both directions.
  A point of \(\PMF(\mathbf{a})\) is determined by the vector of intersection numbers with all curves (up to scaling).
  The map
  \[ |K^{*}(\mathbf{a})| \to \PMF(\mathbf{a})\]
  is defined by
  \[ x \mapsto [I_{\delta}(x)],\]
  where \(I_{\delta}\) is the piecewise-linear function from~\cref{prop:pl-intersections}.

  The map
  \[ \PMF(\mathbf{a}) \to |K^{*}(\mathbf{a})|\]
  is defined by
  \[ x \mapsto [\occ_{s}(x) \mid s \in \PS(\mathbf{a})],\]
  where \(\occ_s\) is the piecewise-linear function from~\cref{prop:pl-occ}.

  To check that these maps are mutual inverses, we may restrict to (multi)-curves, which form a dense subset of \(|K^{*}(\mathbf{a})|\) and \(\PMF(\mathbf{a})\).
  On curves, the fact that these are inverse maps follows from their definitions.
\end{proof}

\section{Categorical setup and background}\label{sec:categorical-setup}
This section contains some general categorical background we need for the remainder of the paper.
In~\cref{sec:stability-conditions}, we briefly discuss Bridgeland stability conditions on general triangulated categories, first introduced in~\cite{bri:07}.

In~\cref{sec:the-ks-category} we consider 2-Calabi-Yau (2-CY) categories associated to undirected graphs, or more generally, Coxeter diagrams.
For our applications, we will focus on the categories of type \(A\), versions of which have been extensively studied over a number of years~\cite{kho.sei:02,tho:06,sei.tho:01,bap.deo.lic:20,bap.deo.lic:22,bap.bec.lic:23}.
We will also consider the categories associated to associated to rank two Coxeter systems, which have been studied, e.g., in~\cite{hen.lic:24,hen:24,del.hen.lic:23,hen:22}.

In~\cref{sec:arcs-to-objects}, we recall a result of Khovanov--Seidel~\cite{kho.sei:02} for the 2-CY category of type \(A\).
This result associates a spherical object of this category (up to shift) to every arc on a linear point configuration.
We also generalise this procedure for point configurations lying in a fundamental domain for the symmetric group action on configurations.

\subsection{Bridgeland stability conditions}\label{sec:stability-conditions}
We assume familiarity with the notion of a Bridgeland stability condition on a triangulated category, which was introduced and developed in~\cite{bri:07}.
We only recall the notation and key results.
Let \(\mathcal{C}\) be a triangulated category with Grothendieck group \(K_{0}(\mathcal{C})\).
A stability condition \(\tau\) on \(\mathcal{C}\) consists of 
\begin{enumerate}
\item a \emph{central charge} \(Z_{\tau}\), which is a group homomorphism \(K_0(\mathcal{C}) \to \mathbf{C}\); and
\item a \emph{slicing} \(\mathcal{P}_{\tau}\), which is a collection of full abelian subcategories \(\mathcal{P}_{\tau}(\phi)\) of \(\mathcal{C}\), index by \(\phi \in \mathbf{R}\).
\end{enumerate}
The slicing and the central charge are required to satisfy several compatibility conditions, which we omit.
For each \(\phi \in \mathbf{R}\), the objects of \(\mathcal{P}_{\tau}(\phi)\) are called \emph{semistable} of phase \(\phi\).
An object of \(\mathcal{P}_{\tau}(\phi)\) is said to be \emph{stable} of phase \(\phi\) if it is a simple object of the category \(\mathcal{P}_{\tau}(\phi)\).
Every object \(X\) in \(\mathcal{C}\) admits a unique filtration whose factors are \(\tau\)-semistable objects of strictly decreasing phase.
This filtration is called the \emph{Harder--Narasimhan (HN) filtration}. 
For any interval \(I \subset \mathbf{R}\), we let \(\mathcal{P}_{\tau}(I)\) be the full subcategory of \(\mathcal{C}\) consisting of objects whose HN factors have phase in \(I\).
For any \(\phi \in \mathbf{R}\), the subcategories \(\mathcal{P}_{\tau}([\phi,\phi+1))\) and \(\mathcal{P}_{\tau}((\phi, \phi+1])\) are abelian, and are moreover hearts of bounded \(t\)-structures on \(\mathcal{C}\).
In particular, \(\heart = P([0,1))\) is the heart of a bounded \(t\)-structure on \(\mathcal{C}\) and \(Z_{\tau}\) sends non-zero-objects of \(\heart\) to the upper half plane \(\mathbf{H}\).

Conversely, \cite[Proposition 5.3]{bri:07} says that a stability condition \(\mathcal{C}\) is uniquely specified by
\begin{enumerate}
\item the heart \(\heart\) of a bounded \(t\)-structure on \(\mathcal{C}\); and 
\item a central charge \(Z \colon K_0(\heart) \to \mathbf{C}\) that sends each object of \(\heart\) to \(\mathbf{H}\),
\end{enumerate}
subject to the condition that every object of \(\heart\) has an HN filtration.
The above data gives a stability condition \(\tau\) with \(P_{\tau}([0,1)) = \heart\) and \(Z_{\tau} = Z\).

Let \(\Stab(\mathcal{C})\) be the set of Bridgeland stability conditions on \(\mathcal{C}\).
By~\cite[Corollary 1.3]{bri:07}, every connected component of \(\Stab(\mathcal{C})\) has the structure of a complex manifold such that the forgetful map \(\tau \mapsto Z_{\tau}\)  is a local isomorphism to a linear subspace of \(\Hom(K_0(\mathcal{C}), \mathbf{C})\).

\subsection{2-Calabi--Yau categories associated to Coxeter systems}\label{sec:the-ks-category}
We recall the definitions of 2-CY categories associated to simply laced Coxeter systems.
These categories have many equivalent definitions, especially in the case of finite (ADE) types.
In these types, they can be defined using minimal resolutions of Kleinian surface singularities or using pre-projective algebras (see, e.g.~\cite[\S 1]{bri:09}).
We give a formulation using zig-zag algebras.

Fix a field \(\k\) of characteristic zero.
Simply laced Coxeter systems can be described in terms of unlabelled undirected graphs.
To this end, let \(\Gamma\) be an undirected graph that has no self-loops or multiple edges.
Let \(\Gamma^{\dbl}\) be its doubled quiver, which has two directed edges \((i,j)\) and \((j,i)\) for every undirected edges \(\{i,j\}\) in \(\Gamma\).
Let \(k\Gamma^{\dbl}\) denote the path algebra of \(\Gamma^{\dbl}\), graded by path length.
For concreteness, choose an ordering \((1, \ldots,n)\) of the vertices of \(\Gamma\), and set the \emph{sign} of a directed edge \((i,j)\) in \(\Gamma^{\dbl}\) to be \[s_{ij} = \begin{cases}
  1,&i < j;\\
  -1,&i > j.
\end{cases}
\]
\begin{definition}
  The zig-zag algebra \(A(\Gamma)\) of \(\Gamma\) is the quotient of \(k\Gamma^{\dbl}\) by the following relations.
  \begin{enumerate}
  \item All length-3 paths are set to zero.
  \item All length-2 paths whose source and target are different are set to zero.
  \item For a vertex \(i\) and two vertices \(j, k\) that are connected to \(i\) by edges of \(\Gamma\), set
    \[s_{ij}(i,j)(j,i) = s_{ik}(i,k)(k,i).\]
  \end{enumerate}
\end{definition}
\begin{remark}
  The definition above is of the \emph{signed zig-zag algebra}, considered in~\cite[Definition 6.4]{bap.deo.lic:20}, as opposed to the un-signed zig-zag algebra that has appeared earlier in the literature (see, e.g.~\cite{hue.kho:01}).
  The signed zig-zag algebra gives better categorical properties, as explained in~\cite[Section 6]{bap.deo.lic:20}, and different choices of signs give isomorphic algebras~\cite[Proposition 6.5]{bap.deo.lic:20}.
  However, the signed and the un-signed zig-zag algebras are isomorphic in type \(A_n\), which is the main example for this paper.
\end{remark}
The length-zero paths \((i)\) in \(A(\Gamma)\) are primitive idempotents, and give rise to indecomposable projective modules \(P_i = A(\Gamma)(i)\).
Regard \(A(\Gamma)\) as a differential graded (dg) algebra where all differentials are zero.
Let \(K(\dgmod A(\Gamma))\) be the category of finite-dimensional dg modules over the dg algebra \(A(\Gamma)\).
Set \(\mathcal{C}_{\Gamma}\) to be the smallest full and strict triangulated subcategory of \(K(\dgmod A(\Gamma))\) containing the objects \(P_i\).
As explained earlier, the category \(\mathcal{C}_n\) is the category \(\mathcal{C}(\Gamma)\) in the case where \(\Gamma\) is a graph of type \(A_n\).

The category \(\mathcal{C}_n\) is a \(\k\)-linear, triangulated, and (strongly) 2-Calabi--Yau category.
The extension closure of the objects \(P_i\) in \(\mathcal{C}_n\) is an abelian category \(\heart_{\std}\), which is the heart of a bounded \(t\)-structure on \(\mathcal{C}_n\).
We call \(\heart_{\std}\) the \emph{standard heart} of \(\mathcal{C}_n\), and the induced bounded \(t\)-structure the \emph{standard \(t\)-structure} on \(\mathcal{C}_n\).
The objects \(P_i\) are simple objects of \(\heart_{\std}\).
\begin{remark}\label{rem:general-coxeter-category}
  The 2-CY category associated to a more general Coxeter system is defined in~\cite[Section 4.7]{hen.lic:24}, and has been studied further in~\cite{del.hen.lic:23,hen:22}.
  In this case, the correct analogue of the zigzag algebra is a particular algebra object in a fusion category associated to the Coxeter system.
  The associated 2-CY category is again defined as a subcategory of the homotopy category of dg modules over this zigzag algebra object.
  In~\cref{sec:rank-2-category} we give a few more details about the case we are most interested in, namely the case of rank two Coxeter systems.
\end{remark}

\subsection{The Khovanov--Seidel reconstruction}\label{sec:arcs-to-objects}
Let \(\mathcal{C}_n\) be the 2-CY category associated to a graph of type \(A_n\).
We recall a construction due to Khovanov--Seidel~\cite{kho.sei:02} that associates a spherical object of \(\mathcal{C}_n\) to an arc on a configuration of \(n+1\) points on a line.

Let \(\mathbf{a} = (a_0, \ldots, a_n)\) be a configuration such that all points lie on a single line.
For each \(i \in \{1,\ldots,n\}\) consider the line \(\ell_i\) perpendicular to the segment \(a_{i-1}a_i\) and passing through its midpoint.
We always draw arcs on \(\mathbf{a}\) such that they have minimal intersections with the lines \(\ell_i\).
Let \(\alpha\) be an arc.
We can reconstruct a dg module \(P_{\alpha}\) (up to triangulated shift) from this arc, as explicitly described below.
In fact, the segment \(a_ia_{i+1}\) is associated to the object \(P_i\) (up to triangulated shift), and the construction is compatible with the action of the braid group as proved in~\cite{kho.sei:02}.
So we can uniquely determine the object (up to shift) associated to any arc.

To describe the reconstruction explicitly, orient the arc \(\alpha\), and follow it from its start to its end point.
Suppose that \(\alpha\) successively intersects the lines \((\ell_{i_1}, \ldots, \ell_{i_m})\).
As a vector space, the dg module \(P_{\alpha}\) is a direct sum of objects \(P_{i_j}\langle g_j \rangle\), as \(j\) ranges over \(1\) to \(m\).
The integers \(g_j\) are grading shifts, which are determined inductively as described below.
Simultaneously, we inductively describe the differential on \(P_{\alpha}\), whose terms always map between \(P_{i_j}\langle g_j \rangle\) and \(P_{i_{j+1}}\langle g_{j+1} \rangle\) for some \(j\).
\begin{enumerate}
\item As the base case, set \(g_1 = 0\). (Note that this is an arbitrary choice, and varying it shifts the whole resulting object \(P_{\alpha}\).)
\item For the induction step, assume that we have already described \(g_{i_j}\) as well as the component of the differential between \(P_{i_{j-1}}\langle g_{j-1} \rangle\) and \(P_{i_j}\langle g_j \rangle\).
Consider the portion of \(\alpha_{i_j}\) of \(\alpha\) that starts from its intersection with \(\ell_{i_j}\) and ends at its intersection with \(\ell_{i_{j+1}}\).
  This portion is contained in the region either before the first line \(\ell_1\), between two successive lines, or beyond the last line \(\ell_n\).
  Each such region contains a single point \(a_m\) of the configuration, and we can detect whether the portion \(\alpha_{i_j}\) travels around the point \(a_m\) clockwise or anti-clockwise in this region.
  \begin{enumerate}
  \item If \(\alpha_{i_j}\) travels clockwise around \(a_m\), the corresponding differential maps \(P_{i_j}\langle g_j \rangle\) to \(P_{i_{j+1}}\langle g_{j+1}  \rangle\) by the unique non-trivial map between them.
  \item If \(\alpha_{i_j}\) travels anti-clockwise around \(a_m\), the corresponding differential maps \(P_{i_{j+1}}\langle g_{j+1} \rangle\) to \(P_{i_j}\langle g_j  \rangle\) by the unique non-trivial map between them.
  \end{enumerate}
  In either case, the grading shift \(g_{j+1}\) is  uniquely determined by \(g_j\).
\end{enumerate}

Now let us consider more general configurations.
Recall that \(S_{n+1}\) acts on \(\BigConf_{n+1}\) by permuting the coordinates.
Denote by \(\mathbf{H}\) the half-closed upper-half plane in \(\mathbf{C}\), defined as
\[\mathbf{H} = \{z \in \mathbf{C} \mid z \neq 0\text{ and }\arg(z) \in [0, \pi)\}.\]
Let \(\BigFDomain_{n+1} \subset \BigConf_{n+1}\) be the subset of configurations \((a_0, \ldots, a_n)\) where \(a_i - a_{i-1} \in \mathbf{H}\) for \(1 \leq i \leq n\).
Then \(\BigFDomain_{n+1}\) is a fundamental domain for the action of \(S_{n+1}\).
Fix \(\mathbf{a} = (a_0, a_1, \ldots, a_n) \in \BigFDomain_{n+1}\).

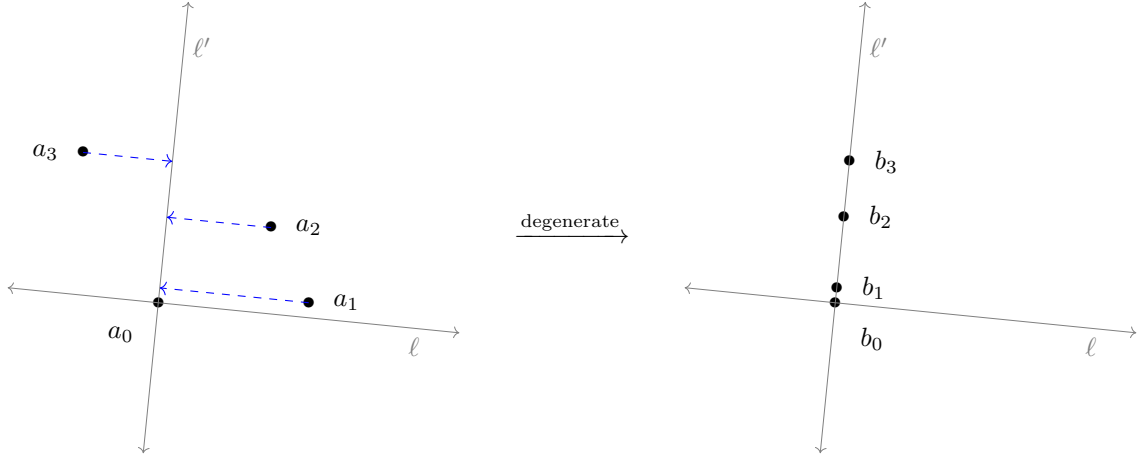
\begin{figure}[ht]
  \centering
  \begin{tikzpicture}
    \node[outer sep=0mm] (p0) at (0,0) {\(\bullet\)};
    \node[below left=0mm of p0] {\(a_0\)};
    \node[outer sep=0mm] (p1) at (2,0) {\(\bullet\)};
    \node[right=0mm of p1] {\(a_1\)};
    \node[outer sep=0mm] (p2) at (1.5,1) {\(\bullet\)};
    \node[right=0mm of p2] {\(a_2\)};        
    \node[outer sep=0mm] (p3) at (-1,2) {\(\bullet\)};
    \node[left=0mm of p3] {\(a_3\)};
    \coordinate (B) at (0.4,4);
    \coordinate (A) at (-0.2,-2);
    \draw[very thin, gray, <->] (-2,0.2) -- node[below, pos=0.9]{\(\ell\)} (4,-0.4);
    \draw[very thin, gray, <->] (A) -- node[right, pos=0.9]{\(\ell'\)} (B);
    \draw[blue, dashed,->] (p1.center) -- ($(A)!(p1)!(B)$);
    \draw[blue, dashed,->] (p2.center) -- ($(A)!(p2)!(B)$);
    \draw[blue, dashed,->] (p3.center) -- ($(A)!(p3)!(B)$);

    \node at (5.5,1) {\(\xrightarrow{\text{degenerate}}\)};

    \begin{scope}[xshift=9cm]
    \coordinate (A) at (-0.2,-2);
    \coordinate (B) at (0.4,4);
    \coordinate (p0) at (0,0);
    \coordinate (p1) at (2,0);
    \coordinate (p2) at (1.5,1);
    \coordinate (p3) at (-1,2);
    
    \node (q0) at (p0) {\(\bullet\)};
    \node (q1) at ($(A)!(p1)!(B)$) {\(\bullet\)};
    \node (q2) at ($(A)!(p2)!(B)$) {\(\bullet\)};
    \node (q3) at ($(A)!(p3)!(B)$) {\(\bullet\)};
    \node[below right=0mm of q0] {\(b_0\)};
    \node[right=0mm of q1] {\(b_1\)};
    \node[right=0mm of q2] {\(b_2\)};
    \node[right=0mm of q3] {\(b_3\)};            
    \draw[very thin, gray, <->] (-2,0.2) -- node[below, pos=0.9]{\(\ell\)} (4,-0.4);
    \draw[very thin, gray, <->] (A) -- node[right, pos=0.9]{\(\ell'\)} (B);
    \end{scope}
  \end{tikzpicture}
  \caption{A degeneration of the configuration \(\mathbf{a}\) to the linear configuration \(\mathbf{b}\).}
  \label{fig:shearing}
\end{figure}

Consider an almost-horizontal line \(\ell\) through the point \(a_0\) with slope \(-\epsilon\) for \(0 \ll \epsilon < 1\), and let \(\ell'\) be the almost-vertical line perpendicular to \(\ell\) passing through \(a_0\).
Choose \(\epsilon\) sufficiently small so that none of the lines \(a_ia_j\) are parallel to the line \(\ell\).
Consider a path \(\mathbf{a}(t)\) in \(\BigConf_{n+1}\) that starts at \(\mathbf{a}(0) = \mathbf{a}\), and translates every point \(a_i\) along the direction of \(\ell\) so that all points of \(\mathbf{b} = \mathbf{a}(1)\) lie on \(\ell'\).
Such a degeneration is illustrated in~\cref{fig:shearing}.

Such a degeneration gives a bijection from arcs on \(\mathbf{a}\) to arcs on \(\mathbf{b}\).
Since the points of \(\mathbf{b}\) all lie on a line, any arc on \(\mathbf{b}\) corresponds to an object of \(\mathcal{C}_{n}\) up to shift via the Khovanov--Seidel construction.
We also set this to be the object associated to the corresponding arc on \(\mathbf{a}\).
\begin{proposition}\label{prop:generalised-ks}
  Let \(\mathbf{a} = (a_0, a_1, \ldots, a_n) \in \BigFDomain_{n+1}\).
  Then the construction described above gives a map from arcs on \(\mathbf{a}\) to spherical objects of \(\mathcal{C}_n\) up to triangulated shift, independent of the choice of \(\epsilon\).
  Furthermore, the straight-line segment from \(a_i\) to \(a_{i+1}\) corresponds to the object \(P_i\) up to triangulated shift.
\end{proposition}
\begin{proof}
  The bijection from arcs on \(\mathbf{a}\) to arcs on \(\mathbf{b}\) is independent of \(\epsilon\), provided it is sufficiently small.
It is immediate from the construction that the straight line segment from \(a_i\) to \(a_{i+1}\) corresponds to the object \(P_i\) up to triangulated shift.
\end{proof}

\section{Stability conditions on \(\mathcal{C}_n\) via point configurations}\label{sec:stab-via-point-configurations}
Let $\mathcal{C}_n$ be the 2-CY triangulated category associated to the $A_n$ graph.
\cref{sec:arcs-to-objects} gave a correspondence from arcs on a point configuration to spherical objects in \(\mathcal{C}_n\).
In this section, we enhance this correspondence to Bridgeland stability conditions and HN filtrations.
First, we recall the explicit description of the space of stability conditions on \(\mathcal{C}_n\) as the universal cover of configuration space.
Next we identify the semistable objects of a given stability condition as special kinds of curves.
The main result of this section is~\cref{thm:hn-multiplicity-count}, which describes the HN filtration of any spherical object in terms of curves on point configurations.

\begin{remark}\label{rem:thomas-work}
  Some of the ideas of this section originate in the older paper~\cite{tho:06} by Thomas.
  The paper~\cite{tho:06} considers configurations of points that form the vertices of a convex polygon.
  In such cases,~\cite[Section 5]{tho:06} constructs corresponding standard stability conditions, and shows that the semistable (= stable) spherical objects correspond to straight-line segments.
  Our results extend this further to HN filtrations and arbitrary point configurations.
\end{remark}

\subsection{Standard stability conditions}
Recall that the Grothendieck group \(K(\mathcal{C}_n)\) is the lattice of type \(A_n\).
So \(\Hom(K(\mathcal{C}_n), \mathcal{C})\) is the complexified lattice of type \(A_n\), which we identify with \(\mathbf{C}^n = \mathbf{C}^{n+1}/ \mathbf{C}\). 
Recall that \(\BigConf_{n+1}\) is the complement of the root hyperplanes in the same space.
By the discussion around Theorem 6.4 of~\cite{tho:06}, the central charge of any stability condition on \(\mathcal{C}_n\) takes non-zero values on the roots.
As a result, we have a map
\[
  \begin{split}
    \pi \colon \Stab(\mathcal{C}_n) &\to \BigConf_{n+1} \\
    \tau &\mapsto Z_{\tau}.
  \end{split}
\]
Explicitly, we can write the map as follows.
Recall the indecomposable objects $P_1,P_2, \ldots, P_n$ of $\mathcal{C}_n$.
Then
\begin{equation}\label{eq:covering-map}
  \pi(\tau) = \left(0,Z_{\tau}(P_1), \dots, \sum_{j \leq i}Z_{\tau}(P_j), \ldots, \sum_{j \leq n}Z_{\tau}(P_j)\right).
\end{equation}
The following proposition is a combination of~\cite[Theorem 6.4]{tho:06} and the fact that \(\Stab(\mathcal{C}_n)\) is connected (see, e.g.,\cite{ish.ued.ueh:10,ish.ueh:05,ada.miz.yan:19}).
\begin{proposition}\label{prop:universal-cover}
  The map $\pi \colon \Stab(\mathcal{C}_n) \to \BigConf_{n+1}$ is a universal covering map.
\end{proposition}
The covering map \(\pi\) intertwines the \(B_{n+1}\) action on \(\Stab(\mathcal{C}_n)\) and the \(S_{n+1}\) action on \(\BigConf_{n+1}\).
Denote by \(\mathbf{H}\) the half-closed upper-half plane in \(\mathbf{C}\), defined as
\[\mathbf{H} = \{z \in \mathbf{C} \mid z \neq 0\text{ and }\arg(z) \in [0, \pi)\}.\]
Let \(\BigFDomain_{n+1} \in \BigConf_{n+1}\) be the fundamental domain for the action of \(S_{n+1}\) that consists of configurations \((a_0, \ldots, a_n)\) where \(a_i - a_{i-1} \in \mathbf{H}\) for \(1 \leq i \leq n\).
We restrict our attention to a particularly accessible subset of \(\Stab(\mathcal{C}_n)\), constructed via the following proposition.
\begin{proposition}\label{prop:lifting}
  Let \(\mathbf{a} = (a_0,\ldots,a_n)\) be a point of \(\BigFDomain_{n+1} \subset \BigConf_{n+1}\).
  Then there is a unique stability condition \(\tau \in \pi^{-1}(\mathbf{a})\) satisfying the following conditions.
  \begin{enumerate}
  \item The objects $P_1$, $P_2$, \ldots, $P_n$ are \(\tau\)-stable.
  \item The \(\tau\)-phases of the objects $P_1$, $P_2$, \ldots, $P_n$ lie in the interval $[0,1)$.
  \end{enumerate}
Conversely, for any stability condition \(\tau\) satisfying these two properties, its image \(\pi(\tau) \in \BigConf_{n+1}\) under the covering map~\eqref{eq:covering-map} lies in \(\BigFDomain_{n+1}\).
\end{proposition}
We say that a stability condition is \emph{standard} if it satisfies the conditions of~\cref{prop:lifting}.
\begin{proof}
  The category $\mathcal{C}_n$ has a \emph{standard heart} \(\heart_{\std}\), namely the extension closure of the objects $P_1$, \ldots, $P_n$.
  The Grothendieck group \(K_0(\mathcal{C}_n) = K_0(\heart_{\std})\) is the free abelian group generated by the classes of the \(P_i\).
  We specify \(\tau\) by writing down the obvious central charge \(Z_{\tau}\) on \(K_0(\heart_{\std})\):
  \[Z_{\tau} \colon P_i\mapsto (a_i - a_{i-1}).\]
  Since \(\heart_{\std}\) has finite length, this central charge function automatically enjoys the HN property, and therefore (uniquely) specifies the stability condition \(\tau\).

  The converse is clear.
\end{proof}
\subsection{Semistable objects and Harder--Narasimhan filtrations via curves}
Fix a standard stability condition \(\tau\) on \(\mathcal{C}_n\).
In this case, \(\mathcal{P}_{\tau}([0,1))\) is exactly \(\heart_{\std}\).
Recall that \(\heart_{\std}\) is the extension closure of the objects \(P_1,\ldots,P_n\) in \(\mathcal{C}_n\), and also these objects are simple in \(\heart_{\std}\).
This section has two aims.
\begin{enumerate}
\item Describe all ``\(\tau\)-stable curves'' and ``\(\tau\)-semistable curves'' in \(\heart_{\std}\).
  These are curves that correspond to \(\tau\)-stable and \(\tau\)-semistable spherical objects lying in \(\heart_{\std}\).
\item Given a spherical object corresponding to a curve, describe its HN filtration with respect to \(\tau\) in terms of semistable curves.
\end{enumerate}
Towards the first aim, we write down a categorical algorithm to reconstruct all stable and semistable spherical objects in \(\heart_{\std}\), and then apply it to curves.
The algorithm we write here is a slight generalisation of~\cite[Proposition 4.2]{bap.deo.lic:22}.
The statement in~\cite{bap.deo.lic:22} assumes that \(\tau\) is generic, in the sense that it has no indecomposable semi-stable spherical objects that are not already stable.
Our generalisation omits this assumption, but is proved almost identically.

Let us recall some notation.
Let \(\alpha\) be an arbitrary positive root of type \(A_n\).
Regarding \(\alpha\) as an element of the root lattice and hence of the Grothendieck group \(K_0(\mathcal{C}_n)\), the central charge \(Z_{\tau}(\alpha)\) lives in the upper half plane \(\mathbf{H}\).
Consider two subsets of \(\mathbf{H}\), defined as
\begin{align*}
  \mathbf{H}^{\alpha,\tau}_+ &= \{z \in \mathbf{H} \mid \arg z \geq \arg Z_{\tau}(\alpha)\};\\
  \mathbf{H}^{\alpha,\tau}_- &= \{z \in \mathbf{H} \mid \arg z \leq \arg Z_{\tau}(\alpha)\},
\end{align*}
where \(\arg\) takes values in \([0,\pi)\).

Consider an expression
\[\alpha = s_{i_1} s_{i_2} \ldots s_{i_{n-1}}(\alpha_{i_n}),\]
where \(\alpha_{i_n}\) is a simple root.
Associate to this expression its \emph{root sequence} \(R = (R_{i_1}, \ldots, R_{i_n})\), where
\[R_{i_k} = s_{i_1} \cdots s_{i_{k-1}}(\alpha_{i_k}).\]
Note that \(R_{i_n} = \alpha\).

\begin{proposition}\label{prop:more-general-algorithm}
  Consider the setup above for a fixed expression
  \[\alpha = s_{i_1} s_{i_2} \ldots s_{i_{n-1}}(\alpha_{i_n}).\]
  Consider a sign vector \(\underline{\epsilon} = (\epsilon_1, \ldots, \epsilon_{n-1})\), where \(\epsilon_k \in \{+1,-1\}\) for each \(k\).
  Let \(X_{\underline{\epsilon}}\) be the object
  \[X_{\underline{\epsilon}} = \sigma_{i_1}^{\epsilon_1} \cdots \sigma_{i_{n-1}}^{\epsilon_{n-1}}(P_{i_n}),\]
  where \(\sigma_k\) is the spherical twist in the object \(P_k\).
  \begin{enumerate}
  \item The object \(X_{\underline{\epsilon}}\) is \(\tau\)-semistable if and only if \(Z_{\tau}(R_{i_k}) \subset \mathbf{H}^{\alpha,\tau}_{\epsilon_k}\) for each \(k\).
    Explicitly, we require
    \begin{align*}
      \arg Z_{\tau}(R_{i_k}) &\geq \arg Z_{\tau}(\alpha) \text{ if }\epsilon_k = +1,\\
      \arg Z_{\tau}(R_{i_k}) &\leq \arg Z_{\tau}(\alpha) \text{ if }\epsilon_k = -1.
    \end{align*}
  \item Furthermore, the object \(X_{\underline{\epsilon}}\) is \(\tau\)-stable if and only the inequalities above are strict for each \(k\).
  \end{enumerate}
  Finally, objects constructed as above are the only possible \(\tau\)-semistable (resp. \(\tau\)-stable) spherical objects of class \(\alpha\) in the Grothendieck group.
\end{proposition}
\begin{proof}
  The proof of~\cite[Proposition 4.2]{bap.deo.lic:22} extends easily to this more general situation to yield proofs of parts (1) and (2).

  Let us prove that the objects constructed as above are the only possible \(\tau\)-semistable (resp.~\(\tau\)-stable) spherical objects of class \(\alpha\) in the Grothendieck group.  In type \(A_n\) (or more generally, any finite type), all spherical objects in the category \(\mathcal{C}_n\) lie in the braid group orbit of one of the basic objects \(P_i\)~\cite{ish.ued.ueh:10,ish.ueh:05,bap.deo.lic:22}.
  Thus an arbitrary semistable object \(X\) of class \(\alpha\) can be written as \(\underline{\beta}(P_i)\) for some braid word \(\beta\).
  Passing to the Grothendieck group, we obtain an expression of the form \(\alpha = \underline{w}(\alpha_i)\) for some sequence of simple reflections \(\underline{w}\).
  Now we apply parts (1) and (2) to see that the signs in the word \(\underline{\beta}\) must satisfy the constraints mentioned in parts (1) and (2).
\end{proof}

Armed with the semistability and stability criterion above, we return to describing indecomposable semistable and stable spherical objects of \(\heart_{\std}\).
We will make use of the following observation.
\begin{lemma}\label{lem:reverse-braid}\mbox{}
  \begin{enumerate}
  \item Let \(\alpha_j\) be a simple root of type \(A_n\).
    Then for any \(i < j\), we have an equality of roots:
    \[s_is_{i+1} \cdots s_{j-1}(\alpha_j) = s_j s_{j-1} \cdots s_{i+1}(\alpha_i).\]
  \item Let \(P_j\) be the \(j\)th standard simple object of \(\heart_{\std}\).
    For any \(i < j\) and a  set of indices \(\epsilon_i, \ldots, \epsilon_{j-1} \in \{+1,-1\}\), we have an equality of objects of \(\mathcal{C}_n\):
    \[\sigma_i^{\epsilon_i} \cdots \sigma_{j-1}^{\epsilon_{j-1}}(P_j) = \sigma^{-\epsilon_{j-1}}_j \sigma^{-\epsilon_{j-2}}_{j-1} \cdots \sigma^{-\epsilon_i}_{i+1}(P_i).\]
  \item Consider a configuration \((a_0, \ldots, a_n)\) in the fundamental domain, and let \(p_j\) be the straight line arc joining \(a_{j-1}\) to \(a_j\).
    Then we have an equality of arcs:
    \[\sigma_i^{\epsilon_i} \cdots \sigma_{j-1}^{\epsilon_{j-1}}(p_j) = \sigma^{-\epsilon_{j-1}}_j \sigma^{-\epsilon_{j-2}}_{j-1} \cdots \sigma^{-\epsilon_i}_{i+1}(p_i).\]
  \end{enumerate}
\end{lemma}
\begin{proof}
  Let us prove the second assertion; the other two are analogous.
  We prove the result by induction.
  As the base case, it is easy to check directly that
  \[\sigma_{j-1}^{\epsilon_{j-1}}(P_j) = \sigma_j^{-\epsilon_{j-1}}(P_{j-1}).\]
  
  For the induction step, consider the left hand side
  \(\sigma_i^{\epsilon_i} \cdots \sigma_{j-1}^{\epsilon_{j-1}}(P_j)\)
  and rewrite it as
  \[\sigma_i^{\epsilon_i} \cdots \sigma_{j-1}^{\epsilon_{j-1}}(P_j) = \left(\sigma_j^{-\epsilon_{j-1}} \sigma_j^{\epsilon_{j-1}}\right)\sigma_i^{\epsilon_i} \cdots \sigma_{j-1}^{\epsilon_{j-1}}(P_j).\]
  Noting that \(\sigma_j^{\epsilon_{j-1}}\) commutes with all factors but the last, we obtain
  \[\sigma_j^{-\epsilon_{j-1}}\left(\sigma_i^{\epsilon_i} \cdots\sigma_{j-2}^{\epsilon_{j-2}}\right)\sigma_j^{\epsilon_{j-1}} \sigma_{j-1}^{\epsilon_{j-1}}(P_j).\]
  It is clear from the base case calculation that \(\sigma_j^{\epsilon_{j-1}} \sigma_{j-1}^{\epsilon_{j-1}}(P_j) = P_{j-1}\).
  Thus we see that
  \[\sigma_i^{\epsilon_i} \cdots \sigma_{j-1}^{\epsilon_{j-1}}(P_j) = {\sigma_j^{-\epsilon_{j-1}}}\sigma_i^{\epsilon_i}\cdots \sigma_{j-2}^{\epsilon_{j-2}}(P_{j-1}).\]
  By applying the inductive hypothesis to the smaller expression
  \(\sigma_i^{\epsilon_i}\cdots \sigma_{j-2}^{\epsilon_{j-2}}(P_{j-1})\), we obtain the desired result.

\end{proof}
\begin{proposition}\label{prop:curve-stables-and-semistables}
  Let \(\mathbf{a} = (a_0, \ldots, a_n) \in \BigFDomain_{n+1}\).
  Let \(\tau_{\mathbf{a}}\) be the stability condition corresponding to \(\mathbf{a}\) described in~\cref{prop:lifting}.
 An indecomposable spherical object is \emph{semistable} (resp.~\emph{stable}) with respect to \(\tau_{\mathbf{a}}\) if and only if can be reconstructed via~\cref{prop:generalised-ks} from a pseudo-straight (resp.~straight) segment on the configuration \(\mathbf{a}\).
\end{proposition}
\begin{proof}
  Fix \(\mathbf{a}\) as in the proposition.
  For brevity, we write \(\tau_{\mathbf{a}}\) as \(\tau\) in the remainder of this proof.
  Let \(X\) be an indecomposable semistable spherical object in \(\heart_{\std}\).
  Its class \(\alpha = [X]\) in the Grothendieck group is a positive root of the \(A_n\) root system.
  Thus we have \(\alpha = \alpha_i + \alpha_{i+1} + \cdots + \alpha_j\) for some \(i \leq j\).
  Also recall that for each \(i\), the object \(P_i\) is stable in \(\heart_{\std}\) of class \(\alpha_i\).

  First suppose that \(i = j\); that is, \([X] = \alpha_i = [P_i]\).
  In this case, let us show that \(X = P_i\).
  The hom-pairing between \(X\) and \(P_i\) can be computed as 
  \[\langle X, P_i \rangle = \sum_i(-1)^j \dim\Hom(X,P_i[j]) = \langle \alpha_i, \alpha_i \rangle = 2.\]
  Since \(X\) and \(P_i\) are both in \(\heart_{\std}\), there are no negative-degree maps between them in either direction.
  The 2-Calabi--Yau property of \(\mathcal{C}_n\) now implies that the only possible non-zero morphisms from \(X\) to \(P_i\) lie in degrees \(0\), \(1\), and \(2\).
  The fact that \(\langle X, P_i \rangle = 2\) implies that either \(\Hom^0(X,P_i) \neq 0\) or \(\Hom^2(X,P_i) = \Hom^0(P_i,X) \neq 0\).
  Since \(P_i\) is a simple object of the standard heart \(\heart_{\std}\), we see that \(P_i\) is either a quotient or a sub-object of \(X\), which means that there exists a short exact sequence in \(\heart_{\std}\) that contains \(P_i\) and \(X\) as adjacent objects.
  It is either of the form
  \[0 \to P_i \to X \to Y \to 0\]
  or of the form
  \[0 \to Y \to X \to P_i \to 0\]
  for some \(Y \in \heart_{\std}\).
  However, since \([X] = [P_i]\) in the Grothendieck group, we have \([Y] = 0\) in the Grothendieck group, which can only happen if \(Y = 0\).
  Therefore \(X = P_i\).

  Otherwise, \([X] = \alpha_i + \cdots + \alpha_j\) for \(i < j\).
  In this case, we have a (minimal) expression
  \[
    \alpha = s_is_{i+1}\cdots s_{j-1}(\alpha_j),
  \]
  where \(s_k \in S_{n+1}\) is the root reflection corresponding to the simple root \(\alpha_k\).
  By~\cref{prop:more-general-algorithm}, we know that there is a sign vector \((\epsilon_i, \ldots, \epsilon_{j-1})\) such that
  \[X = \sigma_i^{\epsilon_i} \cdots \sigma_{j-1}^{\epsilon_{j-1}}(P_j),\]
  and \(Z_{\tau}(R_k) \subset \mathbf{H}^{\alpha,\tau}_{\epsilon_k}\) for each \(k\).
  Let \(p_\ell\) be the straight line segment joining \(a_{\ell-1}\) to \(a_\ell\).
  The object \(X\) is associated to the arc
  \[x = \sigma_i^{\epsilon_i} \cdots \sigma_{j-1}^{\epsilon_{j-1}}(p_j) = \sigma_j^{-\epsilon_{j-1}}\sigma_{j-1}^{-\epsilon_{j-2}} \cdots \sigma_{i+1}^{-\epsilon_i}(p_i),\]
  where the second equality is a consquence of~\cref{lem:reverse-braid}.
  We claim that \(x\) is isotopic to a pseudo-straight segment joining \(a_{i-1}\) to \(a_j\).
  Using that \(x\) is obtained from \(p_i\) by successive twists in \(p_{i+1}, \dots, p_j\), we can describe it as follows.
  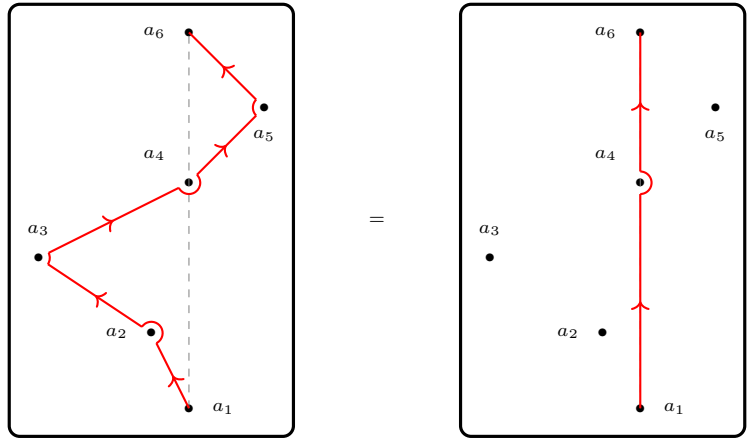
\begin{figure}[ht]
    \centering
    \begin{tikzpicture}[font=\scriptsize]
      \tikzset{
        config/.pic={
          \node[outer sep=0mm] (-p1) at (0,0) {\(\bullet\)};
          \node[right=0mm of -p1] {\(a_1\)};
          \node[outer sep=0mm] (-p2) at (-0.5,1) {\(\bullet\)};
          \node[left=0mm of -p2] {\(a_2\)};        
          \node[outer sep=0mm] (-p3) at (-2,2) {\(\bullet\)};
          \node[above=0mm of -p3] {\(a_3\)};                
          \node[outer sep=0mm] (-p4) at (0,3) {\(\bullet\)};
          \node[above left=0mm of -p4] {\(a_4\)};                        
          \node[outer sep=0mm] (-p5) at (1,4) {\(\bullet\)};
          \node[below=0mm of -p5] {\(a_5\)};                                
          \node[outer sep=0mm] (-p6) at (0,5) {\(\bullet\)};
          \node[left=0mm of -p6] {\(a_6\)};
          \draw[very thin, gray, dashed] (-p1.center) edge (-p6.center);
          \node[draw, solid, very thick, rounded corners, inner sep=2mm, fit=(-p1) (-p3) (-p5) (-p6)] {};
        }}

    \node[matrix] (f1) at (0,0) {\pic (c1) {config};\\};
    \node[matrix] (f2) at (6,0) {\pic (c2) {config};\\};
    \node at ($(f1)!0.5!(f2)$) {\(=\)};

\coordinate (r1) at ($(c1-p1)!1.5mm!(c1-p2)$);
    \coordinate (q2) at ($(c1-p2)!1.5mm!(c1-p1)$);
    \coordinate (r2) at ($(c1-p2)!1.5mm!(c1-p3)$);
    \coordinate (q3) at ($(c1-p3)!1.5mm!(c1-p2)$);
    \coordinate (r3) at ($(c1-p3)!1.5mm!(c1-p4)$);
    \coordinate (q4) at ($(c1-p4)!1.5mm!(c1-p3)$);
    \coordinate (r4) at ($(c1-p4)!1.5mm!(c1-p5)$);
    \coordinate (q5) at ($(c1-p5)!1.5mm!(c1-p4)$);
    \coordinate (r5) at ($(c1-p5)!1.5mm!(c1-p6)$);

\coordinate (x1) at ($(c2-p4)!1.5mm!(c2-p1)$);
    \coordinate (x2) at ($(c2-p4)!1.5mm!(c2-p6)$);    

      \path[thick, red, decoration={markings,mark=at position 0.5 with {\arrow{>}}}]
      (c1-p1.center) edge[postaction=decorate] (q2)
      (r2) edge[postaction=decorate] (q3)
      (r3) edge[postaction=decorate] (q4)
      (r4) edge[postaction=decorate] (q5)
      (r5) edge[postaction=decorate] (c1-p6.center)
      ;
      \pic [draw=red, thick, angle radius=1.5mm]{angle=c1-p1--c1-p2--c1-p3};
      \pic [draw=red, thick, angle radius=1.5mm]{angle=c1-p2--c1-p3--c1-p4};
      \pic [draw=red, thick, angle radius=1.5mm]{angle=c1-p3--c1-p4--c1-p5};
      \pic [draw=red, thick, angle radius=1.5mm]{angle=c1-p6--c1-p5--c1-p4};

      \path[thick, red, decoration={markings,mark=at position 0.5 with {\arrow{>}}}]
      (c2-p1.center) edge[postaction=decorate] (x1)
      (x2) edge[postaction=decorate] (c2-p6.center)
      ;
      \pic [draw=red, thick, angle radius=1.5mm]{angle=x1--c2-p4--x2};      
    \end{tikzpicture}
    \caption{Construction of the semistable curve \(\sigma_1\sigma_2\sigma_3\sigma_4^{-1}(p_5)= \sigma_5\sigma_4^{-1}\sigma^{-1}_3\sigma_2^{-1}(p_1)\).}
    \label{fig:semistable-as-ps-segment}
  \end{figure}
  Start with the curve \(y\) that is the concatenation of \(p_i, \dots, p_j\).
  This is not an arc as it passes through the marked points \(a_i, \dots, a_{j-1}\).
  For \(\ell \in \{i, i+1, \ldots, j-1\}\) in order, we modify \(y\) in a small neighbourhood of \(a_{\ell}\) by swerving it so that it avoids \(a_{\ell}\).
    Instead of following \(p_{\ell}\) all the way to \(a_{\ell}\), we stop close to it, take a circular detour until we reach \(p_{\ell+1}\), and then continue along \(p_{\ell+1}\).
    The detour is counter-clockwise if \(\epsilon_{\ell} = 1\) and clockwise if \(\epsilon_{\ell} = -1\) (see~\cref{fig:semistable-as-ps-segment}).
    Let \(D\) be the closed region in the plane bounded by the curve \(y\) and the segment \(a_ia_j\).
    Note that the following are equivalent:
    \begin{enumerate}
    \item the vector \(a_{\ell}-a_i\) in the upper half plane is strictly to the left (resp.~strictly to the right) of \(a_j-a_i\);
    \item a small counter-clockwise (resp. clockwise) detour around \(a_{\ell}\) lies in the region \(D\);
    \item \(Z_{\tau}(R_{\ell})\) lies in the open part of \(\mathbf{H}^{\alpha,\tau}_{\epsilon_{\ell}}\).
    \end{enumerate}
    Here, the third point follows from~\cref{lem:reverse-braid}, by observing again that
    \[R_{\ell} = s_is_{i+1} \cdots s_{\ell-1}(\alpha_{\ell}) = s_{\ell}s_{\ell-1}\cdots s_{i+1}(\alpha_i).\]
    In particular, if all \(\epsilon_{\ell}\) are chosen as in \cref{prop:more-general-algorithm}, then all detours at points outside the segment \(a_ia_j\) lie in \(R\).
    Therefore, the region of the plane bounded by the segment \(a_ia_{j}\) and \(x\) does not contain any points of the configuration outside the segment \(a_ia_j\).
  As a result, \(x\) can be isotoped to be arbitrary close to the segment \(a_ia_j\).
  That is, \(x\) is pseudo-straight.

  Conversely, any pseudo-straight segment from \(a_i\) to \(a_j\) has a representative \(x\) that is a modification of the curve \(y\) obtained by taking small circular detours around the intermediate points \(a_{\ell}\).
  Since the segment is pseudo-straight, the detours at all points outside of \(a_ia_j\) must lie in \(R\).
  Otherwise, the interior of the region bounded by \(x\) and the segment \(a_ia_{j}\) would contain points \(a_{\ell}\) not on the segment \(a_i a_j\), contradicting that \(x\) is isotopic to curves arbitrary close to the segment \(a_ia_j\).
  By the equivalence of the three conditions above and \cref{prop:more-general-algorithm}, the object \(X = \sigma_i^{\epsilon_i} \cdots \sigma_{j-1}^{\epsilon_{j-1}}(P_j) \) associated to \(x\) is semi-stable.

  Now note that \(X\) is \emph{stable} if and only if \(Z_{\tau}(R_{\ell})\) lies in the open part of \(\mathbf{H}^{\alpha,\tau}_{\epsilon_{\ell}}\) for every \(\ell\).
  This occurs if and only if the vector \(a_l - a_i\) is not parallel to \(a_j-a_i\) for any \(\ell\), which occurs if and only if there is no point \(a_{\ell}\) on the segment \(a_ia_j\).
  Finally, this occurs if and only if the semistable curve \(x\) we construct is a straight segment.
\end{proof}

Let us now describe the HN filtration of an arbitrary arc in terms of semistable (i.e. pseudo-straight) arcs.
We will do this by explicitly using the chain complexes produced by the Khovanov--Seidel reconstruction theorem, and arguing directly with properties of these chain complexes.

First, we write down a straightforward but useful observation, whose proof we omit.
It is convenient to think of a chain complex \(X^{\bullet}\) as a graded vector space with a differential \(\partial \colon X^{\bullet} \to X^{\bullet}\) that raises graded degree by \(1\).
If we have a decomposition of graded vector spaces \(X^{\bullet} = X_1^{\bullet} \oplus X_2^{\bullet}\), then \(\partial\) has components \(\partial_{ij} \colon X^{\bullet}_i \to X^{\bullet}_j\) for \(i,j \in \{1,2\}\).
\begin{lemma}\label{lem:subcomplex}
  Consider a complex \((X^{\bullet}, \partial)\).
  Suppose that \(X^{\bullet} = X_1^{\bullet} \oplus X_2^{\bullet}\) as graded vector spaces, and that \(\partial_{21} = 0\).
  Then \((X_1^{\bullet}, \partial_{11})\) and \((X_2^{\bullet}, \partial_{22})\) are themselves chain complexes, and we have  an exact sequence of chain complexes
  \[0 \to (X_{2}^{\bullet}, \partial_{22}) \to (X^{\bullet}, \partial) \to (X_1^{\bullet},\partial_{11}) \to 0.\]
\end{lemma}
Now consider a configuration \((a_0, \ldots, a_n) \in \FDomain_{n+1}\).
Recall that any arc joining two distinct points among $(a_0, a_1,\ldots, a_n)$ can be ``pulled tight'' around these points.
The arc then breaks up along its spine, as defined in~\cref{def:spine} and illustrated in~\cref{fig:spine}.
Each element of the spine is a pseudo-straight segment, and thus corresponds to a \(\tau\)-semistable object.
In fact, we can uniquely augment the spine to record the relative phase of each element in the spine.
\begin{definition}\label{def:augmented-spine}
  Let \(\mathbf{a} = (a_0, \dots, a_n)\) be a configuration and let \(\gamma\) be an oriented arc from \(a_i\) to \(a_j\), and let
  \((\gamma_1, \ldots, \gamma_k)\) be its spine.
  The \emph{phase-augmented spine} (or just \emph{augmented spine}) of \(\gamma\) is the tuple \(((\gamma_1, \varphi_1), \dots, (\gamma_k, \varphi_k))\) defined as follows.
  \begin{enumerate}
  \item \(\varphi_1 \in [0,2)\) is \((1/\pi)\) times the angle made by the underlying straight vector of \(\gamma_1\) with the horizontal.
  \item For every \(i > 0\), the quantity \(\varphi_i - \varphi_{i-1}\) is defined to be \((1/\pi)\) times the signed angle from the underlying straight vector of \(\gamma_{i-1}\) to the underlying straight vector of \(\gamma_i\).
    We assume that \(\varphi_i - \varphi_{i-1} \in [-1,1)\) for every \(i > 1\).
    Note also that \(\varphi_i \neq \varphi_{i-1}\) for any \(i\).
  \end{enumerate}
\end{definition}
Let \(s\) be a pseudo-straight segment whose associated Khovanov--Seidel complex in \(\heart_{\std}\) is \(X^{\bullet}\).
Given an \emph{phase-augmented} pseudo-straight segment \((s,\varphi)\); that is, an orientation on \(s\) together with an arbitrary real number \(\varphi \in \mathbf{R}\), associate to it the Khovanov--Seidel complex \(X^{\bullet}[n]\), where \(n = \lfloor \varphi \rfloor\).

With the convention above, the complex associated to \(\gamma\) and the complexes associated to the \((\gamma_i,\varphi_i)\) are compatible with each other and also with \(\tau\), in the following sense.
\begin{lemma}\label{lem:cutting-curve-complex}
  Let \(\gamma\) be an oriented arc with augmented spine \(((\gamma_1,\varphi_1), \ldots, (\gamma_r, \varphi_r))\).
  Let \((X^{\bullet}, \partial)\) be a choice of a Khovanov--Seidel complex associated to \(\gamma\).
  For each \(i\), let \((X^{\bullet}_i, \partial_i)\) be the Khovanov--Seidel complex associated to \((\gamma_i,\varphi_i)\).
  \begin{enumerate}
  \item For every \(i\), the \(\tau\)-semistable object \(X_i\) has phase \(\varphi_i\).
  \item\label{cuts}
    For some integer \(m\), we have a decomposition \(X^{\bullet}[m] = \bigoplus_i X_i^{\bullet}\) of graded vector spaces, satisfying the following properties.
    \begin{enumerate}
    \item\label{cut-pieces} For every \(i\), the component \(\partial \colon X^{\bullet}_i \to X^{\bullet}_i\) agrees with \(\partial_i\).
    \item\label{cut-distant} For any \(i,j\) such that \(|i - j| > 1\), the component \(\partial \colon X^{\bullet}_i \to X^{\bullet}_j\) is zero.
    \item\label{cut-neighbour} For any \(i,j\) such that \(|i - j| = 1\), the component \(\partial \colon X^{\bullet}_i \to X^{\bullet}_j\) is non-zero if and only if \(\varphi_i < \varphi_j\).
    \end{enumerate}
  \end{enumerate}
\end{lemma}
\begin{proof}
  This can be deduced by carefully examining the Khovanov--Seidel construction.
\end{proof}
The following theorem states that the indecomposable pieces of the HN factors of the Khovanov--Seidel complex associated to any arc are just the spherical objects associated to its augmented spine.

Fix a configuration \((a_0,\ldots, a_n) \in \FDomain_{n+1}\).
Let \(\gamma\) be an arc on this configuration, and orient \(\gamma\) so that if it has start and end points \(a_i\) and \(a_j\) respectively, then \(i < j\).
Let \(X^{\bullet}\) be a complex that corresponds to \(\gamma\) under the Khovanov--Seidel construction.
\begin{theorem}\label{thm:hn-multiplicity-count}
  Let \(\gamma\) and \(X^{\bullet}\) be as above.
  Let \(((\gamma_1,\varphi_1), \ldots, (\gamma_k,\varphi_k))\) be the augmented spine of \(\gamma\).
  Let \(X_i^{\bullet}\) be the Khovanov--Seidel complex associated to the augmented pseudo-straight segment \((\gamma_i, \varphi_i)\).
  Then there is some integer \(m\) such that for any \(\varphi \in \mathbf{R}\), the HN factor of \(X^{\bullet}[m]\) of phase \(\varphi\) is precisely the direct sum of all \(X_i\) for which \(\varphi_i = \varphi\).
\end{theorem}

\begin{proof}
  We inductively construct the HN filtration of \(X^{\bullet}\).
  By replacing \(X^{\bullet}\) by a shift if necessary, assume that we have a direct sum decomposition of graded vector spaces \(X^{\bullet} = \bigoplus_{i = 1}^r X^{\bullet}_i\) satisfying the properties in~\cref{lem:cutting-curve-complex}\eqref{cuts}.

  Let \(S \subset \{1, \cdots, r\}\) be the set of \(i\) such that \(\varphi_i = \min(\varphi_1, \dots, \varphi_r)\).
  Consider an \(i \in S\).
  By \cref{def:augmented-spine}, we have \(\varphi_{i \pm 1} \neq \varphi_{i}\), and therefore \(\varphi_{i \pm 1} > \varphi_i\).
  In particular, \((i \pm 1) \notin S\).
  By~\cref{lem:cutting-curve-complex}\eqref{cut-neighbour}, the components of \(\partial\) taking \(X^{\bullet}_{i \pm 1}\) to \(X^{\bullet}_i\) are zero.
  By~\cref{lem:cutting-curve-complex}\eqref{cut-distant}, the components of \(\partial\) between \(X^{\bullet}_i\) and \(X^{\bullet}_j\) for \(|i - j| > 1\) are all zero.
  Now apply~\cref{lem:subcomplex} to the decomposition
  \[X^{\bullet} = \left(\bigoplus_{t \notin S} X^{\bullet}_t\right) \oplus \left(\bigoplus_{t \in S}X^{\bullet}_t\right).\]
  We obtain an exact sequence of chain complexes as follows.
  \begin{equation}\label{eq:peel-off-min-phase}
    0 \to \left (\bigoplus_{t \notin S} X^{\bullet}_t, \partial\right ) \to (X^{\bullet},\partial) \to \left(\bigoplus_{t \in S}X^{\bullet}_t, \partial \right) \to 0.
  \end{equation}
  Since the indices in \(S\) are separated by at least 2, \cref{lem:cutting-curve-complex}\eqref{cut-distant} implies that the components of \(\partial\) between distinct summands of \(\bigoplus_{t \in S} X^{\bullet}_t\) are zero.
  Therefore, the last term of the sequence is in fact a direct sum of complexes \(\bigoplus_{t \in S} (X^{\bullet}_t, \partial_t)\).
  Set \(A_1 = \bigoplus_{t \in S} (X^{\bullet}_t, \partial_t)\).
  Note that \(A_1\) is a direct sum of semi-stable sphericals of the same phase. 
  
  The exact sequence~\eqref{eq:peel-off-min-phase} gives rise to a distinguished triangle
  \[\left (\bigoplus_{t \notin S} X^{\bullet}_t, \partial\right ) \to (X^{\bullet},\partial) \to A_1 \xrightarrow{+1}\]
  in \(\mathcal{C}\).
  Now we inductively apply the same procedure to \(\left( \bigoplus_{t \notin S} X^{\bullet}, \partial \right)\).
  Namely, at the \(i\)th step, we peel off the portion of minimum phase and call it \(A_i\).
  After at most \(r\) steps, we obtain a filtration of \(X\) with factors \(A_s, \dots, A_1\) where
  \begin{enumerate}
  \item each \(A_j\) is a direct sum of semistable indecomposable objects \(X^{\bullet}_t\) of the same phase, say \(\varphi(A_j)\).
  \item for each \(j\), we have \(\varphi(A_j) > \varphi(A_{j-1})\).
  \end{enumerate}
  We have thus constructed the HN filtration of \(X^{\bullet}\), and it has the desired properties by construction.
\end{proof}
\begin{remark}
  Let \(X^{\bullet}\) be the spherical object associated to the curve \(\gamma\).
  Choose an orientation of \(\gamma\) and let \((\gamma_1, \dots, \gamma_n)\) be the spine of \(\gamma\) (see~\cref{def:spine}).
  \cref{thm:hn-multiplicity-count} says that the HN factors of \(X^{\bullet}\) are precisely the semi-stable spherical objects corresponding to \(\gamma_1, \dots, \gamma_n\) (up to shifts).
  The phase order of the HN factors, however, is not the natural geometric order of the \(\gamma_i\)'s.
\end{remark}

\begin{corollary}\label{cor:HN-pieces-spherical}
  Let \(X\) be a spherical object of \(\mathcal{C}_n\) and let \(\tau\) be a stability condition on \(\mathcal{C}_n\).
  Then the \(\tau\)-semistable HN factors of \(X\) are direct sums of semistable spherical objects.
\end{corollary}
\begin{proof}
  Any stability conditions is in the braid group orbit of a standard stability condition.
  So we may take \(\tau\) to be standard.
  Then the statement follows from~\cref{thm:hn-multiplicity-count}.
\end{proof}

Let \(\mathbf{a}\) be a configuration and \(\tau = \tau_{\mathbf{a}}\) the corresponding standard stability condition.
Recall that we have the ppt* complex \(K^{*}(\mathbf{a})\) from~\cref{def:pptcomplex} and the complex of semistable sphericals \(\Sigma_{\tau}\) from~\cref{def:tau-simplicial-complex}.
\begin{corollary}\label{cor:K=sigma}
  The bijection between pseudo-straight segments on \(\mathbf{a}\) and \(\tau\)-semistable spherical objects (up to shift) given by \cref{prop:curve-stables-and-semistables} induces an isomorphism between \(K^{*}(\mathbf{a})\) and \(\Sigma_{\tau}\).
\end{corollary}
\begin{proof}
  By~\cref{prop:traintrack,prop:densecurves}, ppt*s are precisely the maximal collections of pseudo-straight segments that appear when an arc is pulled tight.
  So the result follows from \cref{thm:hn-multiplicity-count}.
\end{proof}

\section{The simplicial complex for rank two Coxeter systems}\label{sec:rank-2}
In~\cref{sec:the-ks-category} we gave the construction of the 2-CY category associated to a simply-laced Coxeter system.
As mentioned in \cref{rem:general-coxeter-category}, the construction generalises to give a 2-CY category for every finite rank Coxeter system.
The aim of this section is to explain how to generalise the definition of the simplicial complex \(\Sigma_{\tau}\) and recover analogues of the piecewise-linear homeomorphisms of~\cref{sec:simplicial-complex} in the setting of rank two Coxeter systems.
\subsection{The 2-CY category of a rank two Coxeter system}\label{sec:rank-2-category}
Fix a rank two Coxeter system \(I_2(n)\), with the following Coxeter diagram.
\[
  \begin{tikzcd}
    \bullet \arrow[-]{r}{n} & \bullet
  \end{tikzcd}
\]
We sketch the definition of the associated 2-CY category \(\mathcal{C}_2(n)\) below; for more details see~\cite[\S4]{hen.lic:24} and citations therein.
However, the following points are sufficient to understand the simplicial complex \(\Sigma_{\tau}\).
\begin{itemize}
\item As a \(\mathbf{C}\)-linear triangulated category, the category \(\mathcal{C}_2(n)\) is equivalent to the 2-CY category \(\mathcal{C}_{n-1}\) of type \(A_{n-1}\), studied in~\cref{sec:stab-via-point-configurations}.
  The additional structure on \(\mathcal{C}_2(n)\) is an action of a fusion category called the Temperley--Lieb--Jones category.
\item As a result, it is natural to restrict our attention to fusion-equivariant stability conditions on \(\mathcal{C}_2(n)\).
\end{itemize}
Let \(\TLJ_n\) denote the Temperley--Lieb--Jones category at a $2n$-th root of unity.
$\TLJ_n$ is a semi-simple fusion category.
A \(\mathbf{Z}\)-graded enhancement of \(\TLJ_n\) has a particular algebra object \(A(I_2(n))\) called the \emph{zig-zag algebra of type \(I_2(n)\)}.
We may regard \(A(I_2(n))\) as a dg algebra object with trivial differential.
Now we consider the category of dg modules over this dg algebra object, and proceed as in simply laced type to define the associated 2-CY category \(\mathcal{C}_2(n)\).
This category naturally carries a (right) action of \(\TLJ_n\).

The monoidal structure on \(\TLJ_n\) induces a natural product on the Grothendieck group of \(\TLJ_n\).
Thus \(\TLJ_n\) has a Grothendieck ring \(K_0(\TLJ_n)\), which acts on the Grothendieck group \(K_0(\mathcal{C}_2(n))\) of \(\mathcal{C}_2(n)\).
\begin{remark}
  The construction of \(\mathcal{C}_2(n)\) is a generalisation of the simply-laced construction detailed in~\cref{sec:the-ks-category}.
  Indeed, the (usual) zig-zag algebra is an algebra object in the category of \(\mathbf{Z}\)-graded vector spaces, and the Grothendieck ring of the category of vector spaces (which is isomorphic to \(\mathbf{Z}\)) acts on the Grothendieck ring of the associated 2-CY category.
\end{remark}

\subsection{Fusion-equivariant stability conditions}
The space of stability conditions \(\Stab(\mathcal{C}_2(n))\) has a distinguished closed submanifold consisting of the stability conditions that are \(\TLJ_n\)-equivariant, whose meaning we now explain.

We say that a slicing is \(\TLJ_n\)-equivariant if each subcategory defined by the slicing is closed under the action of \(\TLJ_n\).
To describe \(\TLJ_n\)-equivariant central charges, we choose an action of \(K_0(\TLJ_n)\) on \(\mathbf{C}\), as follows.
There is an isomorphism
\[K_0(\TLJ_n) \cong \mathbf{Z}[x]/\Delta_{n-1}(x),\]
where \(\Delta_{n-1}(x)\) is the \((n-1)\)th normalised Chebyshev polynomial of the second kind.
The assignment \(x \mapsto 2 \cos(\pi/n)\) thus yields a ring homomorphism from this ring to \(\mathbf{R}\).
This homomorphism is known as the Frobenius--Perron dimension \(\FP\).
We let \(K_0(\TLJ_n)\) act on \(\mathbf{C}\) via the map \(\FP\).
Then we say that a central charge is \(\TLJ_n\)-equivariant it intertwines the \(K_0(\TLJ_n)\) actions on \(K_0(\mathcal{C}_2(n))\) and on \(\mathbf{C}\).

We say that a stability condition on \(\mathcal{C}_2(n)\) is \emph{fusion-equivariant} if its slicing and central charge are \(\TLJ_n\)-equivariant as described above.
We refer to~\cite{hen:22,hen.lic:24,del.hen.lic:23} for more details.

Let \(\mathbf{P}(\Stab(\mathcal{C}_2(n)))\) denote the quotient of the space of fusion-equivariant stability conditions by the action of \(\mathbf{C}\).
Let us recall the basic topological and combinatorial structure of the space \(\mathbf{P}(\Stab(\mathcal{C}_2(n)))\).
We summarise the main points in the theorem below, which is a compilation of results from~\cite{bap.deo.lic:20,hen:22}.
\begin{theorem}\label{thm:rank2}
Let \(n \geq 2\) be a positive integer. We have the following statements.
\begin{itemize}
    \item The space \(\mathbf{P}(\Stab(\mathcal{C}_2(n)))\) is homeomorphic to an open 2-dimensional disk.
    \item The disk \(\mathbf{P}(\Stab(\mathcal{C}_2(n)))\) is tiled by $n$-gons without their vertices, such that the edges of the \(n\)-gon consist of stability conditions that lie on a wall.
  \item The vertices of the \(n\)-gons lie on the boundary of the disk, and are in bijection with the orbits of the spherical objects under the action of \(\TLJ_n\) as well as the shift functor.
\end{itemize}
\end{theorem}
The theorem above is proven for $n=3$ in~\cite{bap.deo.lic:20}, where the relevant fusion category is just the category of vector spaces, and in~\cite{hen:22} for other values of $n$, where the action of \(\TLJ_n\) plays an important role.
The case $n=\infty$ is also worked out in~\cite{bap.deo.lic:20}: in this case, almost all statements of the the above theorem hold, except that $\infty$-gons that tile the disk all have an accumulation point on the boundary, which is a vertex corresponding to a semi-rigid (rather than spherical) object.

The theorem above has several nice consequences, illustrated in~\cref{fig:polygonal-tiling}.
In particular, recall that the vertices of the \(n\)-gons that tile \(\mathbf{P}(\Stab(\mathcal{C}_2(n)))\) are labelled by spherical objects of \(\mathcal{C}_2(n)\), modulo shifts and the action of \(\TLJ_n\).
\begin{proposition}\label{prop:stability-polygons}
We have the following properties.
\begin{itemize}
\item  If \(\tau\) is an off-the wall stability condition, then it lies in the interior of an \(n\)-gon.
  The \(\tau\)-stable objects are the vertices of this \(n\)-gon.
\item If \(\tau'\) is an on-the-wall stability condition, then it lies on an edge that borders two adjacent \(n\)-gons.
  The \(\tau'\)-semistable objects are the vertices of the \((2n-2)\)-gon formed as the union of these two \(n\)-gons.
  Among these, the \(\tau'\)-stable objects are the two vertices of the edge containing \(\tau'\).
\end{itemize}
\end{proposition}
In order to study the simplicial complex \(\Sigma_{\tau}\), we need to understand the HN factors of spherical objects.
The following results achieve this goal, and are proved using a tool called HN automata constructed in~\cite{bap.deo.lic:20,hen:22}.
Fix a fusion-equivariant stability condition on \(\mathcal{C}_2(n)\).
In what follows, we consider semistable objects up to the action of \(\TLJ_n\) and the shift functor.
\begin{theorem}\label{thm:two-hn-factors}
  In the setup above, any spherical object in \(\mathcal{C}_2(n)\) has at most two distinct HN factors.
\end{theorem}
\begin{theorem}\label{thm:polygonal-hn-decomposition}
  In the setup above, suppose that \(A\) and \(B\) are two distinct semistable spherical objects that appear together in the HN filtration of a spherical object \(X\).
  We have the following.
  \begin{enumerate}
  \item If the chosen stability condition \(\tau\) is off-the-wall, then \(A\) and \(B\) are connected by an edge of the \(n\)-gon containing \(\tau\) in its interior.
  \item If the chosen stability condition \(\tau'\) is on-the-wall, then \(A\) and \(B\) are connected by an edge of the \((2n-2)\)-gon containing \(\tau'\) in its interior.
  \item In each case, \(A\) and \(B\) are precisely the endpoints of the edge of the polygon that separates \(X\) and \(\tau\).
  \end{enumerate}
\end{theorem}
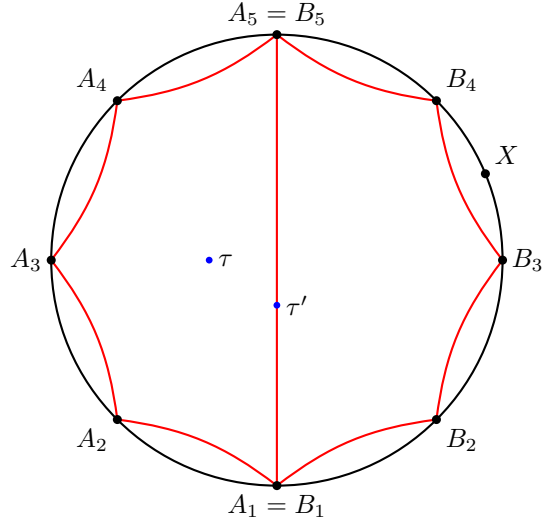
\begin{figure}[h]\label{fig:tiling}
  \centering
  \begin{tikzpicture}[scale=3]
\draw[thick, black] (0,0) circle (1);
    
\coordinate (A1) at (0,-1);           \coordinate (A2) at (-0.707,-0.707);  \coordinate (A3) at (-1,0);           \coordinate (A4) at (-0.707,0.707);   \coordinate (A5) at (0,1);            \coordinate (B4) at (0.707,0.707);    \coordinate (B3) at (1,0);            \coordinate (B2) at (0.707,-0.707);   

\coordinate (X) at (0.924,0.383);     

\draw[red, thick] (A1) to[bend right=15] (A2);
    \draw[red, thick] (A2) to[bend right=15] (A3);
    \draw[red, thick] (A3) to[bend right=15] (A4);
    \draw[red, thick] (A4) to[bend right=15] (A5);
    \draw[red, thick] (A5) to[bend right=15] (B4);
    \draw[red, thick] (B4) to[bend right=15] (B3);
    \draw[red, thick] (B3) to[bend right=15] (B2);
    \draw[red, thick] (B2) to[bend right=15] (A1);
    
\draw[red, thick] (A1) -- (A5);
    
\fill[black] (A1) circle (0.02);
    \fill[black] (A2) circle (0.02);
    \fill[black] (A3) circle (0.02);
    \fill[black] (A4) circle (0.02);
    \fill[black] (A5) circle (0.02);
    \fill[black] (B4) circle (0.02);
    \fill[black] (B3) circle (0.02);
    \fill[black] (B2) circle (0.02);
    \fill[black] (X) circle (0.02);
    
\node[below] at (A1) {$A_1 = B_1$};
    \node[below left] at (A2) {$A_2$};
    \node[left] at (A3) {$A_3$};
    \node[above left] at (A4) {$A_4$};
    \node[above] at (A5) {$A_5=B_5$};
    \node[right] at (B3) {$B_3$};
    \node[above right] at (B4) {$B_4$};
    \node[below right] at (B2) {$B_2$};
    \node[above right] at (X) {$X$};
    
\coordinate (tau) at (-0.3,0);        \coordinate (tauprime) at (0,-0.2);   

    \fill[blue] (tau) circle (0.015);
    \fill[blue] (tauprime) circle (0.015);
    
\node[right] at (tau) {$\tau$};
    \node[right] at (tauprime) {$\tau'$};
    
  \end{tikzpicture}
  \caption{Pentagonal tiling of \(\mathbf{P}(\Stab(\mathcal{C}_2(5)))\).
    \(\tau\) is an off-the-wall stability condition, and \(\tau'\) is an on-the-wall stability condition.
    The five objects \(A_i\) are \(\tau\)-stable.
    The objects \(A_i\) together with the objects \(B_i\) (altogether eight) are \(\tau'\)-semi-stable.
    Among these, the objects \(A_1 = B_1\) and \(A_5 = B_5\) are \(\tau'\)-stable.
    The \(\tau\)-HN filtration of the object \(X\) contains \(A_5\) and \(A_1\).
    The \(\tau'\)-HN filtration of \(X\) contains \(B_4\) and \(B_3\).}\label{fig:polygonal-tiling}
\end{figure}

\subsection{The simplicial complex \(\Sigma_\tau\) in rank 2}
The tiling of \(\mathbf{P}(\Stab(\mathcal{C}_2(n)))\) from Theorem \ref{thm:rank2}, together with~\cref{thm:two-hn-factors,thm:polygonal-hn-decomposition}, immediately gives a description of the simplicial complexes \(\Sigma_\tau\) for different stability conditions, as well as an intuitive picture of the piecewise-linear homeomorphisms between them.

By Theorem \ref{thm:rank2}, the spherical objects of \(\mathcal{C}_2(n)\) lie on the \(S^1\) boundary of \(\mathbf{P}(\Stab(\mathcal{C}_2(n)))\).  In particular, the semi-stable spherical objects of a fusion-equivariant stability condition \(\tau\) correspond to points on this \(S^1\).  Let \(P_{\tau}\) be the 1-dimensional simplicial complex 
whose vertices are these $\tau$-semistable spherical points, and whose edges are the arcs of the boundary \(S^1\) connecting adjacent points.  To avoid confusion: note that both the vertices and edges of \(P_{\tau}\) lie completely on the boundary; this is in contrast to the polygons which tile 
\(\mathbf{P}(\Stab(\mathcal{C}_2(n)))\), which have their vertices on the boundary and their edges in the interior.

\begin{theorem}\label{thm:simplicialrank2}
There is an isomorphism of simplicial complexes from \(\Sigma_{\tau}\) to \(P_{\tau}\).
\end{theorem}
\begin{proof}
  The vertices of \(\Sigma_{\tau}\) are the \(\tau\)-semistable spherical objects, by definition.
  By~\cref{prop:stability-polygons}, these are in bijection with the vertices of \(P_{\tau}\).
  The edges of \(\Sigma_{\tau}\) are pairs of vertices that can appear together as HN factors of a single spherical object.
  By~\cref{thm:polygonal-hn-decomposition}, these correspond precisely to the edges of \(P_{\tau}\).
  By~\cref{thm:two-hn-factors}, there are no higher-dimensional cells in \(\Sigma_{\tau}\).
\end{proof}
Let \(\tau\) be an off-the-wall stability condition, and let \(\tau'\) be an adjacent on-the-wall stability condition.  That is, if \(\tau\) lies in the interior of a tiling polygon, then \(\tau'\) lies on one of the edges of that polygon.
\begin{proposition}\label{prop:PLrank2}
In the setup above, there is a piecewise-linear homeomorphism from \(\Sigma_{\tau}\) to \(\Sigma_{\tau'}\), which is compatible with the HN supports of spherical objects.
\end{proposition}
\begin{proof}
  We describe the map from \(\Sigma_{\tau'}\) to \(\Sigma_{\tau}\).
  We write the vertices of \(\Sigma_\tau\) as $A_1,\hdots,A_n$, clockwise.  Assume \(\tau'\) lies on the wall joining \(A_n\) and \(A_1\).  We denote the vertices of 
  \(\Sigma_{\tau'}\) by \(A_1=B_1,A_2,\hdots,A_{n-1},A_n=B_n, B_{n-1},\hdots, B_2\) clockwise.
  The map sends a vertex of \(\Sigma_{\tau'}\) to its \(\tau\)-HN support.  We extend this linearly to edges.  Note that the map is the identity on the vertices that are common to both \(\Sigma_{\tau}\) and \(\Sigma_{\tau'}\); moreover, one can check that the \(\tau\)-HN support of \(B_i\) is 
  \[
    \frac{n-i}{\gcd(n-i,i-1)}A_1 + \frac{i-1}{\gcd(n-i,i-1)}A_n.  
  \]
  The only thing remaining to check is that this map is compatible with HN supports, meaning that the $\tau'$ support vector of a spherical object $X$ is taken to the $\tau'$ support of $X$; this is a consequence of the fact that the \(\tau\)-HN filtration of a spherical object $X$ refines the $\tau'$-HN filtration of $X$.
  This refinement fact is proven using the machinery of HN automata in \cite{bap.deo.lic:20,hen:22}.  Now is easy to check that this map is a piecewise-linear isomorphism.
\end{proof}
As a result of Theorem \ref{thm:simplicialrank2}, the simplicial complexes $\Sigma_\tau$, for different stability conditions $\tau$, may be thought of as different simplicial decompositions of a single circle, which is the boundary of the stability manifold.
For this reason, Theorem \ref{thm:simplicialrank2} and Proposition \ref{prop:PLrank2} are perhaps best understood visually.  For example, the image below continues the example from~\cref{fig:polygonal-tiling}, and draws each of the simplicial complexes \(\Sigma_\tau\) and \(\Sigma_{\tau'}\) as simplicial decompositions of the boundary circle.
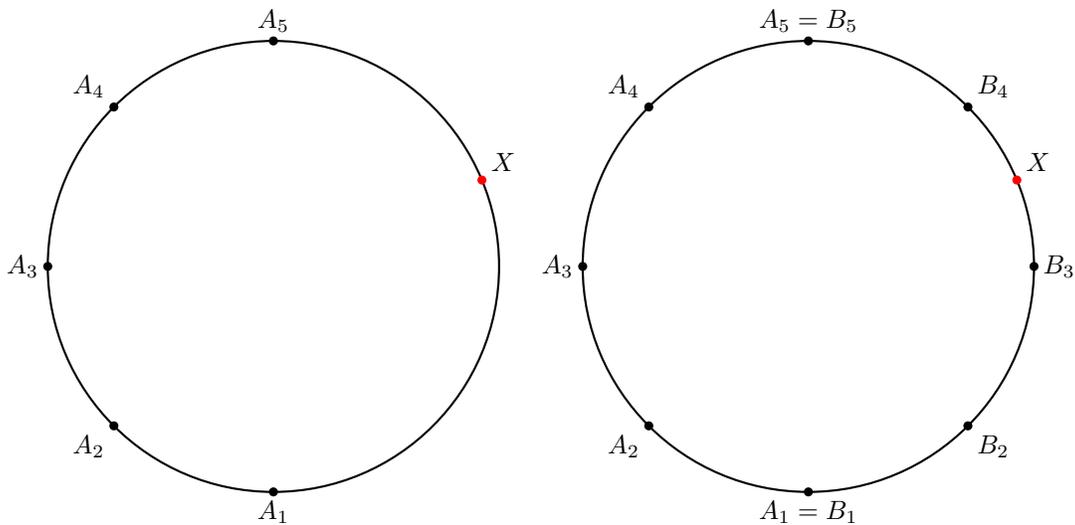
\begin{figure}[h]
  \centering
  \begin{tikzpicture}[scale=3]
\draw[thick, black] (0,0) circle (1);
    
\coordinate (A1) at (0,-1);           \coordinate (A5) at (0,1);            \coordinate (A4) at (-0.707,0.707);   \coordinate (A3) at (-1,0);           \coordinate (A2) at (-0.707,-0.707);  \coordinate (X) at (0.924,0.383);     

\fill[black] (A1) circle (0.02);
    \fill[black] (A2) circle (0.02);
    \fill[black] (A3) circle (0.02);
    \fill[black] (A4) circle (0.02);
    \fill[black] (A5) circle (0.02);
    \fill[red] (X) circle (0.02);
    
\node[below] at (A1) {$A_1$};
    \node[below left] at (A2) {$A_2$};
    \node[left] at (A3) {$A_3$};
    \node[above left] at (A4) {$A_4$};
    \node[above] at (A5) {$A_5$};
    \node[above right] at (X) {$X$};
\end{tikzpicture}
\begin{tikzpicture}[scale=3]
\draw[thick, black] (0,0) circle (1);
    
\coordinate (A1) at (0,-1);           \coordinate (A2) at (-0.707,-0.707);  \coordinate (A3) at (-1,0);           \coordinate (A4) at (-0.707,0.707);   \coordinate (A5) at (0,1);            \coordinate (B4) at (0.707,0.707);    \coordinate (B3) at (1,0);            \coordinate (B2) at (0.707,-0.707);   \coordinate (X) at (0.924,0.383);     

\fill[black] (A1) circle (0.02);
    \fill[black] (A2) circle (0.02);
    \fill[black] (A3) circle (0.02);
    \fill[black] (A4) circle (0.02);
    \fill[black] (A5) circle (0.02);
    \fill[black] (B2) circle (0.02);
    \fill[black] (B3) circle (0.02);
    \fill[black] (B4) circle (0.02);
    \fill[red] (X) circle (0.02);
    
\node[below] at (A1) {$A_1 = B_1$};
    \node[below left] at (A2) {$A_2$};
    \node[left] at (A3) {$A_3$};
    \node[above left] at (A4) {$A_4$};
    \node[above] at (A5) {$A_5=B_5$};
    \node[right] at (B3) {$B_3$};
    \node[above right] at (B4) {$B_4$};
    \node[below right] at (B2) {$B_2$};
    \node[above right] at (X) {$X$};
\end{tikzpicture}
\caption{Simplicial complexes in type $I_2(5)$.  The first image is the pentagon $\Sigma_\tau$ corresponding to the off-the-wall stability condition $\tau$, whose stable objects correspond to the boundary points $A_1,A_2,A_3,A_4,A_5$.
  The spherical object $X$ corresponds to a point on the $A_5A_1$ edge of that pentagon.
  The second image is the octagon $\Sigma_{\tau'}$ associated to the on-the-wall stability condition $\tau'$.
  The object $X$ is now a point on the $B_4B_3$ edge of that octagon.
  The PL homeomorphism \(\Sigma_{\tau'}\rightarrow \Sigma_\tau\) is linear on each edge of \(\Sigma_{\tau'}\); for example, it takes the edge $B_4B_3$ into a portion of the edge $A_5A_1$.}
\end{figure}

\begin{remark}
  The complexes $\Sigma_\tau$ described above are in the case of 2-CY categories of finite rank 2 Coxeter systems $I_2(n)$ for $n<\infty$.
  However, essentially the same analysis applies in the study of the 2-CY category of type $I_2(\infty)$.  See~\cite{bap.deo.lic:20}, where \(\mathbf{P}(\Stab(\mathcal{C}_2(\infty)))\) and its compactification are described.  In particular, the tiling of \(\mathbf{P}(\Stab(\mathcal{C}_2(\infty)))\) is now by $\infty$-gons, and the vertices of each individual $\infty$-gon have an accumulation point---which is not a vertex corresponding to a spherical object---on the boundary.
  Using this description, one can show that the simplicial complex $\Sigma_\tau$ is homeomorphic to $\mathbf{R}^1\cong S^1-\{\operatorname{pt}\}$.
\end{remark}

\section{Determining objects via their Harder--Narasimhan factors}\label{sec:graded-spherical-support}
In this section, we prove that in a 2-CY category, a spherical object is determined by the ordered sequence of its HN factors.
Towards this goal, we prove results about the factors in a filtration of an object that may be of independent interest.
Fix throughout a \(\k\)-linear hom-finite triangulated category \(\mathcal{C}\), that admits a dg-enhancement.

We recall the \(\Hom\) spectral sequence developed in Appendix~\ref{appendixA} for the particular case of a two-step filtration.
Suppose that we have a two-step filtration of an object \(X\), as follows.
\[
  \begin{tikzcd}
    0 \arrow{rr} && X_1 \arrow{rr} \arrow{dl} && X_2 = X \arrow{dl}\\
    & A_1\arrow[dashed]{ul}{+1} && A_2\arrow[dashed]{ul}{+1}[swap]{\alpha}
  \end{tikzcd},
\]
By~\cref{thm:sstriangulated}, we have a spectral sequence with \(E_1\) page
\[ E_1(p,q) = \bigoplus_{\k} \Hom^{p+q}(A_{p+i}, A_{i})\]
that converges to \(\Hom^{p+q}(X,X)\).
Note that \(E_1(p,q)\) must be zero if \(p \notin \{-1,0,1\}\).
Thus the \(E_1\) page is concentrated in three vertical columns indexed by \(p \in \{-1,0,1\}\), with differential mapping horizontally: \(E_1(p,q) \to E_1(p+1,q)\).
Abbreviating \(\Hom^j(-,-)\) as \((-,-)^j\), the \(E_1\) page looks as follows.
\begin{equation}\label{eq:two-step-ss}
  \begin{tikzcd}[column sep=2em,row sep=1em]
    \vdots & \vdots & \vdots & \vdots & \vdots \\    
    0\arrow{r} & (A_1,A_2)^q\arrow{r} & (A_1,A_1)^{q+1}\oplus (A_2,A_2)^{q+1}\arrow{r} & (A_2,A_1)^{q+2}\arrow{r} & 0\\
    0\arrow{r} & (A_1,A_2)^{q-1}\arrow{r} & (A_1,A_1)^q\oplus (A_2,A_2)^q\arrow{r} & (A_2,A_1)^{q+1}\arrow{r} & 0 \\
    0\arrow{r} & (A_1,A_2)^{q-2}\arrow{r} & (A_1,A_1)^{q-1}\oplus (A_2,A_2)^{q-1}\arrow{r} & (A_2,A_1)^q\arrow{r} & 0\\
    \vdots & \vdots & \vdots & \vdots &\vdots\\        
  \end{tikzcd}
\end{equation}
An immediate consequence of the spectral sequence is the ``Mukai lemma''.
The following is a strengthening of the version due to Huybrechts, Macr\'i, and Stellari as~\cite[Lemma 2.7]{huy.mac.ste:08}, but the main idea has been attributed in that paper to Mukai~\cite{muk:87}.
For variants of this lemma, see also~\cite[Lemma 12.2]{bri:08} and~\cite[Corollary 2.3]{huy:12}.
\begin{lemma}\label{lem:mukai-lemma}
  Consider a distinguished triangle
  \[A_1 \to X \to A_2 \xrightarrow{+1}.\]
  Suppose that \(\Hom^0(A_1,A_2) = 0\) and \(\Hom^2(A_2,A_1) = 0\).
  Then we have
  \[\dim \Hom^1(A_1,A_1) + \dim \Hom^1(A_2,A_2) \leq \dim \Hom^1(X,X).\]
\end{lemma}
\begin{proof}
  Regard the triangle as a two-step filtration of \(X\) with factors \(A_1\) and \(A_2\).
  By~\cref{thm:sstriangulated} we have a spectral sequence for \(\Hom^{\bullet}(X,X)\) whose \(E_1\) page is~\eqref{eq:two-step-ss}.
  We show that the incoming and outgoing differentials at \((p,q) = (0,1)\) vanish on the \(E_1\) page as well as all subsequent pages.
  On the \(E_1\) page, the incoming differential at \((0,1)\) is
  \[ \Hom^0(A_1, A_2) \to \Hom^1(A_1, A_1) \oplus \Hom^1(A_2, A_2),\]
  and the outgoing differential is
  \[ \Hom^1(A_1, A_1) \oplus \Hom^1(A_2, A_2) \to \Hom^2(A_2,A_1).\]
  By our hypothesis, we have \(\Hom^0(A_1, A_2) = \Hom^2(A_2,A_1) = 0\), and so both differentials are zero.
  On the subsequent pages, the incoming and outgoing differentials at \((0,1)\) have source or target beyond the range \(-1 \leq p \leq 1\), and hence are automatically zero.

  We conclude that \(\Hom^1(A_1, A_1) \oplus \Hom^1(A_2, A_2)\) survives until the \(E_{\infty}\) page.
  The statement follows.
\end{proof}
The lemma above applies most readily to Harder--Narasimhan filtrations in a 2-CY category.
We deduce the following corollary.
\begin{corollary}\label{cor:hn-blocks-have-no-hom-1}
  Let \(\mathcal{C}\) be a \(\k\)-linear hom-finite triangulated category with a dg enhancement.
  Suppose also that \(\mathcal{C}\) is 2-Calabi--Yau, and fix a stability condition on \(\mathcal{C}\).  
  Let \(X\) be an object such that \(\Hom^1(X,X) = 0\).
  Consider its HN filtration
  \[0 = X_0 \to X_1 \to \cdots \to X_n = X,\]
  with factors \(A_i\).
Then for each \(i\), we have \[\Hom^1(X_i,X_i) = \Hom^1(A_i,A_i) = 0.\]
\end{corollary}
\begin{proof}
  Consider the triangles \(X_{i-1} \to X_i \to A_i \xrightarrow{+1}\).
  We have \(\Hom^0(X_{i-1},A_i) = 0\) for phase reasons.
  Thus by the 2-CY property, we also have \(\Hom^2(A_i,X_{i-1}) = 0\).
  Now we can apply~\cref{lem:mukai-lemma} iteratively to the triangles above, starting at \(i = n\) and going down to \(i = 1\).
\end{proof}

The two-step spectral sequence has applications beyond the Mukai lemma, which we now explore.
For these applications, we need a more refined understanding of the differentials.
The next few lemmas provide this understanding.
\begin{lemma}\label{lem:bimodule-generator}
  Let \(\mathcal{C}\) be a \(\k\)-linear hom-finite triangulated category with a dg enhancement.
  Consider a distinguished triangle
  \[A_1 \to X \to A_2 \xrightarrow{\alpha} A_1[1]\]
  such that  \(\Hom^0(A_1,A_2) = 0\) and 
  \(\Hom^1(X,X) = 0\).
  Then the element \(\alpha \in \Hom^1(A_2,A_1)\) generates \(\Hom^1(A_2,A_1)\) as a bimodule for the left action of \(\Hom^0(A_1,A_1)\) by post-composition and right action of \(\Hom^0(A_2,A_2)\) by pre-composition.
\end{lemma}
\begin{proof}
  Again, regard the given distinguished triangle as a two-step filtration with factors \(A_1\) and \(A_2\), and consider the \(E_1\) page~\eqref{eq:two-step-ss} of the spectral sequence.
  The outgoing differential at \((p,q) = (0,1)\) is
  \begin{align*}
    \delta \colon (A_1,A_1)^0 \oplus (A_2,A_2)^0 &\to (A_2,A_1)^1,\\
    (f,g) &\mapsto f \alpha - \alpha g.
  \end{align*}
  Since \((A_1,A_2)^0 = 0\) by assumption, the cokernel of \(\delta\) survives on the \(E_{\infty}\)-page.
  However, by assumption, \((X,X)^1 = 0\), so the cokernel of \(\delta\) must be zero.
  The image of \(\delta\) is contained in the bimodule generated by \(\alpha\), so the proof is complete.
\end{proof}

We briefly digress to prove a technical fact about endomorphisms of cyclic modules over finite dimensional algebras.  It is used in the proof of our next lemma about the uniqueness of triangles.
\begin{lemma}\label{lem:cyclic-generators-related-by-invertible}
  Let \(A\) be a finite-dimensional algebra over an infinite field \(\k\).
  Let \(V\) be a cyclic left \(A\)-module.
  Let \(r \in A\) such that the linear map \(v \mapsto rv\) from \(V \to V\) is invertible.
  Then there is an invertible \(\widetilde{r} \in A\) such that
    \(\widetilde{r}v = r v\) for every \(v \in V\).
\end{lemma}
\begin{proof}
  Given a left \(A\)-module \(M\) and \(r \in A\), we let \(r_M \colon M \to M\) be the \(\k\)-linear map \(m \mapsto rm\).
  
  Fix a surjection of left \(A\)-modules \(A \to V\) and let \(I\) be the kernel.
  Then we have a filtration of left \(A\)-modules \(0 \subset I \subset A\) with factors \(I\) and \(V\).
  This filtration is preserved by left-multiplication by any element of \(A\).
  
  Let \(r \in A\) be as in the statement.
  Then \(r_V\) is an isomorphism.
  If \(r_I\) is also an isomorphism, then \(r_A\) is an isomorphism, and thus \(r \in A\) is invertible.
  If this is not already the case, we modify \(r\) to \(\widetilde r\) such that \(\widetilde r _V = r_V\) and \(\widetilde r_{I}\) is invertible.

  To construct \(\widetilde r\), we let \(f(x)\) be the minimal polynomial of the linear map \(r_V\).
  We claim that \(\widetilde r = r + \lambda f(r)\) for a sufficiently generic \(\lambda\) satisfies the desired properties.
  Observe that for any \(\lambda\), we have \(\widetilde r_V = r_{V}\).
  So we only need to check that \(\widetilde r_{I}\) is invertible.
    Let \(\{\mu_i\}\) be the eigenvalues of \(r_I\) in an algebraic closure of \(\k\).
  Then the eigenvalues of \(r_I + \lambda f(r)_{I}\) are \(\mu_i + \lambda f(\mu_{i})\).
  Note that \(f(0) = \det(r_V) \neq 0\); so \(\mu_i\) and \(f(\mu_i)\) are never both zero.
  We simply choose \(\lambda\) in the infinite field \(\k\) so that \(\mu_i + \lambda f(\mu_i) \neq 0\) for any \(i\).
\end{proof}

\begin{lemma}\label{lem:left-module-generator}
Let \(\mathcal{C}\) be a \(\k\)-linear hom-finite triangulated category with a dg enhancement.
Let \(\alpha \colon A_2 \to A_1[1]\) and \(\beta \colon A_2 \to A_1[1]\) be morphisms that fit into exact triangles
  \[A_2 \xrightarrow{\alpha} A_1[1] \to X \xrightarrow{+1} \text{ and } A_2 \xrightarrow{\beta} A_1[1] \to X' \xrightarrow{+1}.\]
  Suppose that both \(\alpha\) and \(\beta\) generate \(\Hom^1(A_2,A_1)\) as a left \(\Hom^0(A_1,A_1)\) module.
  Then the two triangles are isomorphic.
\end{lemma}
\begin{proof}
  Since \(\alpha\) is a generator of \(\Hom^1(A_2,A_1)\), and \(\beta \in \Hom^1(A_2,A_1)\), we can find an element \(g \in \Hom^0(A_1,A_1)\) such that \(\beta = g \alpha\).
  By applying~\cref{lem:cyclic-generators-related-by-invertible} to the algebra \(\Hom^0(A_1,A_1)\), we may assume that \(g\) is invertible.
  We then have a commuting square
  \[
    \begin{tikzcd}
      A_2\arrow{r}{\alpha}\arrow{d}{\id} & A_1[1] \arrow{d}{g}\\
      A_2\arrow{r}{\beta} & A_1[1]
    \end{tikzcd},
  \]
  where the vertical maps are isomorphisms.
  It follows that we have an isomorphism \(X \xrightarrow{\cong} X'\) which makes the triangles isomorphic.  
\end{proof}

We now have the tools to prove the main result of this section.
The setup is as follows: \(\k\) is an infinite field; \(\mathcal{C}\) is a \(\k\)-linear 2-CY hom-finite triangulated category with a dg enhancement;  \(\tau\) is a stability condition on \(\mathcal{C}\) in which every \(\tau\)-semistable spherical object is \(\tau\)-stable.
\begin{theorem}\label{thm:object-determined-by-hn-filtration}
  Consider the setup above.
  Let \(X\) be a spherical object with HN filtration
  \[0 \to X_0 \to X_1 \cdots \to X_n = X,\]
  with factors \(A_i\).
  Let \(W\) be another spherical object with HN filtration
  \[0 \to W_0 \to W_1 \cdots \to W_n = W,\]
  which also has factors \(A_i\).
  Then \(X \cong W\).
\end{theorem}
\begin{proof}
  We prove the result by induction on the length of the filtration.
  Let us first tackle the base case of a 2-step filtration.
  Note that \(X_0 = A_0 = W_0\).
So we have distinguished triangles
  \[A_0 \to X_1 \to A_1 \xrightarrow{\alpha} A_0[1]\quad \text{and} \quad A_0 \to W_1 \to A_1 \xrightarrow{\beta} A_0[1].\]
  Both of these triangles satisfy the hypotheses of~\cref{lem:bimodule-generator}, and hence \(\alpha\) and \(\beta\) both generate \(\Hom^1(A_1,A_0)\) as a \(\Hom^0(A_0[1],A_0[1])-\Hom^0(A_1,A_1)\) bimodule.
  However, \(\Hom^0(A_1,A_1) = \k\) because \(A_1\) is spherical.
  Therefore both \(\alpha\) and \(\beta\) generate \(\Hom^1(A_1,A_0)\) as a left \(\Hom^0(A_0[1],A_0[1])\) module.
  We are now in the setting of~\cref{lem:left-module-generator}, and we see that the two distinguished triangles above are isomorphic.
  In particular, there is an isomorphism \(X_1 \xrightarrow{\cong} W_1\) compatible with the filtration.

  For the induction step, consider the triangles
  \[X_{i-1} \to X_i \to A_i \xrightarrow{+1} \quad \text{and}\quad W_{i-1} \to W_i \to A_i \xrightarrow{+1}.\]
  By induction, we may identify \(W_{i-1}\) with \(X_{i-1}\) and obtain distinguished triangles
  \[X_{i-1} \to X_i \to A_i \xrightarrow{\alpha} X_{i-1}[1] \quad \text{and}\quad X_{i-1} \to W_i \to A_i \xrightarrow{\beta} X_{i-1}[1].\]
  
  By~\cref{cor:hn-blocks-have-no-hom-1}, we know that \(\Hom^1(X_{i-1},X_{i-1}) = \Hom^1(X_{i-1}[1], X_{i-1}[1]) = 0\).
  Thus the triangles above satisfy the hypotheses of~\cref{lem:bimodule-generator}.
  Moreover, since \(A_i\) is spherical, we know that \(\Hom^0(A_i,A_i) = \k\), and thus the maps \(\alpha\) and \(\beta\) satisfy the hypotheses of~\cref{lem:left-module-generator}.
  Therefore we obtain a compatible isomorphism \(X_i \xrightarrow{\cong} W_i\).
  The induction step is complete.
  As a result, \(X \cong W\).
\end{proof}

\appendix

\section{A spectral sequence for homomorphisms between filtered objects}\label{appendixA}
Given two filtered objects of a triangulated category, we construct a spectral sequence that begins with the homomorphisms between the factors of the filtration and converges to the homomorphisms between the objects.
We use this spectral sequence in the proof of~\cref{thm:object-determined-by-hn-filtration}.
This construction is surely known to the experts; for instance, it is similar to the one considered in~\cite[Theorem 2.4]{rud:94}.
Nevertheless, we could not find the precise result in the literature.

We begin with a simpler setting.
Let \(\mathcal{A}\) be an abelian category and \(K(\mathcal{A})\) the category of bounded chain complexes on \(\mathcal{A}\) with the morphisms given by chain maps up to homotopy.
Suppose we are given objects \(A_1, \dots, A_n\) of \(K(\mathcal{A})\).
An \emph{iterated cone} of \(A_1, \dots, A_n\) is an object of \(K(\mathcal{A})\) constructed as follows.
If \(n = 1\), then the only such object is \(A_1\).
If \(n = 2\), then such an object is the cone of a morphism \(A_2[-1] \to A_1\).
If \(n = 2\), then such an object is the cone of a morphism \(A_3[-1] \to A_2\), where \(X_2\) is an iterated cone of \(A_1, A_2\).
In general, an iterated cone of \(A_1, \dots, A_n\) is the cone of a morphism \(A_n[-1] \to X_{n-1}\), where \(X_{n-1}\) is an iterated cone of \(A_1, \dots, A_{n-1}\).

Let \(X\) be an iterated cone of \(A_1, \dots, A_n\).
Unravelling the definition, we see that
\[ X = \operatorname{Cone}(A_n[-1] \to \operatorname{Cone}(A_{n-1}[-1] \to \cdots \to \operatorname{Cone}(A_2[-1] \to A_1) \cdots )).\]
Let \(X_i\) be the object obtained at the \(i\)-th step of the above construction.
That is,
\[ X_i = \operatorname{Cone}(A_i[-1] \to \operatorname{Cone}(A_{i-1}[-1] \to \cdots \to \operatorname{Cone}(A_2[-1] \to A_1) \cdots )).\]
Then \(X_{0} = 0\) and \(X_n = X\) and we have an exact triangle
\[ X_{i-1} \to X_i \to A_i \xrightarrow{+1}.\]
Putting the maps \(X_{i-1} \to X_i\) together gives an ascending filtration
\[ 0 = X_0 \to X_1 \to X_2 \to \cdots \to X_n = X \]
whose associated graded is \(\operatorname{gr}_i(X_{\bullet}) = A_i\).

Let us describe the complex \(X\) more explicitly.
For \(i = 1, \dots, n\) and an integer \(a\), let \(A_{i}^{a}\) denote the \(a\)-th chain group of \(A_i\), so that \(A_i\) is the complex
\[ \cdots \to A_{i}^{a} \to A_{i}^{a+1} \to \cdots .\]
Then the \(a\)-th chain group \(X^{a}\) of \(X\) is the direct sum
\begin{equation}\label{eqn:directsum}
  X^{a} = \bigoplus_{i = 1}^{n} A_{i}^{a}.
\end{equation}
With this direct sum decomposition, the differential \(X^{a} \to X^{a+1}\) is lower triangular---it maps \(A_{i}^{a}\) to \(\bigoplus_{i' \leq i} A_{i'}^{a+1}\).
The diagonal component \(A_{i}^{a} \to A_{i}^{a+1}\) is (up to sign) the differential of \(A_i\);
the other components depend on the morphisms used in taking the cones.

\begin{proposition}\label{thm:sscomplexes}
  Let \(A_1, \dots, A_n\) and \(B_1, \dots, B_m\) be objects of \(K(\mathcal{A})\).
  Suppose \(X\) is an iterated cone of \(A_1, \dots, A_n\) and \(Y\) is an iterated cone of \(B_1, \dots, B_m\).
  Then we have a spectral sequence whose \(E_1\) page is
  \[ E_1(p,q) = \bigoplus_{k}\Hom^{p+q}(A_{k+p},B_{r})\]
  that converges to \(\Hom^{p+q}(X,Y)\).
\end{proposition}
By our convention, the differential on the \(E_r\) page goes from \((p,q)\) to \((p+r,q-r+1)\).
\begin{proof}
  Consider the complex \(C = \Hom(X,Y)\).
  Then \(\Hom^n(X,Y)\) is the \(n\)-th cohomology \(H^n(C)\).

  Using the direct sum decomposition of the chain groups of \(X\) in terms of the chain groups of \(A_i\) from \eqref{eqn:directsum} and a similar description for \(Y\), we get
  \begin{align*} 
    C^{k} &= \bigoplus_{b-a = k} \Hom(X^a,Y^b) \\
    &= \bigoplus_{b-a = k} \bigoplus_{i,j} \Hom(A_{i}^{a}, B_{j}^{b}).
  \end{align*}
  The differential sends \(\Hom(A_{i}^{a}, B_{j}^{b})\) to the direct sum
  \[\bigoplus_{i' \geq i} \Hom(A_{i'}^{a-1}, B_{j}^{b}) \oplus \bigoplus_{j' \leq j} \Hom(A_{i}^{a}, B_{j'}^{b+1}).\]
  For \(\ell \in \mathbf{Z}\), set
  \[ C_{\ell}^{k} = \bigoplus_{b-a = k}\bigoplus_{i-j \geq \ell} \Hom(A_{i}^{a}, B_{j}^{b}).\]
  Then \(C_{\ell}\) is a subcomplex of \(C\).
  In fact, we have a decreasing filtration
  \[ C = C_{-m} \leftarrow C_{-m+1} \leftarrow \cdots \leftarrow C_{n} \leftarrow C_{n+1} = 0.\]
  Its associated graded is the complex whose \(k\)-th chain group is
  \[ \operatorname{gr}_{\ell}(C_{\bullet})^{k} = \bigoplus_{b-a = k}\bigoplus_{i-j = \ell} \Hom(A_{i}^{a}, B_{j}^{b}).\]
  The differential in \(\operatorname{gr}_{\ell}(C_{\bullet})\) maps \(\Hom(A_{i}^{a}, B_{j}^{b})\) to the direct sum
  \[ \Hom(A_{i}^{a-1}, B_{j}^{b}) \oplus \Hom(A_{i}^{a}, B_{j}^{b+1}).\]
  Its first component is obtained by pre-composing with the differential \(A_{i}^{a-1} \to A_{i}^{a}\) of \(A_{i}\) and its second component is obtained by post-composing with the differential \(B_{j}^{b} \to B_{j}^{b+1}\) of \(B_{j}\).
  Thus, we have a direct sum decomposition of complexes
  \[ \operatorname{gr}_{\ell}(C_{\bullet}) = \bigoplus_{i-j = \ell} \Hom(A_i,B_j).\]

  The spectral sequence asserted in the statement is the spectral sequence associated to the filtered complex \(C_{\bullet}\); see \cite[\href{https://stacks.math.columbia.edu/tag/012K}{Tag 012K}]{the:22}.
  Its \(E_1\) page is given by
  \[ E_1(p,q) = H^{p+q}(\operatorname{gr}_{p}(C_{\bullet})) = \bigoplus_{k}\Hom^{p+q}(A_{k+p},B_{k}),\]
  and it converges to \(H^{p+q}(C^{\bullet}) = \Hom^{p+q}(X,Y)\).
\end{proof}

Let us now state the analogue of~\cref{thm:sscomplexes} in the triangulated setting under suitable assumptions.
Let \(\mathcal{C}\) be a triangulated category.
Let \(X\) and \(Y\) be objects of \(\mathcal{C}\) with filtrations
\[ 0 = X_0 \to X_1 \to \cdots \to X_n = X\]
and
\[ 0 = Y_0 \to Y_1 \to \cdots \to Y_m = Y.\]
Let \(A_i\) and \(B_j\) be the objects that complete the triangles \(X_{i-1} \to X_i \to A_i \xrightarrow{+1}\) and \(Y_{j-1} \to Y_j \to B_j \xrightarrow{+1}\).
In our applications, the categories we study all have dg enhancements (see, e.g.~\cite{kel:06}).
In these cases we have the following analogue of~\cref{thm:sscomplexes}.
\begin{proposition}\label{thm:sstriangulated}
  In the setup above, if \(\mathcal{C}\) has a dg enhancement, then 
we have a spectral sequence with \(E_1\) page
  \[ E_1(p,q) = \bigoplus_{k} \Hom^{p+q}(A_{k+p},B_{k})\]
  that converges to \(\Hom^{p+q}(X,Y)\).
\end{proposition}
\begin{proof}
  By definition, if \(\mathcal{C}\) has a dg enhancement, then it is isomorphic as a triangulated category to the homotopy category of a dg category \(\mathcal{D}\).
  The triangulated structure on the homotopy category of \(\mathcal{D}\) consists of triangles that are isomorphic to mapping cones (see, e.g.~\cite[Lemma 3.3]{kel:06}).
  Thus we can reduce to the situation of iterated cones treated in~\cref{thm:sscomplexes}.
\end{proof}

\bibliographystyle{siam}

\end{document}